\documentclass[12pt]{article}
\usepackage[english]{babel}
\usepackage{graphicx}
\usepackage{amsmath,amsthm,amssymb}
\usepackage{times}
\usepackage{enumerate}
\usepackage{enumitem}
\usepackage{amsmath}
\usepackage{amsfonts}
\usepackage{pdfsync}
\usepackage{cleveref}
\usepackage{autonum}
\usepackage[usenames, dvipsnames]{color}

\newcommand{\R}{{\mathbb R}}

\renewcommand{\L}{\Lambda}

\DeclareMathOperator{\dom}{\text{Dom}\,}

\DeclareMathOperator{\Lip}{Lip}
\DeclareMathOperator{\dist}{dist}

\DeclareMathOperator*{\essinf}{ess\,inf}
\DeclareMathOperator{\Dom}{Dom}

\def\titlerunning#1{\gdef\titrun{#1}}
\makeatletter
\def\author#1{\gdef\autrun{\def\and{\unskip, }#1}\gdef\@author{#1}}
\def\address#1{{\def\and{\\}\renewcommand{\thefootnote}{}%
\footnote {#1}}%
\markboth{\autrun}{\titrun}}
\makeatother
\def\email#1{e-mail: #1}
\def\subjclass#1{{\renewcommand{\thefootnote}{}%
\footnote{\emph{Mathematics Subject Classification (2010):} #1}}}
\def\keywords#1{\par\medskip
\noindent\textbf{Keywords.} #1}
\usepackage{amsthm}

\newtheorem*{AssumptionG}{Superlinearity}
\newtheorem*{AssumptionA}{Structure Assumption (A)}
\newtheorem*{ConditionH}{Growth Condition (H$_B^{\delta}$)}
\newtheorem*{ConditionHNEW}{Growth Condition (M$_B^{\delta}$)}

\newtheorem*{basicass}{Basic Assumptions and Notation}
\newtheorem*{conditionS}{Condition (S)}

\newtheorem{theorem}{Theorem}[section]
\newtheorem{lemma}[theorem]{Lemma}

\theoremstyle{definition}
\newtheorem{definition}[theorem]{Definition}
\newtheorem{example}[theorem]{Example}
\newtheorem{corollary}[theorem]{Corollary}
\newtheorem{proposition}[theorem]{Proposition}

\newtheorem*{Notation}{Notation}
\newtheorem*{GrowthGloc}{Growth Condition (G)}

\theoremstyle{remark}
\newtheorem{remark}[theorem]{Remark}

\numberwithin{equation}{section}

\begin{document}
\titlerunning{Equi-Lipschitz minimizing trajectories}
\title{Equi-Lipschitz minimizing trajectories  for non coercive, discontinuous, non convex Bolza controlled-linear optimal control problems}
\author{
Carlo Mariconda
}
\maketitle
\date
\address{
Carlo Mariconda, ORCID $0000-0002-8215-9394$),  Universit\`a
degli Studi di Padova,  Dipartimento di Matematica  ``Tullio Levi-Civita'',   Via Trieste 63, 35121 Padova, Italy;  \email{carlo.mariconda@unipd.it}}
\begin{abstract}This article deals with the Lipschitz regularity of the ''approximate`` minimizers for the  Bolza type control functional of the form
\[J_t(y,u):=\int_t^T\L(s,y(s), u(s))\,ds+g(y(T))\]
among the pairs $(y,u)$ satisfying a prescribed initial condition $y(t)=x$, where the state $y$ is absolutely continuous, the control $u$ is summable and the dynamic is controlled-linear of the form $y'=b(y)u$. For $b\equiv 1$ the above becomes a problem of the calculus of variations. The Lagrangian $\L(s,y,u)$ is assumed to be either convex in the variable $u$ on every half-line from the origin (radial convexity in $u$), or partial differentiable in the control variable and satisfies a local Lipschitz regularity  on the time variable, named Condition (S). It is allowed to be extended valued,  discontinuous in $y$ or in $u$, and non convex in $u$.\\
We assume a very mild growth condition,  actually a violation of the Du Bois-Reymond -- Erdmann equation for high values of the control, that is fulfilled  if the Lagrangian is coercive  as well as in some almost linear cases.  The main result states that, given any admissible pair $(y,u)$, there exists a more convenient  admissible pair $(\overline y, \overline u)$ for $J_t$ where $\overline u$ is  bounded, $\overline y$ is Lipschitz, with   bounds and ranks that are uniform with respect to $t,x$ in the compact subsets of $[0,T[\times\R^n$. The result is new even in the superlinear case. As a consequence, there are minimizing sequences that are formed by pairs of equi-Lipschitz trajectories and equi--bounded controls.\\ A new existence and regularity result follows without assuming any kind of Lipschitzianity in the state variable.\\
We deduce, without any need of growth conditions,  the nonoccurrence of the Lavrentiev phenomenon for a wide class of Lagrangians containing those that satisfy Condition (S), are bounded on bounded sets ``well'' inside the effective domain and are  radially convex  in the control variable.\\
The methods are based on a reparametrization technique and  do not involve the Maximum Principle.
\end{abstract}
\thanks{
I thank Julien Bernis and Piernicola Bettiol for the useful discussions we had in Brest, Marc Quincampoix and the Laboratoire de Math\'ematiques de Bretagne Atlantique for the warm hospitality during the initial part of the work. I am further indebted with Piernicola Bettiol for his useful comments on the structure of the manuscript and with Francis Clarke for his  encouragement.\\
I am grateful to the  referees who read very carefully the work, and for their suggestions that led  me to an improved version of the  paper.\\
This research is partially supported by the  Padua University grant SID 2018 ``Controllability, stabilizability and infimum gaps for control systems'', prot. BIRD 187147}


\setcounter{tocdepth}{1}
\tableofcontents
\tableofcontents
\section{Introduction}
The main object of the article concerns the existence of ``nice'' pairs of approximate solutions to an optimal control problem.  For the sake of clarity, we motivate the  core of the paper by means of the basic problem of the calculus of variations.
\\
The classical problem of the calculus of variations consists of minimizing an integral functional
\[\begin{aligned}\min I(y):=\int_t^T\L(s,y(s), y'(s))\,ds,\,\, &y\in W^{1,1}([t,T];\R^n),\\ y(t)=x\in\R^n,& \, y(T)=y_T \in\R^n,\end{aligned}\]
where $\L$ is a positive, Lebesgue-Borel Lagrangian.
The main ingredients in order to obtain the \textbf{existence} of a solution are summarized in Tonelli's theorem:
\begin{itemize}
\item Lower semicontinuity of $\L(s,y,u)$ with respect to $(y,u)$;
\item Convexity of $\L(s,y,u)$ with respect to $u$;
\item Superlinearity of $\L(s,y,u)$ with respect to $u$:
\[\forall (s,y,u)\in[t,T]\times\R^n\times\R^n\qquad \L(s,y,u)\ge \Phi(|u|),\] where
$\Phi:[0, +\infty[\to [0, +\infty[$ and $\displaystyle\lim_{r\to +\infty}\dfrac{\Phi(r)}{r}=+\infty$.
\end{itemize}
When a  minimizer of $I$ exists, a first step towards regularity is  looking at its \textbf{Lipschitzianity}.\\
When $\L$ is {\em autonomous}, superlinearity alone suffices  to ensure the Lipschitz continuity of the minimizers (see \cite{AAB, CVTrans,  DMF}).
Weaker growth conditions were considered in the last decades, requiring a specific behavior of the Hamiltonian $H(s,y,u,p):=p\cdot u-\L(s,y,u)$ associated with $\L(s,y,u)$, as $p$ belongs to the convex
subdifferential of $u\mapsto\L(s,y,u)$ and $|u|\to +\infty$.
The essential idea of using such indirect growth conditions for the purposes of existence and regularity is due to F. Clarke, who introduced Condition (H) in his seminal paper \cite{Clarke1993} of 1993 for Lagrangians possibly nonautonomous, extended valued, with state and velocity constraints.
Few years later, with different methods, A. Cellina and his school began working around a growth condition (G), formulated first in \cite{CTZ} for  continuous Lagrangians of the form $\L(s,y,u)=f(y)+g(u)$ with $g\in C^1(\R)$ and, in 2003  \cite{Cellina, CF1} for autonomous and continuous Lagrangians.\\
The growth conditions (H) and (G)  will be thoroughly examined below. At this stage we just mention here that if $\L$ is  bounded on bounded sets then superlinearity implies Condition (G) and, in the {\em real valued case},  the validity of Condition (G) implies that of Condition (H); Conditions (G) and (H) are satisfied by some  Lagrangians with almost linear growth, e.g., $\L(u)=|u|-\sqrt{|u|}$ satisfies (G) as well as (H), and  some Lagrangians of the form
$\L(s,y,u)=h(s,y)\sqrt{1+|u|^2}$ -- in particular  $\L(u)=\sqrt{1+|u|^2}$ -- satisfy (H), but not (G).\\
In the real valued autonomous case these {\em weak} growth conditions alone (with no need of convexity or continuity assumptions), instead of superlinearity, ensure the Lipschitz regularity of minimizers, as shown by P. Bettiol and C. Mariconda in \cite{BM2}.
\\
In the {\em nonautonomous} case, there are examples of Lagrangians that satisfy Tonelli's assumptions but whose minimizers are not Lipschitz.
Several  regularity results appeared on the subject (see \cite{ClarkeMemoirs,  CVTrans,  PV}), each requiring some extra  assumptions on the state or velocity variable, e.g., local Lipschitz conditions on the state variable or Tonelli -- Morrey type conditions, mostly motivated by the use of the Weierstrass inequality or Clarke's Maximum Principle. A Lipschitz regularity
result without additional  smoothness or convexity requirements on the state and velocity variables was obtained by P. Bettiol -- C. Mariconda in \cite{BM2, BM1} under the growth condition (H). The price to pay, with respect to Tonelli's assumptions,
is the additional local Lipschitz condition (S) on the time variable, thoroughly examined in \S~\ref{sec:AS},
 requiring that $s\mapsto\L(s,y,u)$ is locally Lipschitz and that, for all $(y,u)$,
\begin{equation}\label{tag:Sintro}|D_s\L(s,y,u)|\le \kappa \L(s,y,u)+A|u|+\gamma(s)\quad \text{a.e. }s\in [t,T],\end{equation}
for some $\kappa, A\ge 0$ and $\gamma\in L^1([t,T])$. Moreover, it turns out without need of any growth condition,  that $\L$ is somewhat \textbf{radially convex} on the velocity variable along any given minimizer $y_*$, in the sense that, for a.e. $s\in [t,T]$, the map
\[0<r\mapsto\L(s,y_*(s), ry_*'(s))\]
has a nonempty subdifferential at $r=1$, in the sense on convex analysis. The role of radial convexity in Lipschitz regularity was prefigured in \cite{Clarke1993} by the fact that the velocity constraint is a cone, and was first explicitly formulated for  {\em autonomous Lagrangians} by C. Mariconda and G. Treu in \cite{MTLip}.
\\
The celebrated example by J. M. Ball and V. J. Mizel in \cite{BM} of a nonautonomous polynomial Lagrangian that is superlinear and convex in the velocity variable shows that, the violation of  Condition (S) may lead not only to minimizers that are not Lipschitz, but even to the \textbf{Lavrentiev phenomenon}, i.e., the fact that
\[\inf_{\substack{y\in W^{1,1}([t,T];\R^n)\\y(t)=x,\, y(T)=y_T }}I(y)<
\inf_{\substack{y\in\Lip([t,T];\R^n)\\y(t)=x,\, y(T)=y_T  }}I(y).\]
Condition (S) is not new and appeared in several results. It is sufficient for the validity of the Du Bois-Reymond -- Erdmann equation (see \cite{Cesari} for smooth Lagrangians, and  \cite{BM2} for a discussion in the general case) and, if one replaces in Tonelli's assumptions the superlinearity condition with the slower growth  (H), it plays an essential role in establishing the Lipschitz continuity of minimizers in F. Clarke -- R. Vinter's \cite[Corollary 3]{CVTrans}. Furthermore, F. Clarke proved in \cite{Clarke1993} that it provides existence of a minimizer, which is actually Lipschitz.\\
The Lavrentiev phenomenon has been widely reconsidered in the  1980s, a long  time after M. Lavrentiev and B. Manià realized (see \cite{Man}) that such a pathology could occur.
Here again, the {\em autonomous} case stands on its own:  G. Alberti -- F. Serra Cassano proved in \cite{ASC} that the Lavrentiev phenomenon never occurs if $\L(y,u)$ is just Borel, possibly extended valued; we refer to \cite{BMT}, \cite{BB}, \cite{GBH} for more insights on Lavrentiev's gap.
More precisely if $y(\cdot)$ is an admissible trajectory and $s\mapsto \L(y(s),y'(s))\in L^1([t,T])$, there  is no {\em Lavrentiev gap} at $y$, i.e., there is a sequence $(y_j)_j$ of Lipschitz functions that share the same boundary values with $y$, converging to $y$ in $W^{1,1}([t,T])$ and in energy, i.e. $\displaystyle\lim_{j\to +\infty}I(y_j)=I(y)$.\\
In the {\em nonautonomous} case some additional conditions have to be added. To the author's knowledge, the  criteria for the avoidance of the Lavrentiev phenomenon either follow trivially from the fact that minimizers exist and are Lipschitz or, as in  \cite{Loewen, TZ, Zas}, they require that $\L$ is locally Lipschitz or H\"older continuous in the state variable.   One of the reasons is that, as was pointed out  by D. Carlson in \cite{Carl}, many of the published results can   actually be obtained as a consequence of   Property (D) introduced by L. Cesari and T.S. Angell in \cite{CeAng}.

Regularity conditions on the state variable or convexity in the velocity variable are not satisfied in several problems arising from real life; \textbf{discontinuous Lagrangians} appear for instance  in models arising  from combustion  in non homogeneous media or light propagation in the presence of layers. Some natural questions arise, and are addressed in the paper:
\begin{itemize}
\item[1.] When existence of a minimizer fails because of the lack of continuity of the Lagrangian with respect to the state or velocity variable, can one at least approach the infimum of the functional $I$ through the values of $I$ along \textbf{``nice'' minimizing sequences} (say equi-Lipschitz)?
\item[2.] May Condition (S) on the time variable replace the customary regularity  assumptions in the state variable, in order to prevent the \textbf{Lavrentiev phenomenon}?
\end{itemize}
Problem 1 was considered by A. Cellina and A. Ferriero in \cite{CF1}, where the authors consider autonomous, continuous Lagrangians that are  convex or differentiable in the velocity variable and satisfy the growth condition (G).
While, in the {\em real valued, convex} case, this result may be seen as a consequence of \cite[Theorem 2]{Clarke1993} ((G) implies (H) so that minimizers exist and are Lipschitz), new results arise in the differentiable case or in the {\em extended valued} framework, when Condition (G) and the Condition (H), in its {\em original} formulation (as in  \cite{BM2, BM3, Clarke1993}),  may even not overlap.
Most of the present work is based on the intuition that some steps of the proof of the main result in \cite{CF1} for the basic problem of the calculus of variations could actually be carried on in a more general setting, namely under a weaker growth condition of type (H) instead of (G),  no more continuity assumptions on the state and velocity, nor convexity in the velocity variable and in the slightly wider framework of optimal control problems with a controlled-linear dynamics.

In this article we consider the more general {\em Bolza optimal control} {\rm (P}$_{t,x}${\rm )} of  minimizing  an integral functional
\[J_t(y,u):=\int_t^T\L(s,y(s), u(s))\,ds+g(y(T))\]
among the absolutely continuous arcs $y:[t,T]\to\R^n$ that have a prescribed value at $t$
\[y(t)=x\in\R^n,\]
that are  subject to a {\em state constraint}
\[y(s)\in \mathcal S\subseteq\R^n\quad\forall s\in [t,T],\]
and to a control-linear differential equation
\begin{equation}\label{tag:controleq}y'=f(y,u)=b(y)u,\end{equation}
with
\[u(s)\in\mathcal U\subseteq\R^m\]
and $g$ is positive, possibly extended valued,  $\mathcal U$ is a cone. If $b$ is the identity matrix, $g$ is the indicator function of a point and $\mathcal S=\R^n$ problem (P$_{t,x}$) is the basic problem of the calculus of variations.
The same Bolza problem was considered in \cite{BM3}; the particular form of the dynamics is motivated by the reparametrization techniques used to obtain the results. The results thus apply for instance to the class of problems called of Grushin type (see \cite{MMMT}) and to control problems related to subriemaniann metrics (see \cite{AMR1}).
 We take nonautonomous Lagrangians which are Lebesgue-Borel measurable and possibly extended valued.\\
We assume that the Lagrangian is  measurable, has at least a \textbf{linear growth from below} and satisfies \textbf{Condition (S)}.
We admit two different types of Lagrangians: those that are  \textbf{radially convex w.r.t. the control variable} {\em or} those that are \textbf{partial differentiable w.r.t the control variable};  no kind of lower semicontinuity nor global convexity in the state or control variable are required.\\
The \textbf{extended valued} case needs some extra assumptions. In this situation we impose, moreover,  that $\L$ tends uniformly to $+\infty$ at the boundary of its  effective domain $\Dom(\L)$, together with some structure conditions on $\Dom(\L)$ that are satisfied if, for instance, $\L(s,y,u)=a(s,y)L(u)$ where $a$ is real valued and $\Dom(L)$ is star-shaped.

\begin{figure}[htbp]
\begin{center}
\includegraphics[width=0.45\textwidth]{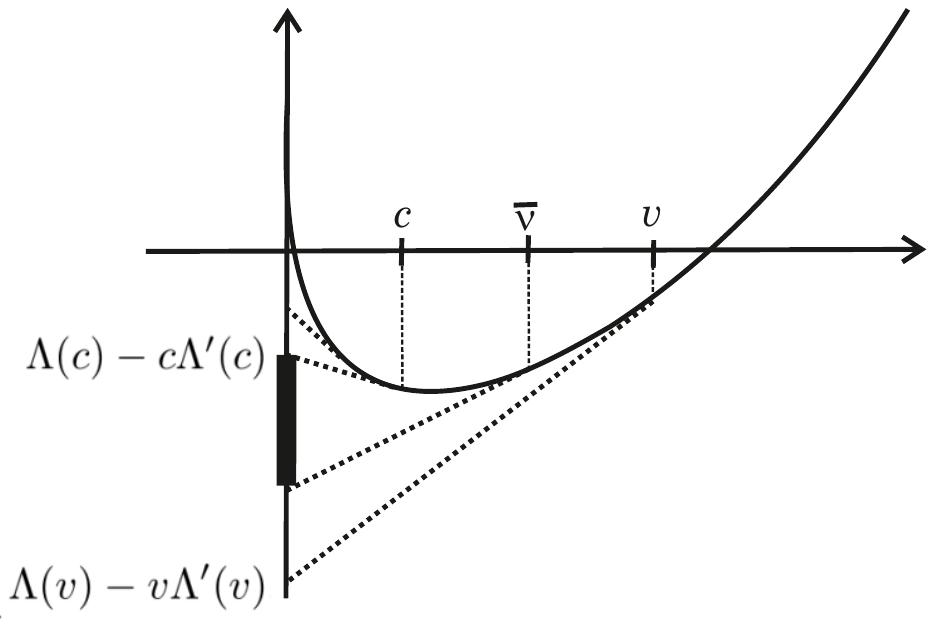}\
\end{center}
\caption{Condition {\rm (H)} in the case of a convex function $\L(u)$ with $\Dom(\L)=[0, +\infty[$: $\displaystyle\inf_{0\le v<c}\L(v)-v \L'(v)<\sup_{v\ge\overline\nu}\L(v)-v \L'(v)$.}
\label{fig:Hintro}
\end{figure}

In Section~\ref{sect:growths} we study various ``slow'' growth conditions  and describe how they are related.
When $\L$ is smooth
Condition (G) imposes that
\[\lim_{|u|\to +\infty}\L(s,y,u)-u\cdot \nabla_u\L(s,y,u)=-\infty\quad \text{uniformly in } s,y.\]
The interpretation of (G) can be easily understood noticing that $ \L(s,y,u)-u\cdot \nabla_u\L(s,y,u)$ is  the value of the intersection with the $w$ axis of the tangent hyperplane to $v\mapsto \L(s,y,v)$ at $v=u$.
Condition (G) has been considered in the framework of the {\em autonomous} case in \cite{Cellina, CF1, CTZ} and extended to the {\em nonautonomous} case in \cite{BM2}.\\
In the smooth setting the original Condition (H), as formulated in \cite{Clarke1993} for the calculus of variations and in \cite{BM3,  FrankTrans} for the optimal control problem considered here, requires that  once $(y(\cdot),u(\cdot))$ is an admissible pair for (P$_{t,x}$) with
\begin{equation}\label{tag:essinf22}\essinf_{s\in [t,T]}|u(s)|<c\end{equation}
then there is $\overline\nu>0$ such that
\begin{align}\label{tag:Hintro22}\sup_{\substack{|v|\ge \overline\nu, v\in \mathcal U\\z\in\R^n}}
\{\L(s,z,v)-v\cdot \nabla_v\L(s,z,v)\}+\Phi(y(\cdot),u(\cdot))<\phantom{AAAAA}
\\ \phantom{AAAAAA}<\inf_{\substack{|v|< {c}, v\in \mathcal U\\ z\in\R^n}}
\left\{\L(s,z,v)-v\cdot \nabla_v\L(s,z,v)\right\}
,\end{align}
where
\[\Phi(y(\cdot),u(\cdot)) :=\int_t^T\left\{\kappa\L(s, y(s), u(s))+A|u(s)|+\gamma(s)\right\}\,ds\]
and the  $\kappa, A, \gamma$ are as in \eqref{tag:Sintro}.
At a first glance, Condition (H) may appear quite involved since it relies on the {\em essinf} of a given admissible pair $(y,u)$ and on $\Phi(y(\cdot),u(\cdot))$, a function depending on  Condition (S). In the autonomous case it appears to be more ductile since in that case $\Phi\equiv 0$ (see Figure~\ref{fig:Hintro} for the interpretation of Condition (H) in the simple case of a Lagrangian of a positive real control variable). However,
as  anticipated in \cite[Theorem 3]{Clarke1993} and proved in \cite{BM3}, conditions \eqref{tag:essinf22} -- \eqref{tag:Hintro22}  represent  merely a  violation of the Du Bois-Reymond -- Erdmann equation for high values of the velocity/control.
With respect to \cite{BM2, BM3, FrankTrans, Clarke1993}  we formulate here Condition (H)  in a slight different way for several reasons. We take into account that  the initial time $t$ and value $x$ may vary. Furthermore, in the {\em extended valued case} the formulation given in \S~\ref{subsect:H} widens the class of functions that satisfy \eqref{tag:Hintro22} (see Remark~\ref{rem:newHwhy} and Example~\ref{ex:Hnew}); as  a byproduct the validity of (G) implies now that of (H) in any ``reasonable'' case (Proposition~\ref{prop:GimpliesH}). At the same time, at least in the {\em real valued and nonautonomous} case, our Condition (H) is slightly more restrictive with respect to the original one due to the presence, in \eqref{tag:Hintro22}, of a technical factor 2 in front of $\Phi(y(\cdot),u(\cdot))$; this does not seem, however,  to have any consequence in concrete applications.\\
A new growth condition (M) is introduced in \S~\ref{sect:M}, so weak that it is in fact satisfied by any Lagrangian that is bounded on the bounded sets ``well-inside'' the effective domain (in the sense of Definition~\ref{def:wellinside}) and radially convex in the control variable. In the real valued, smooth case it simply requires that, for a suitable $c>0$ and $\overline\nu>0$,
 \begin{equation}\label{tag:Mintro}\sup_{\substack{|v|\ge \overline\nu, v\in \mathcal U\\z\in\R^n}}
\{\L(s,z,v)-v\cdot \nabla_v\L(s,z,v)\}<+\infty,\end{equation}
\[
-\infty<\inf_{\substack{|v|< {c}, v\in \mathcal U\\ z\in\R^n}}
\left\{\L(s,z,v)-v\cdot \nabla_v\L(s,z,v)\right\}.\]
\noindent
The main result,  formulated in Theorem~\ref{thm:main2}, considers the two different types of growth (H) or (M):
\begin{itemize}[leftmargin=*]
\item
If Condition (H) is verified, it states that, whenever $(y,u)$ is admissible for (P$_{t,x}$) then there is an admissible pair $(\overline y, \overline u)$ where $\overline y$ is Lipschitz, $\overline u$ is bounded, such that \[J_t(\overline y, \overline u)\le J_t(y,u).\] Moreover, the Lipschitz constant of $\overline y$ and  $\overline u$ are \textbf{uniformly bounded} as $t,x$ vary in compact sets.
\item
If the less restrictive Condition (M) holds, given $\eta>0$ we still get a pair $(\overline y, \overline u)$ with the above regularity properties, and satisfying \[J_t(\overline y, \overline u)\le J_t(y,u)+\eta.\]
Several examples are provided in \S~\ref{sect:examples} to illustrate the growth conditions involved in the article and the applicability of the results.
\end{itemize}
In the proof of Theorem~\ref{thm:main2} the Maximum Principle cannot be invoked,   due to the lack of Lipschitz continuity of the Lagrangian in the state variable. Instead, we extend the method of \cite{CF1} to this more general framework
 in order to build the desired Lipschitz function $\overline y$ via a Lipschitz reparametrization of $y$.
Without entering into the several technical points of the proof, it may be of interest to briefly illustrate the link between reparametrizations and growth conditions. For simplicity, consider the  case of the calculus of variations.
Let $\varphi$ be a smooth, increasing  change of variable on $[t,T]$, $y$ be an admissible trajectory for {\rm (P}$_{t,x}${\rm )}, and set $\overline y(s):=y(\varphi^{-1}(s))$.
Notice that, by taking high values of $\varphi'(\tau)$, one lowers the norm of the derivative of $\overline y(\varphi(\tau))$. The change of variable $s=\varphi(\tau)$ yields
\begin{equation}\label{tag:changevar}I(\overline y)=\int_t^T\L(s, \overline y(s), \overline y'(s))\,ds=
\int_t^T\L\Big(\varphi(\tau), y(\tau), \dfrac{y'(\tau)}{\varphi'(\tau)}\Big)\,\varphi'(\tau)\,d\tau.\end{equation}
Supposing that  $\L$ smooth, the derivative of
$\mu\mapsto \L\Big(\varphi, y,\dfrac{u}{\mu}\Big)\,\mu$ at $\mu=1$ is
\[\L(\varphi, y,u)-u\cdot \nabla_u\L(\varphi, y,u).\]
 The proof consists on finding a suitable increasing and one-to-one change of variable $\varphi:[t;T]\to [t,T]$. By choosing $\overline\nu, c>0$ as in \eqref{tag:Hintro22} (resp.
\eqref{tag:Mintro}),    Conditions of type (H) (resp. (M)) allow to compensate the values of integral in $I$ on the sets  where $|y'|>\overline\nu$ with the ones where $|y'|<c$, up to obtain a lower value than $I(y)$ (resp. $I(y)+\eta$).
 The essential ideas of the multiple step  proof of Theorem~\ref{thm:main2} are described at the beginning of Section~\ref{sect:proofs} for the convenience of the reader.
Many  technical issues  are actually related to the fact that the Lagrangian is allowed to take the value $+\infty$; we invite the reader focused in the real valued case to consult the simplified version in the announcement of the results given in  \cite{BBM3}.
It is worth mentioning that, in  the proof of Theorem~\ref{thm:main2},  the two  growth conditions (H) and (M)  share most of the arguments; their difference play a role just in  few of the many steps. This fact seems to be a byproduct of the care needed to deal with Condition (H) and was unnoticed in \cite{CF1}, where the authors consider the more restrictive (but easier to handle) growth of type (G).

 Theorem~\ref{thm:main2} has several consequences.\\
 Under Condition (H), Corollary~\ref{coro:coro1NEW} 
 yields \textbf{``nice'' minimizing sequences} for (P$_{t,x}$) formed by  \textbf{equi-Lipschitz trajectories and equi--bounded controls} as $t,x$ vary in compact sets. This property, that does not need existence of minimizers, is expected to have a strong impact on the study of the regularity function of the value function
\[V(t,x)=\inf {\rm (P}_{t,x}{\rm )}\]
and will be investigated (\cite{BBM2}, in preparation).\\
As a further consequences of the main result,
the \textbf{existence of a solution} under slow growth conditions to the optimal control problem (P$_{t,x}$) when, in addition to the conditions of Theorem~\ref{thm:main2}, one imposes some standard lower semicontinuity of $\L(s,y,u)$  in $(y,u)$, convexity with respect to $u$, closure of the state constraint set $\mathcal S$, closure and convexity of the control set. Corollary~\ref{coro:coro2} almost overlaps  \cite[Theorem 3]{Clarke1993} when the problem concerns  the calculus of variations, where the major difference relies on the version of Condition (H) mentioned above, but seems to be new in the framework of optimal control problems.
Existence for more general controlled differential equations than \eqref{tag:controleq} was considered in the autonomous case in \cite{FrankTrans}. However, though the  {\em controlled-linear}
structure  of the  system \eqref{tag:controleq} might appear  restrictive,
 the novelty with respect to the known literature is  represented here by the absence of any kind of local Lipschitz condition on the state variable, and by the fact that  $\L$ may  be extended valued.\\
Theorem~\ref{thm:main2} does also provide some answers related to Problem 2. The \textbf{Lavrentiev phenomenon} is excluded in Corollary~\ref{coro:coro2NEW} for a wide class of Lagrangians, assuming a growth condition of type (M). In particular it is avoided (Corollary~\ref{coro:realLav}) when $\L$ is real valued and, moreover:
\begin{itemize}
\item[a)] $\L$ satisfies Condition (S);
\item[b)] $\L$ is radially convex in the control variable;
\item[c)] $\L$ is bounded on bounded sets.
\end{itemize}
We stress again the fact that, differently from other results on the Lavrentiev phenomenon for optimal control problems (see \cite{FrankTrans, Tonon, FMotta}), we do not assume any kind of Lipschitz continuity of $\L$ in the state variable, nor we make use of the Maximum Principle.

\section{Basic setting and notation}\label{sect:notation}
Let $0\le t< T$ and $x\in\R^n$. We consider the Bolza type \textbf{optimal control problem}
\begin{equation}\min J_{t}(y,u):=\int_t^T\L(s,y(s), u(s))\,ds+g(y(T))\tag{P$_{t,x}$}\end{equation}
\text{\textbf{Subject to:} }
\begin{equation}\label{tag:admissible}\begin{cases}
y\in W^{1,1}([t,T];\R^n)\\y'=b(y)u\text{ a.e. } s\in [t,T], \,y(t)=x\in\mathcal S\\u(s)\in \mathcal U\text{ a.e. } s\in [t,T],\, y(s)\in \mathcal S\,\,\forall s\in [t,T],
\end{cases}
\end{equation}
with the following basic assumptions.

\begin{basicass} The following conditions hold.
\begin{itemize}
\item The \textbf{Lagrangian}
$\L:[0, T]\times\R^n\times\R^m\to [0 +\infty[\cup\{+\infty\}$, $(s,y,u)\mapsto\L(s,y,u)$  is Lebesgue-Borel measurable (i.e., measurable with respect to the $\sigma$-algebra generated by products of Lebesgue measurable subsets of $[0,T]$  {\rm (}for $s${\rm )} and Borel measurable subsets of $\R^n\times\R^m$ {\rm (}for $(y,u)${\rm )};
\item $b:\R^n\to L(\R^m,\R^n)$  {\rm(}the space of linear functions from $\R^m$ to $\R^n${\rm)} is a Borel measurable function such that, for some $\theta\ge 0$,
    \begin{equation}\label{tag:assb}|b(y)|\le \theta(1+|y|).\end{equation}  We refer to $y'=b(y)u$ as to the \textbf{controlled differential equation};
\item The \textbf{control} $u:[t,T]\mapsto\R^m$ is measurable;
\item The state  constraint set $\mathcal S$ is a nonempty subset of $\R^n$;
\item The \textbf{control set} $\emptyset\not=\mathcal U\subseteq\R^m$ is a cone, i.e. if $u\in \mathcal U$ then $\lambda u\in \mathcal U$ whenever $\lambda>0$;
    \item \textbf{{\rm (}Linear growth from below{\rm )}} There are $\alpha>0$ and $d\ge 0$ satisfying, for a.e. $s\in[0,T]$ and every $y\in\R^n, u\in \mathcal U$,
\begin{equation}\label{tag:lingrowth}\L(s,y,{u})\ge \alpha|{u}|-d.\end{equation}
\item The  {\bf effective domain} of $\L$ is
   \[\Dom (\L):=\{(s,y,u)\in [0,T]\times\R^n\times\R^m:\, \L(s,y,u)<+\infty\}.\]
We assume that for a.e. $s\in [0, T]$ and every $y\in\R^n$ the set \[\{u\in\R^m:\,(s,y,u)\in\Dom (\L)\}\]    is \textbf{strictly star-shaped on the variable $u$ w.r.t.  the origin}, i.e.,
\begin{equation}\L(s,y,{u})<+\infty, \, 0<r\le 1\Rightarrow \L(s,y,r{u})<+\infty.\label{tag:starshaped}\end{equation}
\item The \textbf{ cost function}
$g:\mathcal S\to [0, +\infty[\cup \{+\infty\}$ is a given positive function, not identically equal to $+\infty$. Notice that we allow  $g$ to take the value $+\infty$ so that the class of problems studied here contains those with a end-point constraint.
\end{itemize}
\end{basicass}\noindent
\begin{remark}
Notice that it is not required that $(s,y,0)\in\Dom(\L)$ for some $(s,y)\in [0,T]\times\R^n$.\end{remark}
An \textbf{admissible pair} for {\rm (P}$_{t,x}${\rm )} is a pair of functions $(y,u):[t,T]\to\R^n\times\R^m$ with $u$ measurable, $(y,u)$ satisfying \eqref{tag:admissible} and such that $J_t(y,u)<+\infty$.
We assume henceforth that, for each $t\in [0, T[$ and $x\in\mathcal S$,  there exists at least an admissible pair  for {\rm (P}$_{t,x}${\rm )}.\\
A \textbf{minimizing sequence} $(y_j,u_j)_j$ for {\rm (P}$_{t,x}${\rm )} is a sequence of admissible pairs such that
\[\lim_{j\to +\infty}J_t(y_j,u_j)=\inf {\rm(}P_{t,x}{\rm )}.\]
\noindent
\\
Notice, that in the particular case where,  the function  $b$ is the identity matrix in the controlled differential equation,  then {\rm (P}$_{t,x}${\rm )} becomes a problem of the \textbf{Calculus of Variations}.\\
If $z\in\R^k$ we shall denote by $B^k_{r}(z)$ (simply $B^k_{r}$ if $z=0$) the closed  ball  of center $z$ and radius $r\ge 0$ in $\R^k$. The norm in $L^1$ is denoted by $\|\cdot\|_1$, and the norm in $L^{\infty}$ by $\|\cdot\|_{\infty}$.\\
If $(s,y,u)\in \Dom(\L)$ we shall denote by $\dist((s,y,u),\partial \Dom(\L))$ the euclidian \textbf{distance} from $(s,y,u)$ to the \textbf{boundary} of  $\Dom(\L)$ in $[0,T]\times \R^n\times \R^m$.\\
We will denote by $|\cdot|$ both the norm in Euclidean spaces and the Lebesgue measure in $\R$; the distinction will be clear  from the context.
\section{Assumption {\rm (A)} and Condition {\rm (S)}}\label{sec:AS}
In what follows, we assume the structure Assumption {\rm (A)} on $\L(s,y,u)$ with respect to $u$, either radial convexity or partial differentiability,  and the local Lipschitz condition {\rm (S)} on $\L(s, y,u)$ with respect to $s$.
\subsection{Assumption {\rm (A)}}
We assume henceforth the following structure condition on $\L(s,y, \cdot)$.\begin{AssumptionA}
At least one  of the two following assumptions holds:
\begin{itemize}
\item[{\rm A}$_c${\rm )}] \textbf{{\rm(}Radial convex case{\rm)}} For a.e. $s\in [0, T]$ and every $y\in\R^n,  u\in \mathcal U$, the map $0<r\mapsto\L(s,y,r{u})$ is convex, {\em or}
\item[{\rm A}$_d${\rm )}] \textbf{{\rm(}Partial differentiable case{\rm)}} For a.e. $s\in [0, T]$ and every $y\in\R^n$, the map $\L(s,y,\cdot)$ has the partial derivative
    \[D_u\L(s,y,u):= \lim_{h\to 0}\dfrac{\L(s,y, u+hu)-\L(s,y,u)}{h} \]
     with respect to $u$  at every point $u\in \mathcal U$ with $(s,y,u)\in \Dom(\L)$.
\end{itemize}
\end{AssumptionA}
\begin{remark} Notice that for $u=0$,  $0<r\mapsto\L(s,y,r{u})=\L(s,y,0)$ is convex and  $D_u\L(s,y,u)=D_0\L(s,y,0)=0$ exists, whenever $(s,y,0)\in\Dom(\L)$.
\end{remark}

\subsection{Condition {\rm (S)}}\label{sect:S}
We will consider the following local Lipschitz condition on the Lagrangian $\L$ with respect to the time variable.
\begin{conditionS}
There are $\kappa,  A\ge 0, \gamma\in L^1[0,T]$,
$\varepsilon_*>0$ satisfying,
for a.e.  $s\in [0,T]$
 \begin{equation}\label{tag:H3}
|\Lambda(s_2,y,u)-\Lambda(s_1,y,u)|\le \big(\kappa\L(s,y,u)+A|u|+\gamma(s)\big)\,|s_2-s_1|
\end{equation}
whenever $s_1,s_2\in [s-\varepsilon_*,s+\varepsilon_*]\cap [0,T]$, $y\in\R^n$, $u\in\mathcal U$, $(s_1,y,u), (s_2,y,u)\in\Dom(\Lambda)$.
\end{conditionS}
Condition {\rm (S)} was  considered  in \cite{BM2, BM1, Clarke1993}. It is a nonsmooth extension of Condition {\rm (S)}, that appears in \cite{Cesari} to establish the validity of the Du Bois-Reymond -- Erdmann equation in the smooth setting.
\begin{remark}\label{rem:auton55}
If $\L(s,y,u)=\L(y,u)$ is {\em autonomous} then Condition (S) is fulfilled with
 $\kappa=A=0$,$\gamma\equiv 0$ and $\varepsilon_*=T$.
\end{remark}
We show now that Condition {\rm (S)}
 is satisfied if $s\mapsto\L(s,y,u)$ fulfills a suitable growth condition. For $y,u\in\R^n$ we denote by $\partial_{s}^P\Lambda(s,y,u)$ the {\em proximal subgradient} of $\tau\mapsto \Lambda(\tau,y,u)$ at $\tau=s$; it coincides with $D_s\L(s,y,u)$ (resp. the convex subgradient of $\tau\mapsto\L(\tau, y, u)$ at $\tau=s$) if $\Lambda(\cdot,y,u)$ is $C^2$ (resp. convex). We refer to \cite{ClarkeBook} for more details on the subject.
\begin{proposition}[\textbf{A  proximal sufficient condition for the validity of Condition {\rm (S)}}]\label{prop:newnewnewnew}
Assume that $\L$ is {\em real valued} and that:
\begin{itemize}
\item[a)]   The map $s\mapsto\L(s,y,u)$ is lower semicontinuous for every $y\in\R^n$, $u\in\mathcal U$;
\item[b)] There is $\beta\ge 0,$  such that, for all  $s\in [0,T]$, $y\in\R^n$, $u\in\mathcal U$:
\begin{equation}\label{tag:H1semplicexxxxynew}
|\partial_{s}^P\Lambda(s,y,u)|\le \beta\big(\Lambda(s,y,u)+ |u|+1\big). \end{equation}
\end{itemize}
Then
$\L$ satisfies Condition {\rm (S)}.
\end{proposition}
\begin{proof}
For $y\in\R^n$ and $u\in\mathcal U$  consider the function
$$
h(s):=\Lambda(s,y,u)+|u|+1,\quad s\in  [0,T].
$$
Condition \eqref{tag:H1semplicexxxxynew} may be rewritten as
$|\partial^Ph(s)|\le \beta\,h(s)$ for all $s\in [0,T]$. By applying
\cite[Proposition 9.1]{BM1}, it follows  that:
\begin{itemize}
\item  If  $\beta=0$, then $h$ is constant on $[0,T]$ and thus $\L$ is autonomous;
\item If $\beta>0$, then
 \begin{equation}\label{tag:two}
|h(s_2)-h(s_1)|\le \beta' \,h(s)\,|s_2-s_1|\quad\forall s, s_1, s_2\in [0,T],
 \end{equation}
where
$$\beta' := (e^{2\beta \,T}+1)\dfrac{e^{2\beta\,T}-1}{2\,T}.$$
\end{itemize}%
In both cases  it turns out that
 $\L$ satisfies  Condition   {\rm (S)}.
 \end{proof}
\section{Growth conditions}\label{sect:growths}
We introduce here the growth Conditions {\rm (G)}, (H$_B^{\delta}$), (M$_B^{\delta}$) that are less restrictive than superlinearity. We leave some examples and proofs to Section~\ref{sect:growthsdeep}.
\subsection{Partial derivatives and subgradients}
In what follows we often deal with  subdifferentials in the sense of convex analysis.
\\
\begin{Notation}
If $(s,y,u)\in \Dom (\L)$, we shall denote by
\begin{itemize}
\item  $\partial_{\mu}\left(\L\Big(s,y,\dfrac{{u}}{\mu}\Big){\mu}\right)_{{\mu}=1}$ the \textbf{convex subdifferential}  of the map  $0<{\mu}\mapsto \L\Big(s,y,\dfrac{{u}}{\mu}\Big){\mu}$ at ${\mu}=1$, given by
 \[
 \left\{p\in\R:\,\forall\mu>0\quad \L\Big(s,y,\dfrac{{u}}{\mu}\Big){\mu}-\L(s,y,u)\ge p(\mu-1)\right\};\]
    \item
 $\partial_r\L\big(s,y,r{u}\big)_{r=1}$ the \textbf{convex subdifferential}  of the map  $0<r\mapsto \L(s,y,r{u})$ at $r=1$ given by
 \[
 \left\{q\in\R:\,\forall r>0\quad \L(s,y,ru)-\L(s,y,u)\ge q(r-1)\right\}.\]
 \item
 $D_v\L(s,y,u)$ the \textbf{partial derivative} of $\L(s,y,\cdot)$ at $u$, with respect to a given vector $v\in\R^m$;
 \item
 $\nabla_u\L(s,y,u)$ the \textbf{gradient} of $\L(s,y,\cdot)$ at  $u$. Notice that if $\L(s,y,\cdot)$ is differentiable then $D_u\L(s,y,u)=u\cdot\nabla_u\L(s,y,u)$.
 \end{itemize}
\end{Notation}

\begin{remark}\label{rem:variediff} Let $(s,y,u)\in \Dom(\L)$. the following points will be used in the sequel.
\begin{enumerate}
\item It is easy to show that if $0<r\mapsto\L(s,x,r u)$ is convex then $0<\mu\mapsto \L\Big(s,y,\dfrac{{u}}{\mu}\Big)\mu$ is convex (see, for instance, \cite[Theorem 11.5.1]{VinterBook}). Furthermore, a simple change of variable $r=\dfrac1{\mu}$ shows in this case that
\[p\in \partial_{\mu}\Big(\L\Big(s,y,\dfrac{{u}}{\mu}\Big){\mu}\Big)_{{\mu}=1}\Leftrightarrow
\Lambda(s,y,{u})-p\in \partial_r\L\big(s,y,r{u}\big)_{r=1}.\]
\item  If $D_u\L(s,y,u)$ exists  then
\[\dfrac{d}{d\mu}\Big[\L\Big(s,y,\dfrac {u}{\mu}\Big)\mu\Big]_{\mu=1}=\L(s,y,u)-D_u\L(s,y,u).\]
In particular, if $\L(s,y,\cdot)$ is differentiable at $u$ then
\[\dfrac{d}{d\mu}\left[\L\Big(s,y,\dfrac u{\mu}\Big)\mu \right]_{\mu=1}=
\L(s,y,u)-u\cdot \nabla_u\L(s,y,u).
 \]
\item If for some $\mu>0$, $D_u\Big(s,y,\dfrac{u}{\mu} \Big)$ exists, then
\[\begin{aligned}\label{tag:changevar2}\dfrac{d}{d\mu}\Big[\L\Big(s,y,\dfrac {u}{\mu}\Big)\mu\Big]_{\mu}&=\L\Big(s,y,\dfrac{u}{\mu}\Big)-D_{\frac u{\mu}}\L\Big(s,y,\dfrac{u}{\mu}\Big)\\
&=\dfrac{d}{d\lambda}\Big[\L\Big(s,y,\dfrac {u/\mu}{\lambda}\Big)\lambda\Big]_{\lambda=1}.\end{aligned}\]
 \end{enumerate}
\end{remark}
\subsection{The Growth Condition {(G)}}
The  growth assumptions introduced below involve some uniform limits.
\begin{definition}If $\phi:\dom(\L)\to\R$ is a function,  $\rho, K\ge 0$ and $E\subseteq \Dom(\phi)$ we write that
\[\lim_{\substack{|{u}|\to +\infty\\ (s,y,u)\in\,E,\, u\in \mathcal U}}\phi(s,y,u)=-\infty\,\,\,\text{unif.  } |y|\le K\]
if  for all $M\in\R$  there exists $R>0$ such that
\[\phi(s,y,u)\le M\qquad \forall (s,y,u)\in\,E,\,s\in [0, T],\,  |y|\le K,\, u\in \mathcal U,\, |u|\ge R.\]
\end{definition}
The growth Condition (G)  was thoroughly studied by Cellina and his school for autonomous Lagrangians of the calculus of variations that are smooth or convex in the velocity variable. The extension to the radial convex case, recalled here, was considered in \cite{MTLip} in the autonomous case and was introduced in \cite{BM1, BM3} for the nonautonomous case.
The growth Condition {\rm (G)} below subsumes the validity of the structure Assumption {\rm (A)}.
The partial differentiable case is new.
\begin{GrowthGloc} We say that $\L$ satisfies {\rm (G)} if (just) one of the following assumptions holds. Either:
\begin{itemize}
\item \textbf{(Radial convex case)}
$\L$ satisfies Assumption {\rm A}$_c${\rm )}.  Moreover,  there is a
 selection $Q(s,y,{u})$ of  the convex subgradient  $\partial_r\L(s, y, r{u})_{r=1}$ such that, for all $K\ge 0$
\begin{equation}\label{tag:Gequiv}
\lim_{\substack{|{u}|\to +\infty\\ (s,y,u)\in\,\Dom(\L),\, u\in \mathcal U}}\L(s,y, {u})-Q(s,y,{u})=-\infty\,\,\,\text{unif.  }  |y|\le K.\end{equation}
Equivalently,  there is a  selection $P(s,y,{u})$ of $\partial_{\mu}\Big(\L\Big(s, y, \dfrac{{u}}{\mu}\Big ){\mu}\Big)_{{\mu}=1}$,
\begin{equation}\label{tag:GequivP}
\lim_{\substack{|{u}|\to +\infty\\ (s,y,u)\in\,\Dom(\L),\, u\in \mathcal U}}P(s,y,{u})=-\infty\,\,\,\text{unif.   }  |y|\le K. \end{equation}
\end{itemize}
Or, alternatively, the following differentiable condition holds.
\begin{itemize}
\item  \textbf{(Partial differentiable case)} $\L$ satisfies Assumption {\rm A}$_d${\rm )}. Moreover,
    \begin{equation}\label{tag:Gdiff}\displaystyle
    \lim_{\substack{|{u}|\to +\infty\\ (s,y,u)\in\,\Dom(\L),\, u\in \mathcal U}}\L(s,y,u)\!-\!D_{u}\L(s,y,u)=-\infty \text{ unif.  }  |y|\le K.\end{equation}
\end{itemize}
\end{GrowthGloc}
\begin{remark}\label{rem:derivative}
\begin{enumerate}
\item  If $u\mapsto\L(s,y,u)$  is convex,  denoting by $\L^*(s,y,p)$ its Legendre transform defined by
    \[\forall p\in\R^m\quad \L^*(s,y,p)=\sup_{v\in\R^m}p\cdot v-\L(s,y,v),\]
then \eqref{tag:Gequiv} is satisfied whenever there is  a selection $p(s,y,u)$ of the convex subdifferential of $v\mapsto \L(s,y,v)$  at $v=u$ satisfying
\[\lim_{\substack{|{u}|\to +\infty\\ (s,y,u)\in\,\Dom(\L),\, u\in \mathcal U}}
\L^*(s,y,p(s,y,u))
=+\infty \text{ unif.  }  |y|\le K.\]
\item If $u\mapsto\L(s,y,u)$ is differentiable, \eqref{tag:Gdiff} becomes
\begin{equation}\label{tag:Gdiff33}\displaystyle\lim_{\substack{|{u}|\to +\infty\\ (s,y,u)\in\,\Dom(\L),\, u\in \mathcal U}}\L(s,y,u)-u\cdot \nabla_{u}\L(s,y,u)=-\infty
\text{ unif. }  |y|\le K.\end{equation}
\item If both assumptions {\rm A}$_c${\rm )} and {\rm A}$_d${\rm )} hold, the radial convex and the partial differential cases in Condition {\rm (G)} are equivalent: indeed in this situation \[\partial_r\L(s,y,ru)_{r=1}=\{ D_u\L(s,y,u)\}.\]
    \item Let $(s,y,u)\in\Dom(\L)$. In the {\em partial differentiable case}, the existence of the partial derivatives $D_u\L(s,y,u)$ at $(s,y,u)\in\Dom(\L)$ implies, by definition, that  there is $\lambda_u>0$ such that $(s,y,\lambda u)\in\Dom(\L)$ whenever $\lambda\in [1-\lambda_u, 1+\lambda_u]$.
\end{enumerate}
\end{remark}
\begin{remark}[\textbf{Interpretation of   {\rm (G)}}]\label{rem:equivG} Consider the {\em radial convex case}, the {\em partial differentiable case} being similar. Let $\L(s,y,{u})<+\infty$
 and $Q(s,y,{u})\in  \partial_{r}\L(s,y,r{u})_{r=1}$.
Then
\[\L(s,y,r{u})\ge\phi(r):=\L(s,y,{u})+ Q(s,y,{u})(r-1)\quad \forall r>0.\]
The value
$\phi(0)=P(s,y,{u}):=\L(s,y,{u})- Q(s,y,{u})$ represents the intersection with the $z$ axis of the ``tangent'' line $z=\phi(r)$ to $0<r\mapsto \L(s,y,r{u})$ at $r=1$.
Condition {\rm (G)}  thus means that the ordinate $P(s,y, {u})$ of the above intersection  tends to $-\infty$ as $|{u}|$ goes to $\infty$.
\begin{figure}[htbp]
\begin{center}
\includegraphics[width=0.65\textwidth]{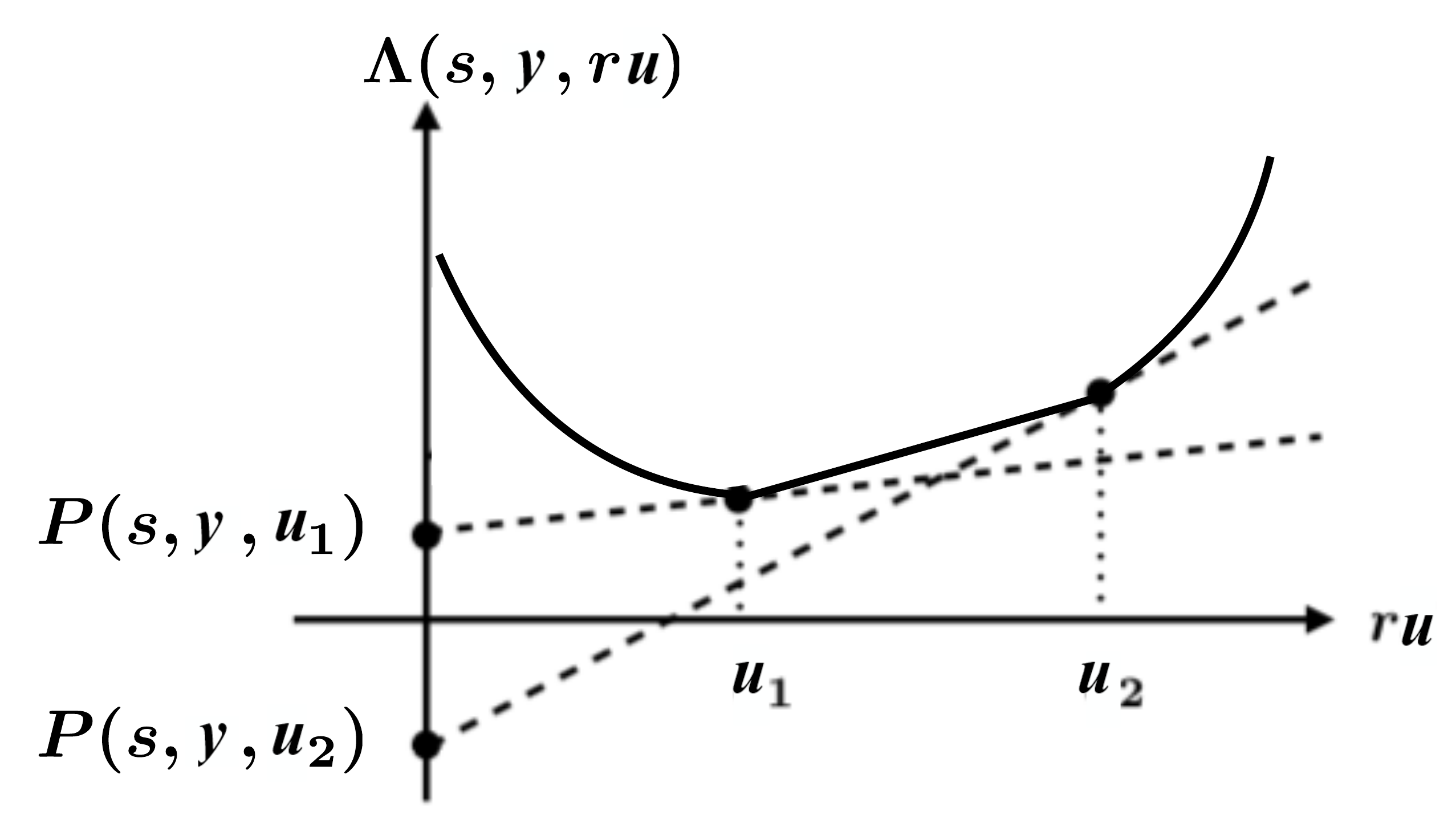}\
\end{center}
\caption{Condition {\rm (G)}.}
\label{fig:G}
\end{figure}
\end{remark}
\begin{example}
Let
\[\forall (s,y,u)\in [0,1]\times\R^2\qquad \L(s,y,u):=h(s,y)\big(|u|-\sqrt{|u|}\big),\]
where $h\ge 0$ is Borel and bounded on bounded sets.
Then $\L$ satisfies Condition (G). Indeed, $0<r\mapsto h(s,y)(r|u|-\sqrt{r|u|})$ is convex for all $u\in\R$ and, for all $u\not=0$,
\[\L(s,y,u)-u \dfrac{d}{du}\L(s,y,u)=-h(s,y)\dfrac{\sqrt{|u|}}{2}\to-\infty\]
as $|u|\to +\infty$ uniformly for $s\in [0,1]$ and $y$ in bounded sets.
\end{example}
The next proposition was formulated for the autonomous case in \cite{CF1} under the stronger assumption that $\L(y,u)$ is either convex or differentiable in $u$. Its proof is postponed to Section~\ref{sect:growthsdeep}.
\begin{proposition}[\textbf{Condition {\rm (G)} implies linear growth}]\label{prop:atleastlinear} Assume that  $\L$ fulfils  Condition {\rm (G)}.
Then $\L$ has a linear growth from below, i.e., there are  $\alpha>0$ and $d\in\R$ such that \eqref{tag:lingrowth} holds  for a.e. $s\in [0,T]$ and every $y\in\R^n, u\in \mathcal U$.
\end{proposition}
Superlinearity plays a key role in Tonelli's existence theorem. It has been widely used as a sufficient condition for Lipschitz regularity of minimizers (see \cite{AAB, CVTrans, DMF}).
\begin{AssumptionG}
There exists $\Theta:[0, +\infty[\to \R$ such that
\[\forall (s,y,u)\in\Dom(\L)\qquad \Lambda(s,y,{u})\ge\Theta(|{u}|)\quad
\displaystyle\lim_{r\to +\infty}\dfrac{\Theta(r)}{r}=+\infty .\tag{G$_\Theta$}\label{tag:superlinearity}
\]
\end{AssumptionG}
If $\L(s,y,\cdot)$ is radially convex then superlinearity, together with some local boundedness condition, imply the validity of the growth Condition {\rm (G)}. We refer to \cite[Proposition 2]{BM2} for the proof of the following result.
\begin{proposition}[\textbf{Superlinearity  $\Rightarrow$ {\rm \textbf{(G)}}}]\label{prop:SuperimpliesG}Let $\L$ be superlinear, radially convex  and assume that
there is $r_0>0$ such that $(s,y,u)\in\Dom(\L)$ whenever $s\in [0,T], y\in\R^n$ and $u\in\R^m$ with $|u|\le r_0$.
 Then
$\L$ satisfies
Assumption {\rm (G)}.
\end{proposition}
\subsection{Growth Condition {(H$_B^{\delta}$)}}\label{subsect:H}
When $B$ is an upper bound of a prescribed family of admissible pairs, with initial time $t$ varying in $[0, \delta]$, the following quantities $c_{t}(B)$ and $\Phi(B)$ will play a role in the proof of the main results.
\begin{definition}[\textbf{$c_{t}(B)$ and $\Phi(B)$}]\label{def:cphi}
Let $t\in [0,T[$, $B\ge  0$ and
assume the linear growth from below \eqref{tag:lingrowth}, i.e., for a.e. $s\in [0,T]$, for all $y\in\R^n, u\in\mathcal U$,
 \[\Lambda(s,y,u)\ge \alpha|{u}|-d\quad (\alpha>0, d\ge 0).\]
 Let
\[c_{t}(B):=\dfrac{B+d(T-t)}{\alpha\,(T-t)}.\]
Moreover, if   Condition {\rm (S)} holds, we  define
\begin{equation}
\Phi(B):=\kappa B+\dfrac A{\alpha}(B+d\,T)+\|\gamma\|_{1},
\end{equation}
where we set $\kappa, A, \gamma$ equal to 0 if $\L$ is {\em autonomous}.
\end{definition}
\begin{remark}\label{rem:cB}
Notice that, in Definition~\ref{def:cphi},  $t\in [0,T[\mapsto c_t(B)$  and $0\le B\mapsto c_t(B)$ are increasing.
\end{remark}
The next result highlights the roles  of $\Phi(B)$ and $c_t(B)$ and is a key tool in the proof of Theorem~\ref{thm:main2}.
\begin{proposition}[\textbf{The roles of $\Phi(B)$ and $c_t (B)$}] \label{lemma:newC} Assume the linear growth from below \eqref{tag:lingrowth} and the validity of Condition {\rm (S)}.
Let $t\in [0,T[$, $x\in\R^n$, $(y,u)$ admissible for {\rm (P}$_{t,x}${\rm )} with $J_t(y,u)\le B$ for some $B\ge 0$.
Then
\begin{enumerate}
\item
\begin{equation}\label{tag:100}\int_t^T|u(s)|\,ds\le \dfrac{B+d(T-t)}{\alpha}=(T-t)c_t(B).\end{equation}
\item For every $\sigma>c_{\delta}(B)$
    $$|\{s\in [t,T]:\, |u(s)|<\sigma\}|\ge \left(1-\dfrac{c_{\delta}(B)}{\sigma}\right)(T-t).$$
\item
$
\displaystyle\int_t^T\big\{\kappa\L(s,y(s),u(s))+A|u(s)|+\gamma(s)\big\}\,ds\le \Phi(B).
$
\end{enumerate}
\end{proposition}
\begin{proof}
1. Condition~\eqref{tag:lingrowth} and the fact that $g\ge 0$  yield
\[\begin{aligned}\int_t^T|u(s)|\,ds&\le \dfrac1{\alpha}\left(\int_t^T\L(s,y(s), u(s))\,ds+d(T-t)\right)\\
&\le \dfrac{B+d(T-t)}{\alpha}.\end{aligned}\]
\noindent
2. Let $\Omega=\{s\in [t,T]:\, |u(s)|<\sigma\}$. It follows from Point (1) that
\[(T-t)c_{\delta}(B)\ge (T-t)c_t(B)\ge\int_{[t,T]\setminus \Omega}|u(s)|\,ds\ge \sigma |[t,T]\setminus \Omega|,\]
implying that
$(T-t)c_{\delta}(B)\ge \sigma(T-t)-\sigma|\Omega|,$
whence the claim.
\\
3. It is enough to notice  that, from \eqref{tag:100},
\[\begin{aligned}\int_t^T\big\{\kappa\L(s,y(s),u(s))+A|u(s)|+\gamma(s)\big\}\,ds&\le
\kappa B+A\,\dfrac{B+d(T-t)}{\alpha}+\|\gamma\|_{1}\\
&\le \kappa B+A\,\dfrac{B+d\,T}{\alpha}+\|\gamma\|_{1}=\Phi(B).\end{aligned}\]
\end{proof}

Given $B\ge 0$ and $\delta\in [0,T[$, the growth Condition (H$_B^{\delta}$) below subsumes the validity of Assumption {\rm (A)} as well as of Condition {\rm (S)}. It will be applied in Theorem~\ref{thm:main2} when $B$ is an upper bound for the values of a given set of admissible pairs for problems {\rm (P}$_{t,x}${\rm )} as $t\in [0, \delta]$. In the autonomous case $\Phi(B)$ is assumed without restriction to be equal to 0, for then we may take $\kappa=A=0$ and $\gamma\equiv 0$ (see Remark~\ref{rem:auton55}).
\begin{ConditionH}
Assume that $\L$ satisfies Condition {\rm (S)}.
Let  $0\le\delta<T$, $B\ge0$ and $t\in [0, \delta]$.
We say that $\L$ satisfies {\rm (H$_{B}^{\delta}$)} if  for all $K\ge 0$ there are $\overline\nu>0$ and ${c}>c_{\delta}(B)$
  satisfying {\rm(}just{\rm)} one of the following assumptions. Either:
\begin{itemize}
\item \textbf{{\rm(}Radial convex case{\rm)}}
$\L$ satisfies Assumption {\rm A}$_c${\rm )}.
   There is a selection $Q(s,y,u)$ of $\partial_r\L(s, y, ru)_{r=1}$
  satisfying, for all  $\rho>0$,
\begin{align}\label{tag:H}
\sup_{\substack{s\in [0,T],|y|\le K\\|u|\ge \overline\nu, u\in \mathcal U\\ \L(s, y, u)<+\infty}}
\{\L(s, y, u)-Q(s,y,u)\}+2\Phi(B)<\phantom{AAAAAAA}\\
\phantom{AAAA}<\inf_{\substack{s\in [0,T],|y|\le K\\|u|<c, u\in \mathcal U\\ \L(s, y, u)<+\infty\\ \dist((s,y,u),\partial\Dom(\L))\ge\rho}}
%
\left\{\L(s, y, u)-Q(s,y,u)\right\}.\end{align}
Equivalently,  there is a  selection $P(s,y,{u})$ of $\partial_{\mu}\Big(\L\Big(s, y, \dfrac{{u}}{\mu}\Big ){\mu}\Big)_{{\mu}=1}$ satisfying
\begin{equation}\label{tag:Hequiv}
\sup_{\substack{s\in [0,T],|y|\le K\\|u|\ge \overline\nu, u\in \mathcal U\\ \L(s, y, u)<+\infty}}
\{P(s,y,u)\}+2\Phi(B)<\inf_{\substack{s\in [0,T],|y|\le K\\|u|<c, u\in \mathcal U\\ \L(s, y, u)<+\infty\\ \dist((s,y,u),\partial\Dom(\L))\ge\rho}}
P(s,y,u).\end{equation}
\end{itemize}
Or, alternatively,
\begin{itemize}
\item \textbf{{\rm(}Partial differentiable case{\rm)}} $\L$ satisfies Assumption {\rm A}$_d${\rm )}.
    For all $\rho>0$,
\begin{align}\label{tag:Hdicff}
\sup_{\substack{s\in [0,T],|y|\le K\\|u|\ge \overline\nu, u\in \mathcal U\\ \L(s, y, u)<+\infty}}
\{\L(s,y,u)-D_{u}\L(s,y,u)\}+2\Phi(B)<\phantom{AAAAA}\\
<\inf_{\substack{s\in [0,T],|y|\le K\\|u|<c, u\in \mathcal U\\ \L(s, y, u)<+\infty\\ \dist((s,y,u),\partial\Dom(\L))\ge\rho}}
\{\L(s,y,u)-D_{u}\L(s,y,u)\}.\end{align}
\end{itemize}
\end{ConditionH}
\begin{remark}\label{rem:HH}
\begin{enumerate}
\item Taking into account Remark~\ref{rem:variediff}, \eqref{tag:Hdicff} is equivalent to
\begin{align}\label{tag:Hdiff3}
\sup_{\substack{s\in [0,T],|y|\le K\\|u|\ge \overline\nu, u\in \mathcal U\\ \L(s, y, u)<+\infty}}
\dfrac d{d\mu}\left[\L\left(s, y, \dfrac{u}{\mu}\right)\mu\right]_{\mu=1}+2\Phi(B)<\phantom{AAAAAAA}\\
<\!\!\!\!\!\!\inf_{\substack{s\in [0,T],|y|\le K\\|u|<c, u\in \mathcal U\\ \L(s, y, u)<+\infty\\ \dist((s,y,u),\partial\Dom(\L))\ge\rho}}
\dfrac d{d\mu}\left[\L\left(s, y, \dfrac{u}{\mu}\right)\mu\right]_{\mu=1}.\end{align}
    \item When $\L(s,y,\cdot)$ is differentiable, in \eqref{tag:Hdicff} we have
    \[\L(s,y,u)-D_{u}\L(s,y,u)=\L(s,y,u)-u\cdot \nabla_{u}\L(s,y,u).\]
\item With the above notation, Condition (H$_B^{\delta}$) is satisfied independently of what the value $\Phi(B)$ is if, for some $\overline\nu>0$ and $c>c_{\delta}(B)$:
    \begin{itemize}
    \item \textbf{Radial convex case:} For all $K\ge 0, \rho>0$,
   \begin{equation}\label{tag:HequivTRIS}
\lim_{\nu\to +\infty}
\sup_{\substack{s\in [0,T],|y|\le K\\|u|\ge \nu, u\in \mathcal U\\ \L(s, y, u)<+\infty}}
P(s,y,u)=-\infty; \,\inf_{\substack{s\in [0,T],|y|\le K\\|u|<c, u\in \mathcal U\\ \L(s, y, u)<+\infty\\ \dist((s,y,u),\partial\Dom(\L))\ge\rho}}
P(s,y,u)\not=-\infty.\end{equation}
or,
    \item \textbf{Partial differentiable case:} For all $K\ge 0, \rho>0$,
 \begin{align}\label{tag:HdiffTRIS}
\lim_{\nu\to +\infty}
\sup_{\substack{s\in [0,T],|y|\le K\\|u|\ge \nu, u\in \mathcal U\\ \L(s, y, u)<+\infty}}
\{\L(s,y,u)-D_{u}\L(s,y,u)\}=-\infty;\\
\inf_{\substack{s\in [0,T],|y|\le K\\|u|<c, u\in \mathcal U\\ \L(s, y, u)<+\infty\\ \dist((s,y,u),\partial\Dom(\L))\ge\rho}}
\{\L(s,y,u)-D_{u}\L(s,y,u)\}\not=-\infty.\end{align}
    \end{itemize}
    \end{enumerate}
    In particular Condition (H$_B^{\delta}$) is  satisfied if, for some $\overline\nu>0$,  $|u|<\overline\nu$ whenever $(s,y,u)\in\Dom(\L)$ and the ``inf'' term in  \eqref{tag:HequivTRIS}, or in \eqref{tag:HdiffTRIS}, is not equal to $-\infty$.
\end{remark}
The next Example~\ref{ex:4567} illustrates the importance of Condition (H$_B^{\delta}$);  it is taken from \cite[Example 4.3]{Clarke1993}.
\begin{example}[The minimal length functional]\label{ex:4567} Let $n, m\ge 1$. The function \[\forall
(s,y,u)\in [0,T]\times\R^n\times \R^m\qquad \L(s,y,u)=L(u):=\sqrt{1+|u|^2}\]
satisfies Condition (H$_B^{\delta}$) for any choice of $\delta\\in [0, T[$ and $B$. Indeed here $\Phi(B)=0$. Moreover,
\[L(u)-u\cdot\nabla L(u)=\dfrac1{\sqrt{1+|u|^2}}\] so that
\[\lim_{\nu\to +\infty}
\sup_{|u|\ge \nu}
\{ L(u)-D_{u} L(u)\}+2\Phi(B)=\lim_{\nu\to +\infty}\dfrac1{\sqrt{1+\nu^2}}=0,\]
whereas, for any $c>0$,
\[\inf_{|u|< c}
\{ L(u)-D_{u} L(u)\}=\inf_{|u|< c}\dfrac1{\sqrt{1+|u|^2}}=\dfrac1{\sqrt{1+c^2}}.\]
Notice that $\L$ does not satisfy Condition {\rm (G)}, since
\[\lim_{|u|\to +\infty}
\{L(u)-D_{u}L(u)\}=\lim_{|u|\to +\infty}\dfrac1{\sqrt{1+|u|^2}}=0.\]
\end{example}
\begin{remark}\label{rem:infsup100} It is useful to clarify the quantities that appear in Condition (H$_B^{\delta}$).
\begin{enumerate}
 \item It follows from Remark~\ref{rem:cB} that if $B>B'\ge 0$ then the validity of Condition (H$_B^{\delta}$) implies that of (H$_{B'}^{\delta}$).
\item Comments on the ``inf'' part of \eqref{tag:H} -- \eqref{tag:Hdicff}:
\begin{itemize}
\item If $\L$ is real valued, the additional condition $\dist((s,y,u),\partial\Dom(\L))\ge\rho$ in  \eqref{tag:H}, \eqref{tag:Hequiv}, \eqref{tag:Hdicff} is trivially fulfilled since  $\Dom(\L)=[0,T]\times\R^n\times\R^m$.
 \item The validity of Condition (H$_B^{\delta}$)  implies that the right-hand side of \eqref{tag:Hequiv}, or \eqref{tag:Hdicff}, is not equal to $-\infty$.
\item (the role of $B$) In view of Proposition~\ref{lemma:newC}, the initial assumption on $B$ ensures that, if $(y,u)$ is admissible for  (P$_{t,x}$) ($t\le \delta$) then,  for $K\ge \|y\|_{\infty}$ and $c>c_{\delta}(B)$,  the set
\[\{(s,z,v)\in \Dom(\L):\, |z|\le K,\,v\in\mathcal U,\,|v|<c\}\]
is non empty, so that, if $\Dom(\L)$ is open in $[0,T]\times\R^n\times\R^m$,  there is $\rho>0$ in such a way that the infimum in the right-hand side of \eqref{tag:H} -- \eqref{tag:Hdicff} is not equal to $+\infty$.
\end{itemize}
\item  Comments on the ``sup'' part of \eqref{tag:H} -- \eqref{tag:Hdicff}:
\begin{itemize}
\item  Unless $|v|<\overline\nu$ whenever $(s,z,v)\in \Dom(\L)$, the sup in the left-hand side of \eqref{tag:H} -- \eqref{tag:Hdicff} is not equal to $-\infty$.
    \end{itemize}
    \end{enumerate}
\end{remark}
\begin{remark}\label{rem:newHwhy} Condition (H$_B^{\delta}$) represents a violation of the  Du Bois-Reymond -- Erdmann condition for high values of the control variable (see \cite{BM3}).  Condition {\rm (H$_{B}^{0}$)} was introduced in \cite{Clarke1993} for a fixed initial time problem $t=0$ in the following setting:
\begin{itemize}
\item  a   convex problem of the calculus of variations with a real valued Lagrangian;
\item $B$ equal to any upper bound of $J_t(y,y')$ for a suitable admissible trajectory $y: [0,T]\to \R^n$.
\end{itemize}
The present  formulated is suitable for  classes of admissible trajectories as the initial time varies in an interval $[0, \delta]$.
We point out that,  with respect to the original version, the term $2\Phi(B)$ in the right-hand side of \eqref{tag:Hequiv} replaces  $\Phi(B)$, so that our condition, in the nonautonomous case, is slightly more restrictive than \cite[Hypothesis (H2)$'$]{Clarke1993}. We do not have, however, an example where this represents a true drawback.\\
At the same time, the present version enlarges the realm of application in several new aspects: it takes into account the partial differentiable case,  variable initial time/position (which involves the choice of $B, c_{\delta}(B)$ and $\Phi(B)$ in Definition~\ref{def:cphi}) and, in \eqref{tag:H}, \eqref{tag:Hdicff}, the additional requirement that the infimum is taken just for the points of the effective domain that satisfy $\dist((s,y,u),\partial\Dom(\L))\ge\rho$ gives more chances for the condition to be satisfied. Indeed, in \cite[Conditions (H1) -- (H2)]{Clarke1993}, the requirements that the domain of $u\mapsto\L(s,y,u)$ is open and star shaped, together with the requirement that $\L$ tends uniformly to $+\infty$ at the boundary of the effective domain seem to fight against the implicit requirement of  \cite[Condition (H1)]{Clarke1993} that
\[\inf_{\substack{s\in [0,T], |y|\le K\\ |v|<c, v\in\mathcal U}}\{\L(s,z,v)-v\cdot\partial_v\L(s,z,v)\}>-\infty.\]
This fact is clarified in Example~\ref{ex:Hnew}.
\end{remark}
Next Proposition~\ref{prop:GimpliesH} shows that the infimum in \eqref{tag:HequivTRIS} -- \eqref{tag:HdiffTRIS}  is finite  under a natural assumption  related to the boundedness of $\L$ on bounded sets that are far away from the boundary of the domain.
In this situation, the validity of Condition {\rm (G)} implies that of Condition (H$_B^{\delta}$), whatever are the choices of $B$ and $0\le\delta<T$.
\begin{definition}\label{def:wellinside}  We say that a subset of $\Dom(\L)$ is
\textbf{well-inside $\Dom(\L)$} if it is contained in  $\{(s,y,u)\in\Dom(\L):\, \dist((s,y,u),\,\partial\Dom(\L))\ge \rho\}$, for a suitable $\rho>0$.
\end{definition}
\begin{remark} If the effective domain of $\L$ is open in $[0,T]\times\R^n\times\R^m$, as is the case under Assumption h$_1$) of Theorem~\ref{thm:main2}, the  notion of "well-inside" coincides  with that of relatively compact subset.
\end{remark}
\begin{proposition}[\textbf{{\rm (G)} implies (H$_B^{\delta}$) for all $B, \delta$}]\label{prop:GimpliesH} Assume that $\L$   is   \textbf{bounded on the bounded sets that are well-inside $\Dom(\L)$} and at least one of the following structure conditions:
\begin{itemize}
\item[a)] $\L$ is \textbf{radially convex in the control} variable as in {\rm A}$_c${\rm )}, {\em or},
\item[b)] $\L(s,y,\cdot)$ is \textbf{partially differentiable} as in {\rm A}$_d${\rm )} and  is  uniformly Lipschitz for $(s,y,u)$ in each bounded set that is well-inside the domain.
\end{itemize}
The following properties hold:
\begin{enumerate}
\item Let $Q(s,y,u)\in \partial_r\L(s, y, ru)_{r=1}$ if a) holds, otherwise set $Q(s,y,u):=D_u\L(s,y,u)$. Then for every $c, \rho>0$ and $K\ge 0$,
\begin{equation}\label{tag:inffinite}\inf_{\substack{s\in [0,T],|y|\le K\\|u|<c, u\in \mathcal U\\ \L(s, y, u)<+\infty\\ \dist((s,y,u),\partial\Dom(\L))\ge\rho}}
\{\L(s,y,u)-Q(s,y,u)\}>-\infty.\end{equation}
\item If $\L$ satisfies Condition (G)  and Condition {\rm (S)}   then
$\L$ satisfies Hypothesis (H$_B^{\delta}$), whatever are the choices of $\delta\in [0,T[$, $c>0$  and $B\ge 0$.
\end{enumerate}
\end{proposition}
\begin{proof}  Fix $K\ge 0$, $c, \rho >0$ and let $W_{K, c, \rho}$ be the subset of $\Dom(\L)$ defined by
\[\{(s,y,u)\in\Dom(\L):\, |y|\le K,\,u\in \mathcal U,\,|u|<c,\dist((s,y,u),\partial\Dom(\L))\ge\rho\}.\]
\\
(1) It is not restrictive to assume that $W_{K, c, \rho}\not=\emptyset$, otherwise the infimum in \eqref{tag:inffinite} equals $+\infty$.
It follows either from Lemma~\ref{tag:boundabovewell}  (radial convex case) or from the  local Lipschitzianity of $\L$ well-inside the domain  (partial differentiable case) that $Q(s,y,u)$ is bounded above on $W_{K, c, \rho}$.
The local boundedness condition on $\L$ implies that
\[\inf_{\substack{s\in [0,T],|y|\le K\\|u|<c, u\in \mathcal U\\ \L(s, y, u)<+\infty\\ \dist((s,y,u),\partial\Dom(\L))\ge\rho}}
\{\L(s,y,u)-Q(s,y,u)\}>-\infty.\]
\noindent
(2) In the partial differentiable case we set $Q(s,y,u):=D_u\L(s,y,u)$; otherwise
let $Q(s,y,u)\in \partial_r\L(s, y, ru)_{r=1}$ be  such that
\[\lim_{\substack{|{u}|\to +\infty\\ (s,y,u)\in\,\Dom(\L),\, u\in \mathcal U}}\L(s,y, {u})-Q(s,y,{u})=-\infty\,\,\,\text{ unif. }  |y|\le K.\]
Then
\begin{equation}\label{tag:iuefgquogf}
\lim_{\nu\to +\infty}\sup_{\substack{s\in [0,T]\\|u|\ge \nu, u\in \mathcal U,\, |y|\le K\\ \L(s, y, u)<+\infty}}
\{\L(s, y, u)-Q(s,y,u)\}=-\infty.
\end{equation}
It follows from \eqref{tag:inffinite} that
 Condition (H$_B^{\delta}$) is valid, for any choice of $B, c>0$, $\delta\in [0,T[$.
\end{proof}
\begin{lemma}[\textbf{Bound of $\partial_r\L(s, y, ru)_{r=1}$ on bounded sets}] \label{tag:boundabovewell}
Assume that $\L(s,y,u)$ is  \textbf{radially convex in the control variable} and \textbf{bounded on the bounded sets} that are \textbf{well-inside $\Dom(\L)$}. Let
\[\forall (s,y,u)\in\Dom(\L)\qquad Q(s,y,u)\in \partial_r\L(s, y, ru)_{r=1}.\]
Then $Q$ is bounded  on the bounded sets that are well-inside $\Dom(\L)$.
\end{lemma}
\begin{proof}
Let $(s,y,u)\in \Dom(\L)$ with $|y|+ |u|\le C$ for some $C>0$ and  \[\dist((s,y,u),\partial\Dom(\L))\ge\rho\] for some $\rho>0$. We have
\[\dist\left(\left(s,y,u+\dfrac{\rho}{2C}u\right), \partial\Dom(\L)\right)\ge \dfrac{\rho}{2}.\] Since
\[\L\left(s,y,u+\dfrac{\rho}{2C}u\right)-\L(s,y,u)\ge\dfrac{\rho}{2C} Q(s,y,u),\]
the  boundedness assumption of $\L$ implies  that $Q(s,y,u)$ is bounded above by a constant  depending only on $C$ and  $\rho$.
Similarly, from
\[\L\left(s,y,u-\dfrac{\rho}{2C}u\right)-\L(s,y,u)\ge -\dfrac{\rho}{2C} Q(s,y,u),\]
we deduce a lower bound for $Q$.
\end{proof}
\begin{remark} At a first glance the conclusion of Lemma~\ref{tag:boundabovewell} is puzzling, if one thinks at a convex function on $[0, +\infty[$ with a vertical tangent in 0. However,  $Q(s,y,u)$ stands for a subgradient of $0<r\mapsto \L(s,y,ru)$ at $r=1$ and not as a subgradient of $v\mapsto\L(s,y,v)$ at $v=u$.
\end{remark}
\begin{remark}
\begin{itemize}
\item The   local boundedness assumption in Proposition~\ref{prop:GimpliesH} is crucial in order to conclude; a counterexample is provided in Example~\ref{ex:GnotH}.
\item Conditions (G) and (H$_{B}^{\delta}$) are not equivalent. Example~\ref{ex:4567} exhibits a real valued Lagrangians $\L(s,y,u)$ that is convex in $u$, bounded on bounded sets, satisfies Condition (H$_B^{\delta}$), for some $\delta, B\ge0$, but nevertheless does not fulfill Assumption {\rm (G)}.
\item Example~\ref{ex:Hnew} illustrates the role of the condition $\dist((s,y,u),\partial\Dom(\L))\ge\rho$ in  \eqref{tag:H} -- \eqref{tag:Hdicff} and shows that in the extended valued case, the original  Condition (H) as defined  in \cite{BM3, FrankTrans, Clarke1993}  may not imply (G).
\end{itemize}
\end{remark}
\subsection{Growth Condition (M$_B^{\delta}$)}\label{sect:M}
We introduce here a new condition,than that turns out to be satisfied by a wide class of Lagrangians, as shown in Proposition~\ref{prop:whenM}.
Whereas, in the nonautonomous case, the growth hypothesis (H$_B^{\delta}$) subsumes the validity of Condition (S) through the definition of $\Phi(B)$, the growth Condition (M$_B^{\delta}$) does not involve it anymore.
\begin{ConditionHNEW}
Let $0\le\delta<T$, $B\ge0$ and $t\in [0, \delta]$.
We say that $\L$ satisfies {\rm (M$_B^{\delta}$)} if, for all $K\ge 0$,
there are $\overline\nu> 0$ and ${c}>c_{\delta}(B)$
  satisfying {\rm(}just{\rm)} one of the following assumptions. Either:
\begin{itemize}
\item \textbf{{\rm(}Radial convex case{\rm)}}
$\L$ satisfies Assumption {\rm A}$_c${\rm )}.
   There is a selection $Q(s,y,u)$ of $\partial_r\L(s, y, ru)_{r=1}$
  satisfying, for all $\rho>0$:
\[\begin{aligned}\label{tag:HNEW}
i) &-\infty<\inf_{\substack{s\in [0,T],|y|\le K\\|u|<c, u\in \mathcal U\\ \L(s, y, u)<+\infty\\ \dist((s,y,u),\partial\Dom(\L))\ge\rho}}
\left\{\L(s, y, u)-Q(s,y,u)\right\},\\
ii)  &\sup_{\substack{s\in [0,T],|y|\le K\\|u|\ge \overline \nu, u\in \mathcal U\\ \L(s, y, u)<+\infty}}
\{\L(s, y, u)-Q(s,y,u)\}<+\infty.\end{aligned}\]
Equivalently,  for a suitable  selection $P(s,y,{u})$ of $\partial_{\mu}\Big(\L\Big(s, y, \dfrac{{u}}{\mu}\Big ){\mu}\Big)_{{\mu}=1}$:

\[\begin{aligned}\label{tag:HequivNEW}
i) &-\infty<\inf_{\substack{s\in [0,T],|y|\le K\\|u|<c, u\in \mathcal U\\ \L(s, y, u)<+\infty\\ \dist((s,y,u),\partial\Dom(\L))\ge\rho}}
\left\{P(s, y, u)\right\},\\
ii) &\sup_{\substack{s\in [0,T],|y|\le K\\|u|\ge \overline\nu, u\in \mathcal U\\ \L(s, y, u)<+\infty}}
\{P(s, y, u)\}<+\infty.\end{aligned}\]
\end{itemize}
Or, alternatively,
\begin{itemize}
\item \textbf{{\rm(}Partial differentiable case{\rm)}} $\L$ satisfies Assumption {\rm A}$_d${\rm )}.
    For all $\rho>0$:
\[\begin{aligned}\label{tag:HdicffNEW}
i) &-\infty<\inf_{\substack{s\in [0,T],|y|\le K\\|u|<c, u\in \mathcal U\\ \L(s, y, u)<+\infty\\ \dist((s,y,u),\partial\Dom(\L))\ge\rho}}
\left\{\L(s,y,u)-D_{u}\L(s,y,u)\right\},\\
ii) &\sup_{\substack{s\in [0,T],|y|\le K\\|u|\ge \overline \nu, u\in \mathcal U\\ \L(s, y, u)<+\infty}}
\{\L(s,y,u)-D_{u}\L(s,y,u)\}<+\infty.\end{aligned}\]

\end{itemize}
\end{ConditionHNEW}
\begin{figure}[htbp]
\begin{center}
\includegraphics[width=0.55\textwidth]{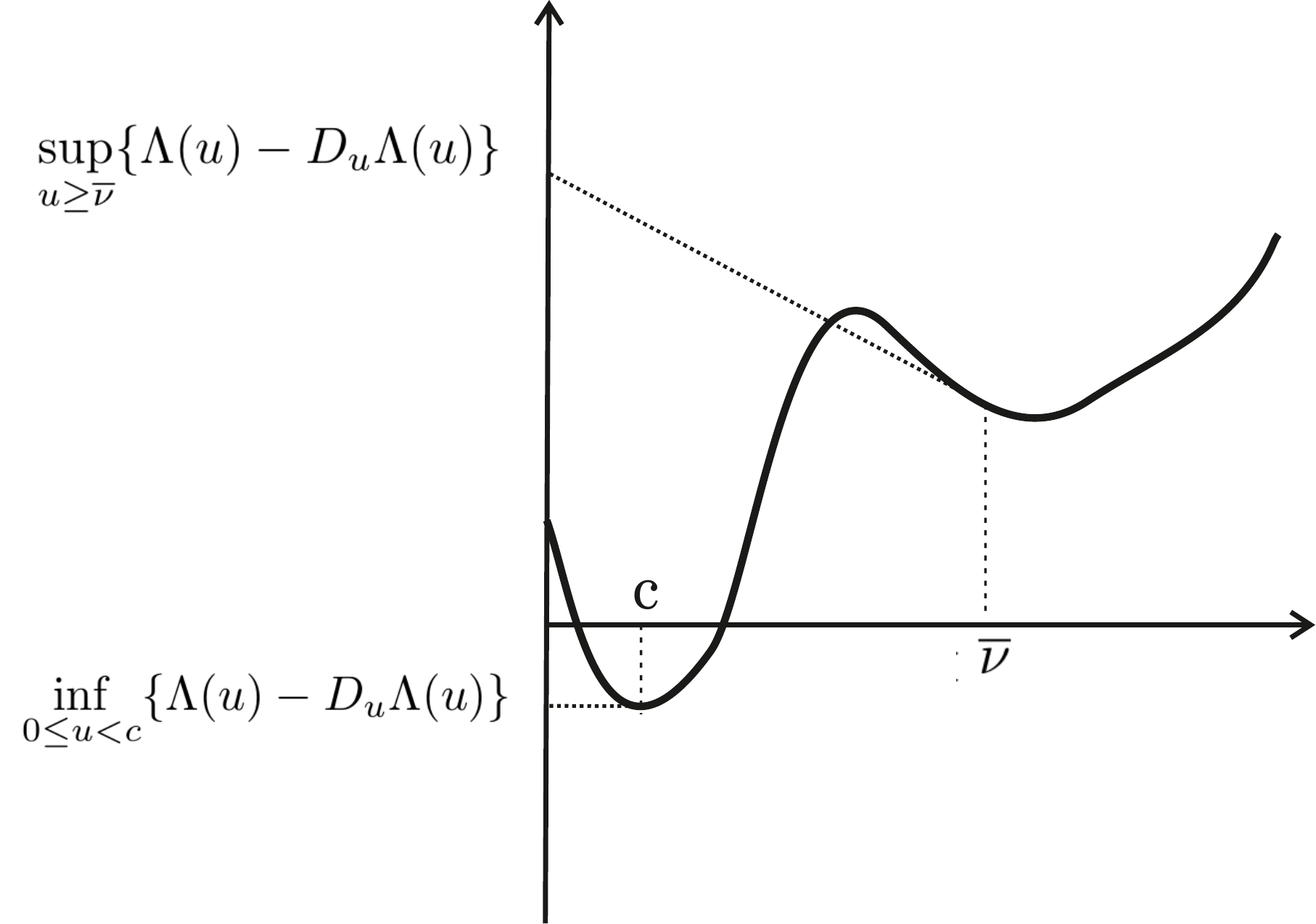}
\caption{The infimum and the sup involved in Condition of type (M): the case of a smooth function  of a positive variable.}
\end{center}
\end{figure}
\begin{remark}\label{rem:nonemptyM} \begin{enumerate}
\item Proposition~\ref{prop:GimpliesH} shows that {\em i}) of Condition (M$_B^{\delta}$)  in \eqref{tag:HNEW}, \eqref{tag:HequivNEW}, \eqref{tag:HdicffNEW} is satisfied under very mild conditions. In particular, in the radial convex case {\rm A}$_c${\rm )}, it is enough that $\L$ is bounded on the bounded sets that are well-inside the domain.
\end{enumerate}
\end{remark}
\begin{remark}\label{rem:Himplies M}  The strict inequality sign in \eqref{tag:H}, \eqref{tag:Hequiv}, \eqref{tag:Hdicff} shows that the validity of Condition (H$_B^{\delta}$) implies that of (M$_B^{\delta}$). The converse is not true, as shown by $\L(u)=|u|$ or by the Lagrangian in Example~\ref{ex:nononlyconvex}.
\end{remark}
\begin{example}[\textbf{A radially concave function that satisfies (M$_B^{\delta}$) but not (H$_B^{\delta}$)} ] \label{ex:nononlyconvex}
Let $$\L(u):=2|u|-\sqrt{1+u^2} \qquad \forall u\in\mathbb R.$$
Then
\[\forall u\in\R\qquad\L(u)-D_u\L(u)=-\dfrac{1}{\sqrt{1+u^2}}.\]
Thus, for any $c>0$ and $\nu>0$ we have
\[\sup_{|u|\ge\nu} \L(u)-D_u\L(u)=0,\quad \inf_{|u|<c} \L(u)-D_u\L(u)=-1.\]
Therefore, for any $B, \delta$, Condition (M$_B^{\delta}$) is satisfied, whereas Condition (H$_B^{\delta}$) is not. Notice that $\L$ is concave on $[0, +\infty[$ and on $]-\infty,0]$.
\begin{figure}[htbp]
\begin{center}
\includegraphics[width=0.4\textwidth]{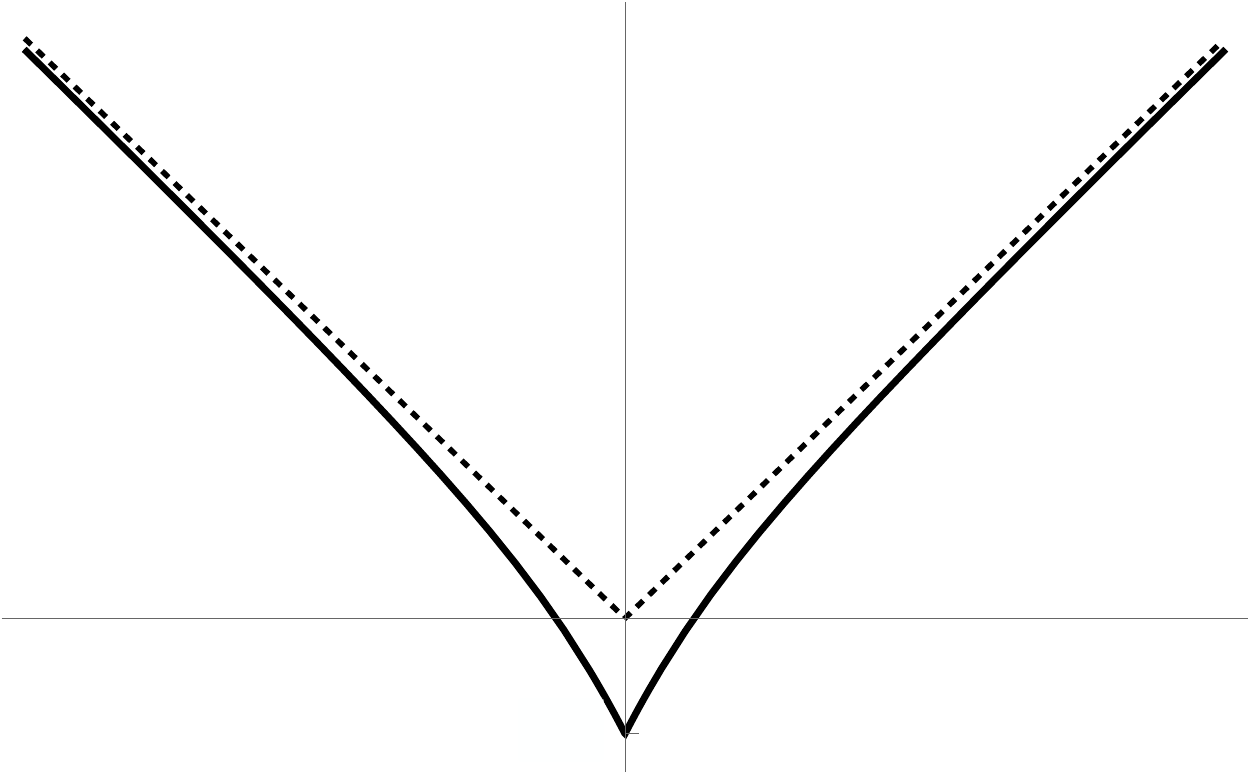}
\end{center}
\caption{The function considered in Example~\ref{ex:nononlyconvex}.}
\end{figure}
\end{example}
Condition (M$_B^{\delta}$) is satisfied by {\em every} real valued radially convex Lagrangians that is bounded on the bounded sets. The proof of Proposition~\ref{prop:whenM} is postponed to Section~\ref{sect:growthsdeep}.
\begin{proposition}[\textbf{Validity of Condition (M$_B^{\delta}$)}]\label{prop:whenM} Assume that $\L$  is \textbf{radially convex in the control} variable as in {\rm A}$_c${\rm )}. Assume, moreover, that:
\begin{itemize}
\item[a)] $\L$  is \textbf{bounded on the bounded sets} that are \textbf{well-inside $\Dom(\L)$};
\item[b)] For all $K\ge 0$ there is $r_K>0$ such that $[0,T]\times B_{K}^n\times B_{r_K}^m\subseteq\Dom(\L)$.
\end{itemize}
Then $\L$ satisfies condition (M$_B^{\delta}$) for any choice of $B\ge0$ and $\delta\in [0,T[$.\\
In particular the condition is satisfied whenever $\L(s,y,u)$ is real valued, continuous and radially convex in the control variable.
\end{proposition}
\begin{remark} Assumption b) in Proposition~\ref{prop:whenM} is a  well known sufficient condition for the nonoccurrence of the Lavrentiev gap for positive autonomous Lagrangians of the calculus of variations (see \cite[Assumption (B)]{ASC}).
\end{remark}
\section{Regularity of minimizing pairs}
 \subsection{Nice admissible pairs}
Theorem~\ref{thm:main2} is the main core of the paper. Some examples that illustrate its applications are postponed to Section~\ref{sect:examples}, whereas Section~\ref{sect:proofs} is entirely devoted to its  proof.
  \begin{theorem}[\textbf{Nice admissible pairs}]\label{thm:main2}
 Suppose that  $\Lambda$  satisfies Assumption  {\rm (A)} and  Condition {\rm (S)}. Let  $\delta\in [0,T[$, $\delta_*\ge 0$, $x_*\in\R^n$,
and  $\mathcal A$ be a family of admissible pairs for {\rm (P}$_{t,x}${\rm )}, for some $t\in [0, \delta]$ and $x\in B_{\delta_*}^n(x_*)$. Assume that $J_t(y,u)\le B$ for some $B\ge 0$, whenever $(y,u)\in\mathcal A$.
 Unless $\L$ is real valued assume, moreover, the following conditions:
 \begin{itemize}
 \item[h$_{1}$)]
    \textbf{The domain is a product}, i.e.,  $\Dom(\L)=[0,T]\times D$ for some $D\subseteq\R^n\times\R^m$.
 \item[h$_{2}$)] $\L$ \textbf{tends uniformly to $+\infty$ at the boundary of the effective domain}, i.e.,  \begin{equation}\label{tag:infinito}\displaystyle\lim_{\substack{\dist((s,z,v),\,\partial\Dom(\L))\to 0\\(s,z,v)\in\Dom(\L),\,v\in\mathcal U}}\L(s,z,v)=+\infty.\end{equation}
\end{itemize}
 The following claims hold:
\begin{enumerate}
\item Suppose that $\L$ satisfies Condition {\rm(H}$^{\delta}_B${\rm )}.
Then there is a constant $K_{\mathcal A}$ such that,
for every admissible pair $(y,u)\in\mathcal A$ for {\rm (P}$_{t,x}${\rm )} there exists an admissible pair $(\overline y, \overline u)$ for {\rm (P}$_{t,x}${\rm )} such that
 \begin{enumerate}
 \item $\overline y=y\circ \psi$, where $\psi$ is a Lipschitz reparametrization of $[t,T]$;
 \item
 $\|\overline u\|_{\infty}\le K_{\mathcal A}, \, \|\overline y\|_{\infty}=\| y\|_{\infty}\le K_{\mathcal A}$, and
$\overline y$ is Lipschitz of rank $K_{\mathcal A}$;
\item
$J_{t}(\overline{y},\overline{u} )\le  J_{t}(y, u)$, the inequality being \textbf{strict} if $u$ is not bounded.
\end{enumerate}
\item Suppose that $\L$ satisfies Condition {\rm(M}$^{\delta}_B${\rm )} and let $\eta>0$. Then the conclusions  (a), (b) of Claim (1) remain valid with $K_{\mathcal A}$ possibly depending on $\eta$ and, moreover,
    \begin{itemize}
\item[(c$'$)]     $J_{t}(\overline{y},\overline{u} )\le  J_{t}(y, u)+\eta$.
    \end{itemize}
\end{enumerate}
\end{theorem}
\begin{remark}\label{rem:propag}
\begin{enumerate}
\item Hypothesis h$_1$) implies that
 \[\forall(s,z,v)\in\Dom(\L)\quad \dist ((s,z,v),\partial\Dom(\L))=\dist((z,v), \partial D);\] in particular if $\dist ((s,y,v),\partial\Dom(\L))\ge\rho$ for some $(s,y,v)\in\Dom(\L)$ then $\dist ((\tilde s,z,v),\partial\Dom(\L))\ge\rho$ for every $\tilde s\in [0, T]$;  this fact is essential in the proof of Theorem~\ref{thm:main2}. Hypothesis h$_{1}$)   is satisfied if, for instance,  $\L(s,z,v)=\phi(s)L(z,v)$ where $\phi:[0,T]\to [0, +\infty[$ and $L:\R^n\times\R^m\to [0, +\infty]$.
\item Hypothesis h$_{2}$) implies that the effective domain $\Dom(\L)$ is open in $[0,T]\times\R^n\times\R^m$.
\item In the case of a family $\mathcal A$ reduced to single admissible pair $(y,u)$ for {\rm (P}$_{t,x}${\rm )}, the requirement of the validity of condition {\rm(H}$^{\delta}_B${\rm )} (resp. {\rm(M}$^{\delta}_B${\rm )}) in Theorem~\ref{thm:main2} may be replaced by that of the validity of (H$^{t}_{J_t(y,u)}$) (resp. (M$^{t}_{J_t(y,u)}$)).
\end{enumerate}
\end{remark}
The existence of an upper bound $B$ in Theorem~\ref{thm:main2} is ensured, obviously, if the family $\mathcal A$ is reduced to a singleton $\{(y,u)\}$, in which case $B=J_t(y,u)$ is a suitable choice. Some sufficient conditions for the existence of $B$ may be obtained when $\L$ is real valued.
\begin{lemma}[\textbf{A uniform upper bound for the infima of {\rm(P}$_{t,x}${\rm )}}]\label{lemma:unifbound} Assume that  $\Lambda$ is \textbf{finite valued} and   bounded on bounded sets. Suppose   that  one of the two assumptions holds:
\begin{enumerate}
\item either  \textbf{$b=1$ in the controlled differential equation, $\mathcal S$ is convex and $\mathcal U=\R^m$}, \emph{or}
    \item the  cost function \textbf{$g$ is real valued, locally bounded} and \textbf{$0\in\mathcal U$}.
        \end{enumerate}
Let $\delta\in [0,T[$, $\delta_*\ge 0$ and $x_*\in\R^n$. There is $B\ge 0$ such that for every $t\in [0, \delta],  x\in B_{\delta_*}^n(x_*)$, there exists an admissible pair $(y,u)$ for {\rm (P}$_{t,x}${\rm )} satisfying $J_t(y,u)\le B$.
\end{lemma}
\begin{proof}
We consider separately the cases 1 and 2.\\
(1)
Let $\xi_*\in\mathcal S$ be such that $g(\xi_*)<+\infty$.\\
Define $y(s):=\dfrac{T-s}{T-t}x+\dfrac{s-t}{T-t}\xi_*$. Then, since $y(t)=x$ and, by  the convexity of $\mathcal S$, $y$ has values in $\mathcal S$ then $(y, y')$ is admissible.
Now,
\[J_t(y,y')=\int_t^T\L(s,y,y')\,ds+g(\xi_*). \]
 Moreover, for every $s\in [t,T]$, $y(s)$ belongs to the segment joining $x$ with  $\xi_*$, and $y'(s)=\dfrac{\xi_*-x}{T-t}$.
 Since $|\xi_*-x|\le |\xi_*-x_*|+\delta_*$ and $t\le\  \delta$ we get
\[ |y(s)|\le |\xi_*-x_*|+\delta_*,\quad |y'(s)|\le \dfrac{|\xi_*-x_*|+\delta_*}{T-\delta}:\]
It follows from the boundedness assumption of $\L$ on bounded sets that there is a constant $B$, depending only on $x_*, \delta, \delta_*$ such that
\[ J_t(y,y')\le B.\]
\\
(2) Assume that $0\in\mathcal U$. Since the constant pair $(y(s):=x,u(s):=0)$ is admissible  we have
\begin{equation}\label{tag:bbound}J_t(y,u)= \int_t^T\L(s,x,0)\,ds+g(x).\end{equation}
The fact that  $\L$ is bounded on bounded sets and  that $g$ is bounded on $B_{\delta_*}^n(x_*)$ imply that the right-hand side of \eqref{tag:bbound} is bounded above by a constant depending only on $x_*, \delta, \delta_*$, whence the claim.
\end{proof}
\begin{remark}\label{rem:estimateB} The proof of  Lemma~\ref{lemma:unifbound} shows that  $B$ can be estimated in both cases from above, respectively, by the following constants:
\\
(1) $B\le  T \sup\left\{\L(s,z,v):\, s\in [0, T], \, |z|\le |\xi_*-x_*|+\delta_*,\, |v|\le \dfrac{|\xi_*-x_*|+\delta_*}{T-\delta}\right\}+g(\xi_*)$, where is an arbitrary element in $\mathcal S$ such that $g(\xi_*)<+\infty$;\\
(2)a) If $0\in\mathcal U$ then  $B\le \sup\{T\,\L(s,x,0)+g(x):\, s\in [0, T],\, x\in B_{\delta_*}^n(x_*)\}$.\\
(2)b) If $\mathcal S=\R^n$ then, given $u_*\in\mathcal U$,
\[B\le \sup\{T\,\L(s,y,u_*)+g(z):\, s\in [0, T],\, y, z\in B_{\delta_*+T|b(u_*)|T}^n(x_*)\}.\]
\end{remark}
\subsection{Nice minimizing pairs}
As a consequence of Theorem~\ref{thm:main2}, the existence of uniformly equi-Lipschitz minimizing sequences under a growth assumption of type (H).
\begin{corollary}[\textbf{Equi Lipschitz/bounded  minimizing pairs}] \label{coro:coro1NEW}
Assume that  $\Lambda$  satisfies Assumption  {\rm (A)} and  Condition {\rm (S)}. Unless $\L$ is real valued assume, moreover, hypotheses h$_{1}$), h$_{2}$) of Theorem~\ref{thm:main2}.
\begin{enumerate}
\item Let $t\in [0,T[$, $x\in \R^n$ and suppose  that $\L$ satisfies Condition $\big($H$^t_{J_t(y,u)}\big)$, for a suitable admissible pair $(y,u)$ for {\rm (P}$_{t,x}${\rm )}.
    Then there are a minimizing sequence  $(\overline y_j, \overline u_j)_j$ for {\rm (P}$_{t,x}${\rm )}, a constant $K_{t,x}$ such that
$\|\overline y_j\|_{\infty}\le K_{t,x}, \|\overline u_j\|_{\infty}\le K_{t,x}$
and each $\overline y_j$ is Lipschitz of rank $K_{t,x}$.
\item (\textbf{Uniformity w.r.t.  $t, x, j$}) Let $0\le \delta<T$, $x_*\in\mathbb R^n$ and $\delta_*\ge 0$. Assume that, for some $B\ge 0$,  and any  $t\in [0, \delta], x\in B_{\delta_*}^n(x_*)$, there is an admissible pair $(y,u)$ for {\rm (P}$_{t,x}${\rm )} satisfying $J_t(y,u)\le B$. Furthermore, suppose that $\L$   satisfies Condition {\rm (H}$^{\delta}_{B}${\rm )}.
Then the constants $K_{t,x}$ in Claim (1) may be chosen  to be uniformly bounded above with respect to $t\in [0, \delta], x\in B_{\delta_*}^n(x_*)$.
\end{enumerate}
\end{corollary}
\begin{proof} 1. Consider any minimizing sequence  $(y_j,u_j)_j$ for  {\rm (P}$_{t,x}${\rm )} with $J_t(y_j, u_j)\le J_t(y,u)$ for every $j\in \mathbb N$.
 The application of Claim (1) of Theorem~\ref{thm:main2} with $\mathcal A=\{(y_j, u_j):\, j\in  \mathbb N\}$, $B=J_t(y,u)$, $\delta=t$, $\delta_*=0$ and $x=x_*$ yields the claim.\\
2. Let $t\in [0, \delta], x\in B_{\delta_*}^n(x_*)$ and consider any minimizing sequence $(y^{t,x}_j,u^{t,x}_j)$ for  {\rm (P}$_{t,x}${\rm )}; we may assume without restriction that \[\forall j\in\mathbb N\quad J_t\big(y^{t,x}_j,u^{t,x}_j\big)\le B.\]
The application of Claim  (1) of Theorem~\ref{thm:main2}  with $\mathcal A=\{(y^{t,x}_j, u^{t,x}_j):\, j\in\mathbb N,\, t\in [0, \delta],\,x\in B_{\delta_*}^n(x_*)\}$  allows to conclude.
\end{proof}
\begin{remark}
The construction of a equi-Lipschitz minimizing sequence was considered under  Condition {\rm (G)} in \cite{CF1} for {\em continuous} and {\em autonomous} Lagrangians of the calculus of variations, under  prescribed boundary data and conditions, assuming either convexity or differentiability of $\L(y,y')$ with respect to the ``velocity'' variable $y'$.
Corollary~\ref{coro:coro1NEW} extends \cite[Theorems 1,2,3,4]{CF1} in several directions.
\end{remark}
\subsection{Avoidance of the Lavrentiev phenomenon}
Another  consequence of Theorem~\ref{thm:main2} is the avoidance of the Lavrentiev phenomenon under  Condition (S) and the weakest growth Condition of type (M). Differently from other results in the literature (see \cite{Tonon, Zas}), we do not assume neither continuity of $\L$,  nor the  local Lipschitz continuity of $\L(s,\cdot, u)$, nor global convexity on the control variable and we do not make use of the Maximum Principle. With respect to Corollary~\ref{coro:coro1NEW} we loose  equi--boundedness of the minimizing sequences  of controls; nevertheless in Claim (2) we still keep some uniformity with respect to the initial points and state, for any given index $j$ of the sequence.
\begin{corollary}[\textbf{Non-occurence of the Lavrentiev phenomenon}] \label{coro:coro2NEW}
Assume that  $\Lambda$  satisfies  Assumption  {\rm (A)} and  Condition {\rm (S)}. Unless $\L$ is real valued assume, moreover,  hypotheses h$_{1}$), h$_{2}$) of Theorem~\ref{thm:main2}.
\begin{enumerate}
\item Let $t\in [0,T[$, $x\in \R^n$ and suppose that $\L$ satisfies Condition $\big($M$^t_{J_t(y,u)}\big)$, for a suitable admissible pair $(y,u)$ for {\rm (P}$_{t,x}${\rm )}.
    Then there is a minimizing sequence $(y^{t,x}_j, u^{t,x}_j)_j$ for {\rm (P}$_{t,x}${\rm )} where, for each $j\in\mathbb N$,  $u_j$  is bounded and  $y_j$  is Lipschitz.
\item (\textbf{Uniformity w.r.t.  $t,x$}) Let $0\le\delta<T$, $x_*\in\mathbb R^n$ and $\delta_*\ge 0$. Assume that there is $B\ge 0$ with $\inf{\rm (P}_{t,x} {\rm )}< B$ for any  $t\in [0, \delta], x\in B_{\delta_*}^n(x_*)$.
Furthermore, suppose that $\L$ is  satisfies Condition {\rm (M}$^{\delta}_{B}${\rm )}. Then in Claim (1), one  may  choose the minimizing sequences in such a way that, for all $j\in\mathbb N$, there is a suitable constant $K_j$ such that $\|u^{t,x}_j\|_{\infty}\le K_j$ and the rank of $y^{t,x}_j$ is less than $K_j$, as  $t$ varies in $[0, \delta]$ and $x$ varies in $B_{\delta_*}^n(x_*)$.
\end{enumerate}
\end{corollary}
\begin{proof}
 The proof follows the lines of that of Corollary~\ref{coro:coro1NEW} making use of Claim (2), instead of Claim (1), of Theorem~\ref{thm:main2}.\\
1. Consider any minimizing sequence  $(y_j,u_j)_j$ for  {\rm (P}$_{t,x}${\rm )} with $J_t(y_j, u_j)\le J_t(y,u)$ for every $j\in \mathbb N$.
 The application of Claim (2) of Theorem~\ref{thm:main2} with $\mathcal A=\{(y_j, u_j):\, j\in  \mathbb N\}$, $\delta=0, B=J_t(y,u)$,  $\delta_*=0$ and $x=x_*$yields the claim.\\
2. Let $t\in [0, \delta], x\in B_{\delta_*}^n(x_*)$ and consider any minimizing sequence $(y^{t,x}_j,u^{t,x}_j)$ for  {\rm (P}$_{t,x}${\rm )}; we may assume without restriction that \[\forall j\in\mathbb N\quad J_t\big(y^{t,x}_j,u^{t,x}_j\big)\le \inf{\rm (P}_{t,x} {\rm )}+\eta_j\le  B,\]
for a suitable $\eta_j\ge 0$ with $\displaystyle\lim_{j\to +\infty}\eta_j=0$.
Fix $j\in\mathbb N$. The application of Claim (2) Theorem~\ref{thm:main2}  with \[\mathcal A:=\mathcal A_j:=\{(y^{t,x}_j, u^{t,x}_j):\,  t\in [0, \delta],\,x\in B_{\delta_*}^n(x_*)\},\quad \eta:=\eta_j\]  allows to conclude.
\end{proof}
\begin{remark} Condition {\rm (S)} plays an essential role here. In the framework of the calculus of variations the celebrated example by
 Ball and Mizel in \cite{BM} exhibits a Lagrangian $\L(s,y,y')$ that is a polynomial, superlinear and convex in $y'$, for which the Lavrentiev phenomenon occurs with some suitable boundary data. As shown in \cite{ASC}, autonomous problems of the calculus of variations (therefore with the dynamics $y'=u$) do not face the Lavrentiev phenomenon, no matter the growth of the Lagrangian is.
\end{remark}
Taking into account Proposition~\ref{prop:whenM}, we deduce as a particular case of Corollary~\ref{coro:coro2NEW}, the nonoccurrence of the Lavrentiev phenomenon for a wide  class of functionals with a real valued Lagrangian. We stress the fact that apart radial convexity in the control variable, no more regularity than measurability is required on the state and control variable.
\begin{corollary}[\textbf{Nonoccurrence of the Lavrentiev phenomenon for real valued Lagrangians that are radially convex in the control variable}]\label{coro:realLav} Assume that
\begin{itemize}
\item $\L$ is real valued and \textbf{bounded on the bounded sets };
\item $\L$ satisfies  Condition {\rm (S)};
\item $\L$  is \textbf{radially convex in the control} variable, i.e., for a.e. $s\in [0, T]$ and every $y\in\R^n,  u\in \mathcal U$, the map $0<r\mapsto\L(s,y,r{u})$ is convex.
    \end{itemize}
Then the Lavrentiev phenomenon for {\rm (P$_{t,x}$)} does not occur.
\end{corollary}
\begin{remark}Example~\ref{ex:nononlyconvex} shows that there are Lagrangians that are not radially convex in the control variable for which the Lavrentiev phenomenon does not occur.
\end{remark}


\section{Existence and regularity of optimal pairs}
The following Lipschitz regularity result somewhat extends and partly overlaps those  formulated in \cite{BM2, BM3} in the extended valued case. The lower semicontinuity assumption of $u\mapsto \L(s,y,u)$ in  \cite[Theorem 4.2]{BM3} is replaced here by hypotheses h$_{1}$) and h$_{2}$). Also, as explained above, the growth condition {\rm (H$_B^{\delta}$)} considered here is somewhat less restrictive than the one considered previously. At the same time we require here the structure assumption (A), not needed in \cite{BM2, BM3}.
\begin{corollary}[\textbf{Lipschitz regularity}]\label{coro:regular666}
Suppose that  $\Lambda$  satisfies Assumption  {\rm (A)} and  Condition {\rm (S)}.
Unless $\L$ is real valued assume hypotheses h$_{1}$), h$_{2}$) of Theorem~\ref{thm:main2}.\\
 Let $t\in [0,T[$, $x\in \R^n$ and $(y_*,u_*)$ be an admissible pair for {\rm (P}$_{t,x}${\rm )} and suppose that $\L$ satisfies Condition $\big($H$_{J_t(y_*,u_*)}^t\big)$. Then $u_*$ is bounded and $y_*$ is Lipschitz.
\end{corollary}
\begin{proof}
If $u_*$ is not bounded then Claim (1)(c) of Theorem~\ref{thm:main2} provides the existence of an admissible pair $(\overline y, \overline u)$ with $J_t(\overline y, \overline u)<J_t(y_*, u_*)$, a contradiction. Therefore $u_*$ is bounded. The controlled differential equation $y_*'=b(y_*)u_*$ implies the Lipschitzianity of $y_*$.
\end{proof}
The following existence result for the optimal control problem {\rm (P}$_{t,x}${\rm )}, follows easily from the previous claims. In the case of the real valued case of the calculus of variations it gives back \cite[Theorem 3]{Clarke1993}, though the Lagrangians with no linear growth from below escape from our method. In the extended valued case of the calculus of variations our Condition (H) brings in some new cases with respect to the one considered in \cite{Clarke1993} (see Remark~\ref{rem:newHwhy}); at the same time the uniform limit hypothesis h$_{2}$) at the boundary of $\Dom(\L)$ is more restrictive than the lower semi-continuity alone required in \cite[Theorem 3]{Clarke1993}.
In the framework of optimal control problems, the existence question under the same slow growth condition was considered for autonomous Lagrangians (where $\Phi(B)=0$) in
\cite{FrankTrans} with a more general type of differential controlled  equation of the form $y'=f(y,u)$, assuming some extra regularity assumptions, e.g., {\em local Lipschitz continuity} on $\L(y,u)$ and $f(y,u)$ with respect to the $y$ variable, not required here.
\begin{corollary}[\textbf{Existence and regularity of a solution to {\rm (P}$_{t,x}${\rm )}}] \label{coro:coro2}
Suppose that  $\Lambda$  satisfies Assumption  {\rm (A)} and  Condition {\rm (S)}.
Unless $\L$ is real valued assume hypotheses h$_{1}$), h$_{2}$) of Theorem~\ref{thm:main2}.
Let $t\in [0,T[$, $x\in \R^n$. Suppose that $\L$ satisfies Condition $\big($H$_{J_t(y,u)}^t\big)$ for some admissible pair $(y,u)$ for {\rm (P}$_{t,x}${\rm )}.
Moreover, suppose the validity of the following structure conditions:
\begin{itemize}
\item For a.e. $s\in [t,T]$ the function $(y,u)\mapsto \Lambda(s,y, u)$ is \textbf{l.s.c.};
\item For a.e. $s\in [t,T]$ and every $y\in\R^n$ the function
 $u\mapsto \Lambda(s,y,u)$ is \textbf{convex};
 \item The  cost function $g$ is \textbf{l.s.c.} and $b$ is \textbf{continuous};
\item The set $\mathcal U\subseteq\R^m$ is \textbf{closed and convex} and $\mathcal S\subseteq\R^n$ is closed.
 \end{itemize}
 Then Problem {\rm (P}$_{t,x}${\rm )} admits a solution $(y_*,u_*)$, with $y_*$ Lipschitz  and $u_*$ bounded.
\end{corollary}
\begin{proof}
Let $(y_j,u_j)_j$ be a minimizing sequence for {\rm (P}$_{t,x}${\rm )}. We may assume, from  Claim (1) of Corollary~\ref{coro:coro1NEW}, that  $y_j$ are equi-Lipschitz, equi-bounded and that the controls $u_j$ are uniformly bounded. Ascoli's theorem implies that, modulo a subsequence,  $y_j$ converges uniformly to a Lipschitz function $y_*$ on $[t,T]$; the closure of $\mathcal S$ implies that $y_*([t,T])\subseteq \mathcal S$.
By the reflexivity of $L^2[t,T]$ we may also assume that $u_j$ converges weakly in $L^2$ to a function $u_*$: Mazur's lemma shows that $u_*$ is bounded and that, due to the closure and convexity of the sets $\mathcal U$, that $u_*(s)\in \mathcal U$ for a.e. $s$. We may then invoke a standard integral semicontinuity theorem (see, for instance, \cite[Theorem 6.38]{ClarkeBook}) to deduce that
\[J_t(y_*,u_*)\le\liminf_{j\to +\infty}J_t(y_j,u_j).\]
\noindent
There remains only to verify that $y_*$ is the state trajectory  corresponding to $u_*$, which is a standard matter. We proceed by following, for instance, the proof of \cite[Theorem 23.11]{ClarkeBook}. It is enough to show that, for any measurable subset $A$ of $[t,T]$, we have
\[\int_Ay_*'(s)-b(y_*(s))u_*(s)\,ds=0.\]
The equality holds when $y_*$ and $u_*$ are replaced by $y_j$ and $u_j$, respectively. To obtain the desired conclusion, it suffices to justify passing to the limit as $j\to +\infty$. By weak convergence, and by dominated convergence theorem, we have
\[\int_Ay_j'(s)\,ds\to\int_Ay_*'(s)\,ds \text { as }j\to +\infty.\]
We also know that, as $j\to +\infty$,
\[\int_Ab(y_*(s))u_j(s)\,ds\to \int_Ab(y_*(s))u_*(s)\,ds,\]
since $b(y_*(s))$ is bounded.
By H\"older's inequality we have that
\begin{align}\left|\int_A\big(b(y_j(s))-b(y_*(s))\big)u_j(s)\,ds\right|\le\phantom{AAAAAAAAAAAAAAAAAAA}\\ \le
\left(\int_A|b(y_j(s))-b(y_*(s))|^2\,ds \right)^{1/2}\left(\int_Au_j^2(s)\,ds\right)^{1/2}\to 0\end{align}
as $j\to +\infty$. Indeed the first factor tends to 0 by dominated convergence, and the second is uniformly bounded since the sequence $(u_j)_j$ is bounded in $L^2[t,T]$.
Therefore
\[\int_Ab(y_j(s))u_j(s)\,ds= \int_Ab(y_j(s))-b(y_*(s))u_j(s)\,ds+
 \int_Ab(y_*(s))u_j(s)\,ds \]
 tends to
$\displaystyle \int_Ab(y_*(s))u_*(s)\,ds$ as $j\to +\infty$. The result follows.\\
\end{proof}
\begin{remark} The assumptions of Corollary~\ref{coro:coro2} imply that the pair $(y_*, u_*)$ satisfies the Du Bois-Reymond -- Erdmann variational inequality established in \cite{BM3}: there is a real valued absolutely continuous function $p(s)$, whose derivative belongs a.e. to Clarke's partial subgradient $\partial^C_s\L(s, y_*(s), u_*(s))$ of $\tau\mapsto\L(\tau,y_*(s),u_*(s))$ at $\tau=s$, satisfying, for a.e. $s\in [0,T]$,
\[\L(s, y_*,ru_*)-\L(s,y_*,u_*)\ge (\L(s, y_*,u_*)-p(s))(r-1)\quad \forall r>0,\]
or, equivalently,
\[\L(s, y_*,u_*)-p(s)\in\partial_r\L(s,y_*,ru_*)_{r=1}.\]
\end{remark}
\section{Examples}\label{sect:examples}
We consider here some examples related to the growth conditions introduced above and some Lagrangians  to which our results may be applied.
\subsection{Growth Conditions}
\begin{example}[Radial functions of one variable] Consider first the case of a function of a positive real variable. Let $L:\R\to\R\cup\{+\infty\}$ be convex with $\Dom(L)=[0, +\infty[$, assume that $L$ is not definitely affine, namely that for every $c\ge 0$  there is $\overline \nu>0$ such that $\max\partial L (c)<\min \partial L (\overline\nu)$.
Then $ L$ satisfies Condition (H$_B^{\delta}$), no matter what are $B, \delta\in [0, T[$.
Indeed
 $\partial L(ru)_{r=1}=u\,\partial L(u)$ for any $u\ge 0$. Let $Q(u)\in\partial L(ru)_{r=1}$, $u\ge 0$ and set $P(u):= L(u)-Q(u)$. The monotonicity of the subdifferential of $ L$ implies that
 \[v\ge u\ge 0\Rightarrow P(v)\le P(u),\]
 the inequality being strict if $\max \partial L(u)<\min \partial L(v)$.
 Fix $c>0$. The fact that $ L$ is not definitely affine implies that there exists $\overline \nu>0$ such that
 $\max \partial L(c)<\min \partial L(\overline\nu)$.
Therefore, we obtain
 \[
\sup_{\substack{|v|\ge \overline\nu\\  L(v)<+\infty}}
P(v)\le P(\overline\nu)<P(c)\le\inf_{\substack{|v|<c\\  L(v)<+\infty}}
%
P(v).\]
More generally if $ L:\R\to \R$ is radially convex and   $c>0$, \eqref{tag:H} is fulfilled whenever (see Figure~\ref{fig:Hnewreal}) \begin{equation}\label{tag:Hconvex}\max\{P(\overline\nu), P(-\overline\nu)\}<\min \{P(c), P(-c)\}.\end{equation}
\begin{figure}[htbp]
\begin{center}
\includegraphics[width=0.6\textwidth]{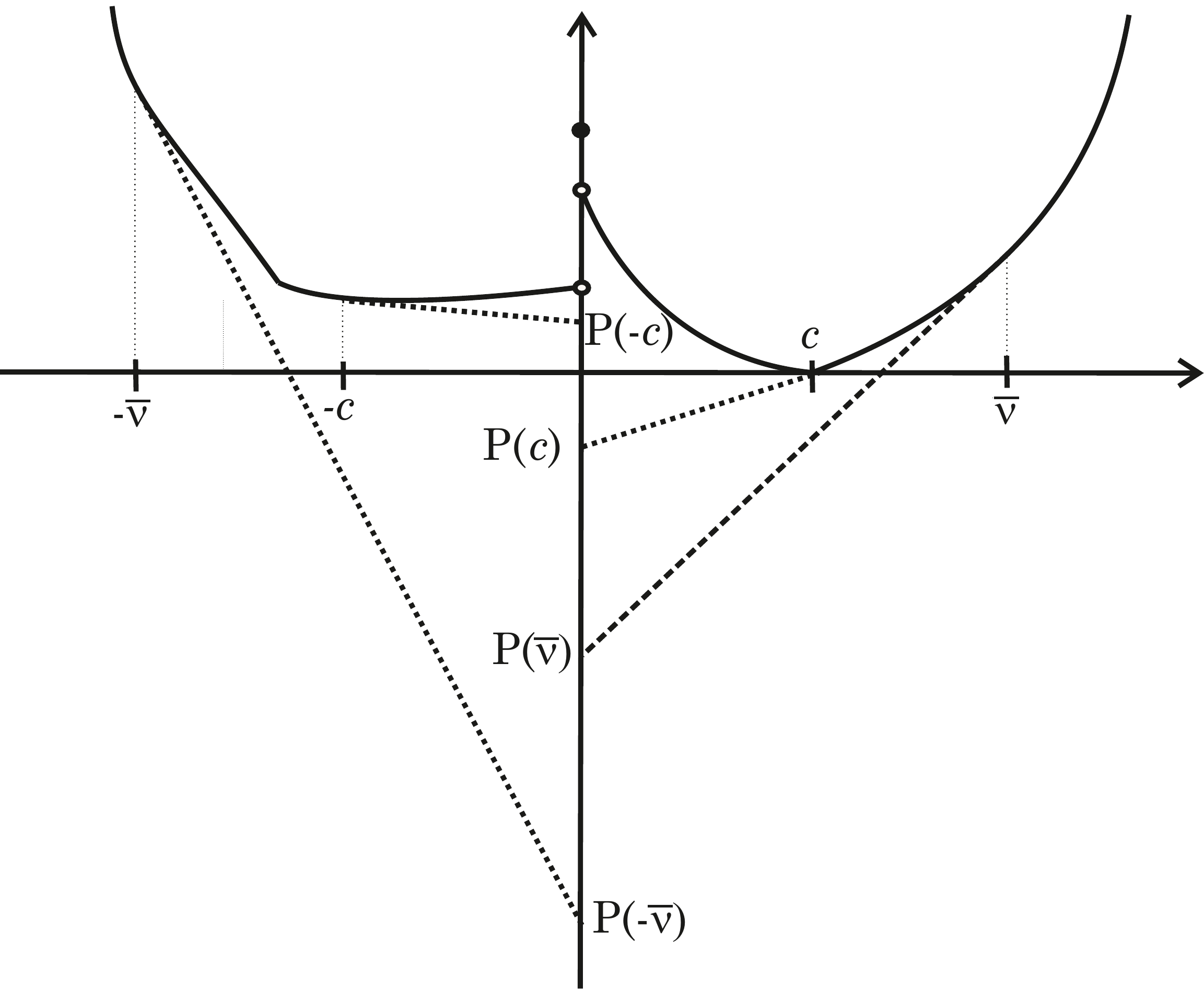}\
\end{center}
\caption{Condition \eqref{tag:Hconvex} in the radial convex case.}
\label{fig:Hnewreal}
\end{figure}
Notice that Condition (M$_B^{\delta}$) is fulfilled for every $B\ge 0, \delta\in [0, T[$ (see Proposition~\ref{prop:whenM}).
\end{example}
\begin{example}\label{ex:Hnew} The example illustrates the role of the condition $\dist((s,y,u),\partial\Dom(\L))\ge \rho$ in \eqref{tag:H} -- \eqref{tag:Hdicff}.
Consider the  function
\[\L(s,y,u)=L(u)=\begin{cases}+\infty\, &u\le -1\\
\dfrac1{1-u^2}\, &-1<u\le 0,\\
u^2+1\, &u>0.\end{cases}\]
If $u>-1$,
\[L(u)-uL'(u)=\begin{cases}\dfrac{1-3 u^2}{\left(u^2-1\right)^2}\, &-1<u<0\\
-u^2\,&u\ge 0.\end{cases}\]
In particular
\[
\lim_{u\to +\infty}L(u)-uL'(u)=-\infty,\]
so that $\L$ satisfies Condition (G).\\
Consider the problem of minimizing $I(z):=\displaystyle\int_0^1\L(z'(s))\,ds$ among the absolutely functions $z:[0,1]\to\R$ satisfying $z(0)=0, z(1)=1$: we thus consider here problem (P$_{t,x}$) with $t=0$ and $x=0$.
An admissible pair is $y(s):=s$; we set
\[B:=I(y)=\displaystyle\int_0^1(1^2+1)\,ds=2.\]
Since $L(u)\ge 2|u|$ for all $u\in\R$ $d>0$,
from  Definition~\ref{def:cphi} we get $c_0(B)=1$.
Now, since $\displaystyle\lim_{u\to (-1)^+}L(u)-uL'(u)=-\infty$, then for any $c>1$ one gets
\[\inf_{\substack{|u|<c\\ u>-1}}\left\{L(u)-uL'(u)\right\}=-\infty,\]
so that $\L$ does not satisfy \cite[Condition (H2)]{Clarke1993}.
Instead, for
any $0<c$ and $0<\rho<1$,
\[\inf_{\substack{|u|<c\\ u>-1+\rho}}\left\{L(u)-uL'(u)\right\}\ge\min\left\{-c^2, \dfrac{1-3 (-1+\rho)^2}{\left((-1+\rho)^2-1\right)^2}\right\}.\]
Therefore $\L$ fulfills Condition (H$_2^{\delta}$) for all $ \delta\in [0,1[$.
\end{example}
\begin{example}[\textbf{Lack of boundedness on bounded sets: a case where {\rm (G)} holds but (H$_B^{\delta}$) does not}]\label{ex:GnotH}
We show here that the  local boundedness assumption in Proposition~\ref{prop:GimpliesH} is crucial in order to obtain the conclusion.
Let, for any $s\in [0,1], y\in\R^2$, $u=(u_1, u_2)\in\R^2$
\[\L(s,y,u)=L(u):=\begin{cases}
(u_1^2+u_2^2)\dfrac{u_2}{u_1}&\text{ if }0<u_1\le u_2,\\
u_1^2+u_2^2&\text{ otherwise}.
\end{cases}\]
\begin{figure}[htbp]
\begin{center}
\includegraphics[width=0.4\textwidth]{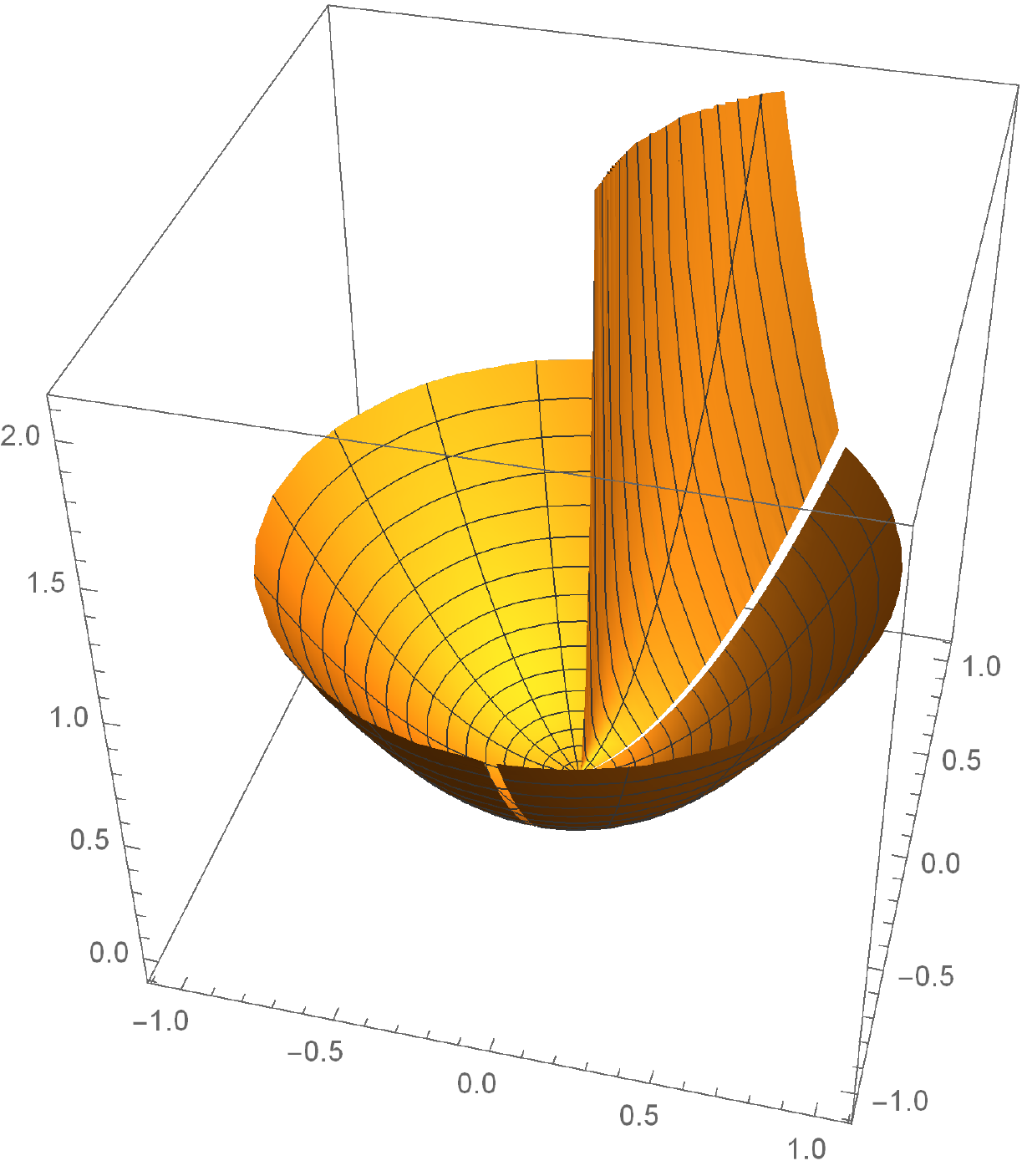}
\includegraphics[width=0.4\textwidth]{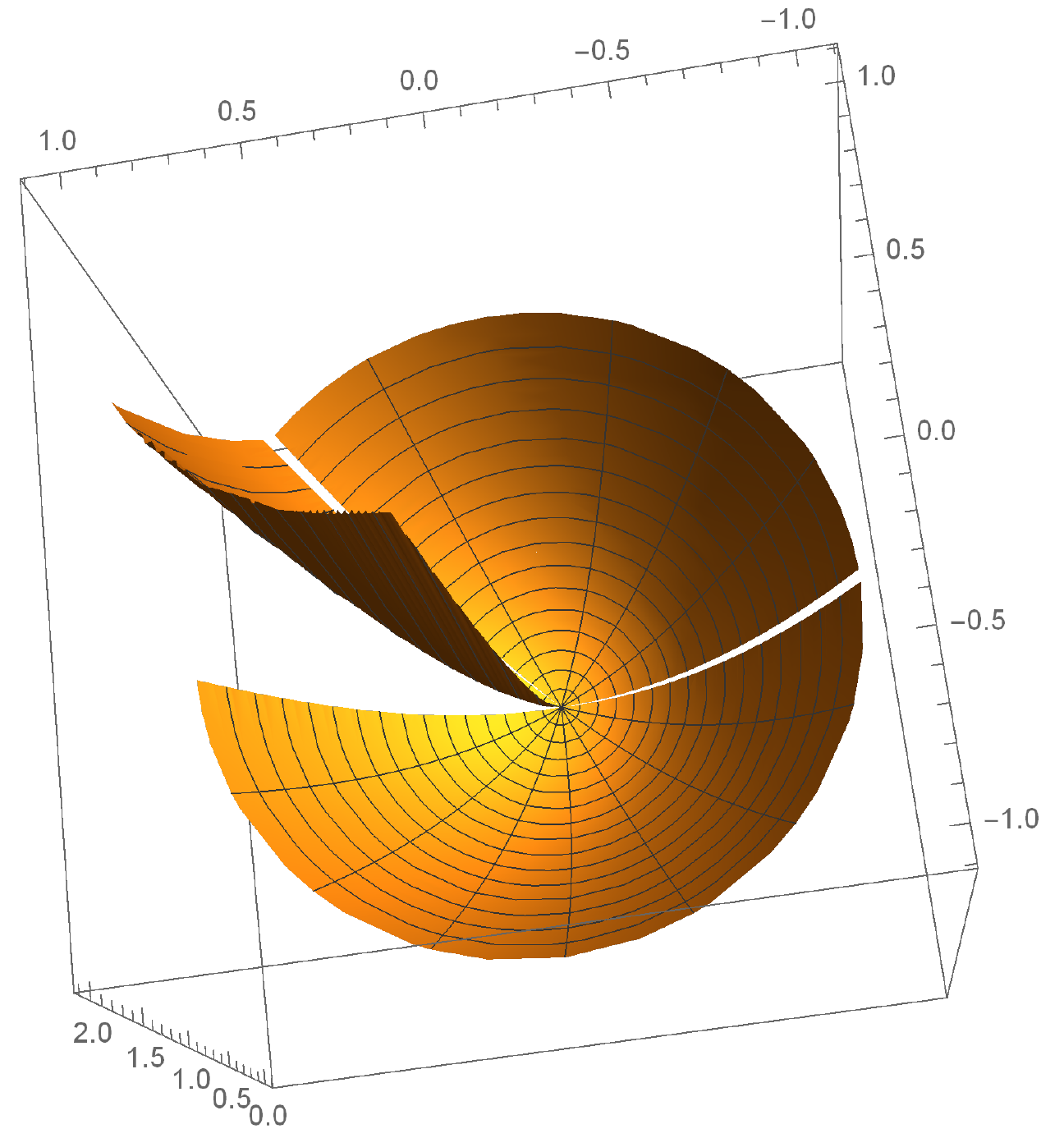}
\end{center}
\caption{The function $L$ in Example~\ref{ex:GnotH}.}
\label{fig:figurexxyy}
\end{figure}
Then
\[ L(u)-D_u L(u)=\begin{cases}-L(u)&\text{ if }0<u_1\le u_2,\\
-|u|^2&\text{ otherwise}.
\end{cases}\]
\begin{itemize}
\item $\L$ satisfies the growth Condition  {\rm (G)}. Indeed
for each $u\in\R^2$ we have
\[ L(u)-D_u L(u)\le -|u|^2\to -\infty\]
uniformly as $|u|\to +\infty$.
\item
Condition (H$_B^{\delta}$) does not hold, no matter  what are the choices of $\delta\in [0,1[, B\ge 0$.
Indeed, fix $c>0$,  $\alpha\in [\pi/4, \pi/2[$ and $\overline u:=\dfrac c2 (\cos\alpha, \sin\alpha)$.  Then
\[ L(\overline u)-D_{\overline u} L(\overline u)=-\dfrac{c^2}4\tan(\alpha)\to -\infty\text{ as } \alpha\to\left(\dfrac{\pi}2\right)^-.\]
%
Therefore
\[\inf_{\substack{|u|< {c}}}\{ L(u)-D_u L(u)\}=-\infty,\]
proving the claim (see  (3) of Remark~\ref{rem:HH}).
\end{itemize}
\end{example}
\begin{example}[\textbf{Violations of (M$_B^{\delta}$)}] Proposition~\ref{prop:whenM} exhibits
a wide variety of Lagrangians that satisfy Condition (M$_B^{\delta}$). Here are some pathological examples where its validity fails.
\begin{itemize}[leftmargin=*]
\item Let, for $s\in [0, T], y\in\R^2, u=(u_1, u_2)\in \R^2$,  \[\L(s,y, u)=L(u):=\begin{cases}\beta|u|^2-1&\text{ if } u_2=\beta u_1\text{ for some }\beta\in\R;\\
    0&\text{ if } u_1=0.
\end{cases}\]
Then, for every $c>0$,
\[\inf_{|u|<c}L(u)-D_uL(u)=-1-\beta c^2\to -\infty\text{ as }\beta\to +\infty,\]
thus violating  {\em i)} in \eqref{tag:HequivNEW} of Condition (M$_B^{\delta}$).
\item Let $\L(s,y,u)=L(u):=\sqrt{|u|}$ for every $s\in [0, T],y, u\in\R$. Then, for every $\nu>0$,
\[\sup_{|u|\le\nu}L(u)-u L'(u)=\sup_{|u|\le\nu}\dfrac{\sqrt{|u|}}{2}=+\infty,\]
thus violating {\em ii)} in \eqref{tag:HequivNEW} of Condition (M$_B^{\delta}$).
\end{itemize}
\end{example}
\subsection{Application of the main results: real valued case}
The next Example~\ref{ex:discont}, a discontinuous version of \cite[Example 4.3]{Clarke1993}, exhibits a {\em real valued} Lagrangian of the calculus of variations that satisfies the assumptions of Theorem~\ref{thm:main2}, to which the previous results of the literature (e.g., \cite{CFM}, \cite{CF1}, \cite{Clarke1993}) do not apply.
\begin{example}
\label{ex:discont}
Let $\phi:[0,1]\to\R$  be a $C^1$ function whose minimum is $m_{\phi}>0$. For $y,u\in\R^2$ let
\[\L(s,y,u):=\phi(s)a(y)\sqrt{1+|u|^2}
\]
where $a$ is any {\rm measurable} and  bounded function with $\inf a=1, a\le 2$: we stress the fact that $a$ might be discontinuous.
We consider the problem  (P) of minimizing
\[J_0(y,y')=\displaystyle\int_0^1\L(s, y(s), y'(s))\,ds\] among the functions $y\in W^{1,1}([0,1];\R^2)$ satisfying
\[y(0)=(0,0), y(1)=\zeta=(\zeta_1,\zeta_2), \zeta_1>0, \zeta_2>0.\]
\begin{itemize}
\item We have $\L(s,y,u)\ge m_{\phi}|u|$ for all $s\in [0,1]$ and $y,u$ in $\R^2$.
\item {\em $\L$ satisfies Condition {\rm {\rm (S)}} with $\kappa=\dfrac{\|\phi'\|_{\infty}}{m_{\phi}}$, $A=\gamma=0$.}
Indeed, if $s,s_1,s_2\in [0,1]$ and $y,u\in\R^2$ then
\begin{equation}\begin{aligned}\label{tag:condS2}|\L(s_2,y,u)-\L(s_1,y,u)|&\le \|\phi'\|_{\infty}a(y)\sqrt{1+|u|^2}|s_2-s_1|\\
&\le \dfrac{\|\phi'\|_{\infty}}{m_{\phi}}\L(s,y,u)|s_2-s_1|.\end{aligned}\end{equation}
\item $0<r\mapsto \L(s,y,ru)=\phi(s)a(y)\sqrt{1+r^2|u|^2}$ is convex for all $s\in [0,1]$ and $y,u\in\R^2$.
\item {\em For a suitable choice of $\phi$, $\L$ satisfies assumption {\rm (H$_B^0$)}, }
where
\[B:=J_0(y_*,  y_*')=\sqrt{1+|\zeta|^2}\int_0^1\phi(s)a(s\zeta )\,ds, \quad y_*(t):=t\zeta, \, t\in [0,1].\]
From Definition~\ref{def:cphi}
we get
\[\Phi(B)=\dfrac{\|\phi'\|_{\infty}}{m_{\phi}}B,\quad c_0(B)=\dfrac{B}{m_{\phi}}.\]
For all $(s,y,u)\in [0,1]\times\R^2\times\R^2$,
\[\begin{aligned}\L(s,y,u)-u\cdot \nabla_u\L(s,y,u)&=
\phi(s)a(y)\sqrt{1+|u|^2}-\phi(s)a(y)\dfrac{|u|^2}{\sqrt{1+|u|^2}}\\
&=\dfrac{\phi(s)a(y)}{\sqrt{1+|u|^2}}\end{aligned}
.\]
Therefore, for any $c\in\R$ and $K\ge 0$ we have
\[\inf_{\substack{|u|< c, u\in \R^2\\  s\in [0,1], |y|\le K}}
\{\L(s,y,u)-u\cdot \nabla_{u}\L(s,y,u)\}=\dfrac{m_{\phi}}{\sqrt{1+c^2}},\]
and for $\nu\in\R$,
\[\sup_{\substack{|u|\ge \nu, u\in \R^2\\  s\in [0,1], |y|\le K}}
\{\L(s,y,u)-u\cdot\nabla_u\L(s,y,u)\}+2\Phi(B)=
\dfrac{\|a\|_{\infty}\|\phi\|_{\infty}}{\sqrt{1+\nu^2}}+2\Phi(B),\]
so that
\[\lim_{\nu\to +\infty}\sup_{\substack{|u|\ge \nu, u\in \R^2\\  s\in [0,1], |y|\le K}}
\{\L(s,y,u)-u\cdot\nabla_u\L(s,y,u)\}+2\Phi(B)=2\Phi(B).\]
Therefore (H$_B^0$) is satisfied whenever, for some \begin{equation}\label{tag:c0big}c>c_0(B)=\dfrac{B}{m_{\phi}}
\end{equation}
we get
\begin{equation}\label{tag:eqc}2\Phi(B)=\dfrac{2\lambda B}{m_{\phi}}<\dfrac {m_{\phi}}{\sqrt{1+c^2}},\quad \lambda:=\|\phi'\|_{\infty}.\end{equation}
Either  $\phi'\equiv 0$, in which case any $c$ satisfies \eqref{tag:eqc},  or \eqref{tag:c0big} and \eqref{tag:eqc} are satisfied whenever
\begin{equation}\label{tag:c0small}\dfrac{B}{m_{\phi}}<c<2\sqrt{1+c^2}<\dfrac{m_{\phi}^2}{\lambda B}.\end{equation}
Now, $B$ is bounded above by $2(m_{\phi}+\|\phi'\|_{\infty})\sqrt{1+|\zeta|^2}$.
Thus, a sufficient condition for the validity of \eqref{tag:c0small} is
\[\dfrac {2\sqrt{1+|\zeta|^2}(m_{\phi}+\lambda)}{m_{\phi}}
<c<2\sqrt{1+c^2}<
\dfrac{m_{\phi}^2}{2\lambda(m_{\phi}+\lambda)
\sqrt{1+|\zeta|^2}}.\]

Since
\[\lim_{\lambda\to 0^+}\dfrac{m_{\phi}+\lambda}{m_{\phi}}
=1,\quad \lim_{\lambda\to 0^+}
\dfrac{m_{\phi}^2}{\lambda(m_{\phi}+\lambda)}=+\infty,
\]
the existence of $c$ is ensured at least for sufficiently small values of $\lambda=\|\phi'\|_{\infty}$ i.e., whenever $\phi$ is sufficiently close for being a constant.\\
\end{itemize}
Note, however, that $\L$ does not satisfy the assumptions required in the results of \cite{CF1} or those of \cite{Clarke1993}:
\begin{itemize}
\item  Depending on the regularity of $a$, $\L$ {\em might not be lower semicontinuous};
\item {\em $\L$ does not satisfy the growth Condition}  {\rm (G)}. Indeed, for $s\in [0,1]$ and $y\in \R^2$,
    \[\lim_{|u|\to +\infty}\L(s,y,u)-u\cdot \nabla_u\L(s,y,u)=\lim_{|u|\to +\infty}\dfrac{\phi(s)a(y)}{\sqrt{1+|u|^2}}=0.\]
\end{itemize}
\end{example}
\subsection{Applications of the main result: extended valued case}
In Example~\ref{ex:unbounded} we exhibit a nonautonomous, {\em extended valued  Lagrangian} $\L(s,y,u)$ that satisfies the assumptions of Theorem~\ref{thm:main2} without being regular in the state variable, nor convex in the control variable.
\begin{example}\label{ex:unbounded}
 For $s\in [0,1], y\in\R^2, u=(u_1, u_2)\in\R^2$ let
 \begin{equation}\label{tag:Ldefi}L(u_1, u_2):=\begin{cases}
 \dfrac1{1-u_1^2-u_2^2}& \text{ if }u_2\le |u_1|, \, u_1^2+u_2^2<1,\\
 \\
 \dfrac{u_1^2+u_2^2}{1-2|u_1|u_2}& \text{ if }u_2>|u_1|, \, u_2\le \dfrac1{2 |u_1|},\, |u_1|\le \dfrac1{\sqrt 2},\\
 \\
 +\infty&\text{ otherwise};
 \end{cases}\end{equation}
 the graph of $L$ is depicted in Figure~\ref{fig:figureL}.
 \begin{figure}[htbp]
\begin{center}
\includegraphics[width=0.5\textwidth]{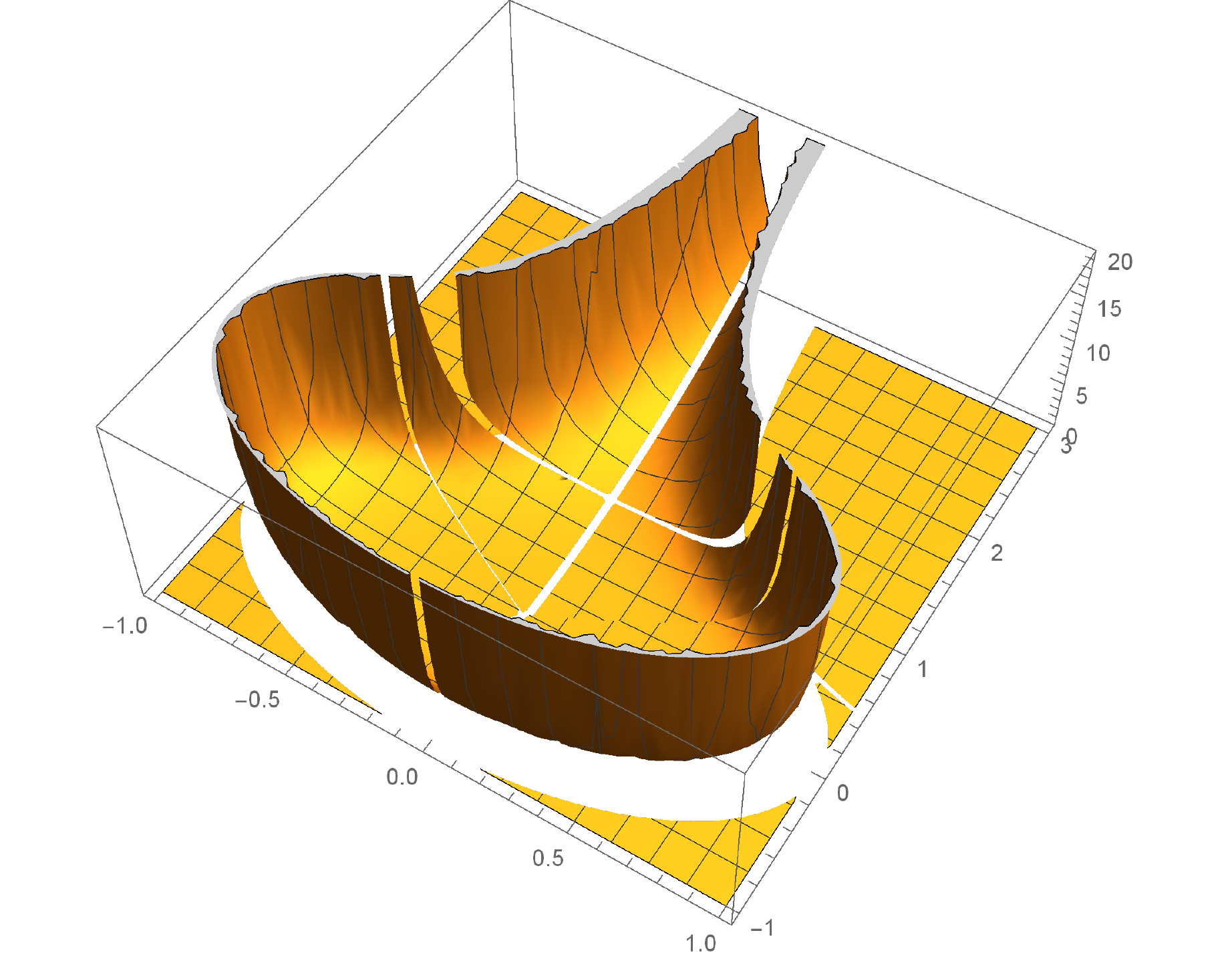}
\end{center}
\caption{The function $L(u_1, u_2)$ in Example~\ref{ex:unbounded}.}
\label{fig:figureL}
\end{figure}
The domain of $L$ is star-shaped with respect to the origin.
The function $L$ has the following properties.
\begin{itemize}[leftmargin=*]
\item {\em $L$ tends uniformly to $+\infty$ at the boundary of its effective domain.}\\
Fix $0<\varepsilon<\dfrac1{\sqrt 2}$ and consider the $\varepsilon$-neighbourhood of $\partial\Dom(L)$ given by
\[U_{\varepsilon}:=\{u\in\Dom(L):\, \dist(u, \partial \Dom(L))\le \varepsilon\}.\]
Then $U_\varepsilon\subseteq U^1_{\varepsilon}\cup U^2_{\varepsilon}\cup U^3_{\varepsilon}$, where (see Figure~\ref{fig:unboundedX})
\[U^1_{\varepsilon}:=\{(u_1,u_2)\in\Dom(L):\, 1-\varepsilon^2\le u_1^2+u_2^2<1,\, u_2\le |u_1|\},\]
\[U^2_{\varepsilon}:=\left\{(u_1,u_2)\in\Dom(L):\, \,  |u_1|\le\varepsilon,\, \dfrac1{2\varepsilon}\le u_2<\dfrac1{2|u_1|}\right\},\]
\begin{align}U^3_{\varepsilon}:=\left\{u=(u_1,u_2)\in\Dom(L):\, \,  \dfrac1{2|u_1|}> u_2\ge |u_1|\ge\varepsilon,\, \dist(u, \partial \Dom(L))<\varepsilon\right\}.\end{align}
\begin{figure}[htbp]
\begin{center}
\includegraphics[width=0.23\textwidth]{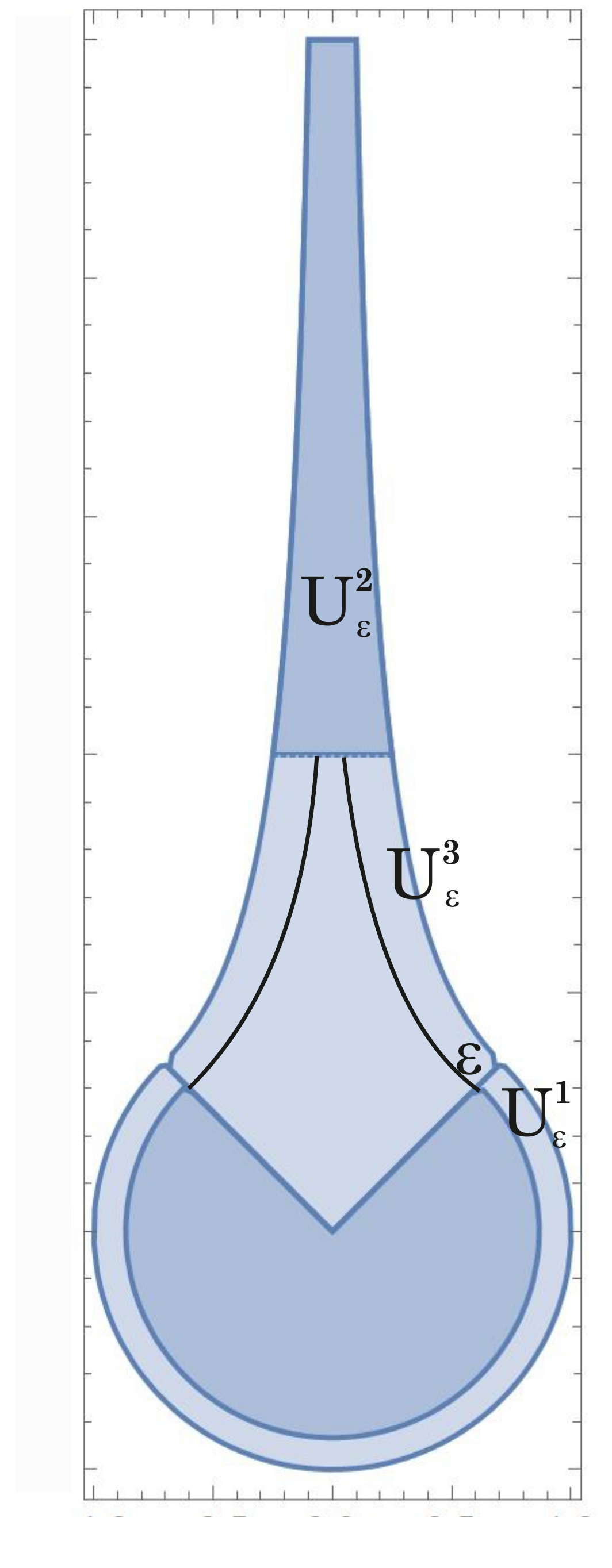}\
\end{center}
\caption{The sets $U^1_{\varepsilon}, U^2_{\varepsilon}, U^3_{\varepsilon}$ in Example~\ref{ex:unbounded}.}
\label{fig:unboundedX}
\end{figure}

On each of the above sets we obtain the following estimates from below:
\begin{itemize}
\item On $U^1_{\varepsilon}$ we have \begin{equation}\label{tag:A2}\forall u\in U^1_{\varepsilon}\quad L(u)\ge\dfrac{1}{\varepsilon^2};\end{equation}
\item On $U^2_{\varepsilon}$, $u_1^2+u_2^2\ge u_2^2\ge \dfrac{1}{4\varepsilon^2}$ and $1-2|u_1|u_2<1$ so that \begin{equation}\label{tag:B2}\forall u\in U^2_{\varepsilon}\quad L(u)\ge\dfrac{1}{4\varepsilon^2};\end{equation}
\item  On $U^3_{\varepsilon}$ we have
\begin{equation}\label{tag:C2}\forall u\in U^3_{\varepsilon}\quad L(u)\ge \dfrac{1}{2\sqrt 2\varepsilon}+o_{\varepsilon\to 0^+}(1).\end{equation}
Indeed, let $(u_1,u_2)\in U^3_{\varepsilon}$ and $|a|\in \left[\varepsilon, \frac1{\sqrt 2}\right]$  be such that
\[|u_1-a|\le\varepsilon,\, \left\vert u_2-\frac1{2|a|}\right\vert\le\varepsilon.\]
We assume that $a>0$, the case $a<0$ being similar.
Then \[0<1-2|u_1|u_2\le \varepsilon\left( \dfrac1a+2a\right)-2\varepsilon^2,\] and therefore
\begin{equation}\label{tag:estim100}L(u_1,u_2)\ge  m(a,\varepsilon):= \dfrac{(a-\varepsilon)^2+\left(\dfrac1{2a}-\varepsilon\right)^2}{\varepsilon\left(\dfrac1a+2a\right)
-2\varepsilon^2}\end{equation}
If $a\le \varepsilon |\log\varepsilon|$ then
\[\begin{aligned}m(a,\varepsilon)&\ge
\dfrac{\left(\dfrac1{2\varepsilon |\log\varepsilon|}-\varepsilon\right)^2}{\varepsilon\left(\dfrac1{\varepsilon}+
2\varepsilon\right)
-2\varepsilon^2}\\&=\left(\dfrac1{2\varepsilon |\log\varepsilon|}-\varepsilon\right)^2\sim \left(\dfrac1{2\varepsilon |\log\varepsilon|}\right)^2\quad \varepsilon\to 0^+.\end{aligned}\]
\noindent
If $a\ge \varepsilon |\log\varepsilon|$ then, using the fact that $a\le \dfrac1{\sqrt 2}$ we obtain
\[\begin{aligned}m(a,\varepsilon)&=
\dfrac{\left(a^2+\dfrac1{4a^2}\right)-2\varepsilon\left(a+\dfrac1{2a}\right)+2\varepsilon^2}{\varepsilon\left(\dfrac1a+2a\right)
-2\varepsilon^2}\\
&\ge \dfrac{1-2\varepsilon\left(\varepsilon|\log \varepsilon| +\dfrac1{2\varepsilon|\log \varepsilon|}\right)
+2\varepsilon^2}{\varepsilon\left(\dfrac1a+2a\right)
-2\varepsilon^2}\\
&\ge \dfrac{1+o_{\varepsilon\to 0^+}(1)}{2\sqrt 2\varepsilon
-2\varepsilon^2}\sim\dfrac{1}{2\sqrt 2\varepsilon}\quad \varepsilon\to 0^+,
\end{aligned}\]
\end{itemize}
which, together with \eqref{tag:estim100}, gives \eqref{tag:C2}.

From \eqref{tag:A2}, \eqref{tag:B2} and \eqref{tag:C2} we conclude that
\[\lim_{\substack{\varepsilon\to 0^+\\ \dist(u,\partial\Dom( L))\le\varepsilon}}L(u)=+\infty,\]
which proves the claim.
\item {\em $L$ is radially convex.} Indeed, Fix $u=(u_1, u_2)\in\Dom(\L)$ and $r\in ]0,1]$.
     Then
\[L(tu_1, tu_2):=\begin{cases}
 \dfrac1{1-r^2(u_1^2+u_2^2)}& \text{ if }u_2\le |u_1|, \, u_1^2+u_2^2<1,\\
 \\
 \dfrac{r^2(u_1^2+u_2^2)}{1-2r^2|u_1|u_2}& \text{ if }u_2>|u_1|, \, u_2\le \dfrac1{2 |u_1|},\, |u_1|\le \dfrac1{\sqrt 2},\\
 \\
 +\infty&\text{ otherwise}.
 \end{cases}\]
We recognize easily that $r\in ]0,1]\mapsto L(ru_1, ru_2)$ is convex.
\item {\em $L$ is superlinear.} Indeed if $u=(u_1, u_2)\in\Dom(\L)$ and $u_2>|u_1|$ then
    \[L(u_1,u_2)=\dfrac{u_1^2+u_2^2}{1-2|u_1|u_2}\ge |u|^2.\]
\item {\em $L$ is bounded on the bounded sets that are well-inside the domain}. Indeed the domain of $L$ is open and $L$ is bounded on the relatively compact subsets of the domain.
\item {\em $L$ is not continuous in the interior of its domain and thus $L$ is not convex.}
Indeed, for instance,
\[\lim_{\substack{(u_1,u_2)\to \frac12\big(\frac1{\sqrt{2}},\frac1{\sqrt{2}}\big)\\u_2<|u_1|}}L(u_1,u_2)
=\dfrac43,\]
whereas
\[\lim_{\substack{(u_1,u_2)\to \frac12\big(\frac1{\sqrt{2}},\frac1{\sqrt{2}}\big)\\u_2>|u_1|}}L(u_1,u_2)
=\dfrac13.\]
\end{itemize}
Now, let $\phi:[0,1]\to ]0, +\infty[$ be Lipschitz, $a:\R^2\to\R$ be {\em any} measurable function that is bounded below by a positive constant and bounded on bounded sets. Define
\[\forall s\in [0,1],\, \forall (y,u)\in\R^2\times\R^2\quad \L(s,y,u):=\phi(s)a(y)L(u),
\]
where $L$ is defined in \eqref{tag:Ldefi}.
The above discussion shows that
\begin{itemize}
\item {\em The domain of $\L$ is a product}, in the sense of Hypothesis h$_1$). Indeed, $\Dom(\L)=[0,1]\times\R^2\times\Dom(L)$;
\item {\em $\L$ tends uniformly to $+\infty$ at the boundary of the effective domain}, in the sense of Hypothesis h$_2$);
\item $0<r\mapsto \L(s,y,ru)$ {\em is convex} for every $(s,y,u)\in\Dom(\L)$;
\item {\em $\L$ is superlinear} and $(s,y,u)\in\Dom(\L)$ whenever $s\in [0,1], y\in\R^n$ and $u\in\R^m$ with $|u|\le 1/2$. Proposition~\ref{prop:SuperimpliesG} then implies that $\L$ satisfies the growth Condition {\rm (G)};
\item {\em $\L$ is bounded on the bounded sets that are well-inside the domain} $\Dom(\L)$ so that, from Proposition~\ref{prop:GimpliesH}, it satisfies Condition {\rm (H$_B^{\delta}$)} for all $\delta\in [0,1[$ and $B\ge 0$;
\item {\em $\L$ satisfies Condition {\rm {\rm (S)}}.} Indeed if $c_{\phi}$ is the Lipschitz constant of $\phi$, for all $s,s_1, s_2\in [0,1]$, $y\in\R^2$ and $u\in\Dom(L)$,
    \[|\L(s_2,y,u)-\L(s_1,y,u)|\le c_{\phi}|s_2-s_1|a(y)L(u)\le c\L(s,y,u)|s_2-s_1|,\]
    where $c:=\dfrac{c_{\phi}}{\min\phi}$;
\end{itemize}
Therefore $\L$  fulfills the assumptions of Theorem~\ref{thm:main2},  might be discontinuous in the variables $y,u$, and is nonconvex in $u$.
\end{example}
\section{Growth conditions in more depth}\label{sect:growthsdeep}
We give  here the proofs of some results formulated  in Section~\ref{sect:growths}.
\subsection{Proof of Proposition~\ref{prop:atleastlinear}}
\begin{proof}[Proof of Proposition~\ref{prop:atleastlinear}]
In the convex case let $Q(s,y,{u})\in\partial_r\L(s,y,r{u})_{r=1}$ for every ${u}$ such that $(s,y,u)\in\Dom(\L)$; in the partial differentiable case we set $Q(s,y,{u}):=D_{u}\L(s,y,u)$.
From \eqref{tag:Gequiv} we may choose ${R}>0$  in such a way that for $s\in [0,T]$, every $y$  and $\L(s,y,w)<+\infty$, $w\in\mathcal U$,
\begin{equation}\label{tag:8989} Q(s,y,w)-\L(s,y,w)\ge 1\quad\forall |w|\ge{R}.\end{equation}
Fix $u\in\mathcal U$ with $|{u}|>{R}>0$, $s\in [0,T]$  and $(s,y, u)\in\Dom(\L)$. It is enough to consider the case where $\L(s,y,{u})<+\infty$. The  assumption on $\Dom(\L)$ implies that   $\L\left(s,y,{r }\dfrac{{u}}{|{u}|}\right)<+\infty$ for all $0<r \le |u|$.
We consider separately the radial convex and the partial differentiable case.
\begin{itemize}[leftmargin=*]
\item[a)] {\em Radial convex case.}
Let $0<R<r \le |{u}|$.
  We have
\[\L\left(s,y, r \dfrac{{u}}{|{u}|}\right)-\L\left(s,y, R\dfrac{{u}}{|{u}|}\right)\ge Q\left(s,y,\rho\dfrac{{u}}{|{u}|}\right)\left(\dfrac{r }{R}-1\right)\]
whence
\begin{equation}\label{tag:ling}\begin{aligned}\L\left(s,y, r \dfrac{{u}}{|{u}|}\right)&\ge \left(Q-\Lambda\right)\left(s,y,{R}\dfrac{{u}}{|{u}|}\right)\left(\dfrac{r }{{R}}-1\right)+\Lambda \left(s,y,{R}\dfrac{{u}}{|{u}|}\right)\dfrac{r }{{R}}\\
&\ge\dfrac{r }{{R}}-1.
\end{aligned}\end{equation}
By choosing $r =|{u}|$ in \eqref{tag:ling} we obtain
\begin{equation}\label{tag:eq55}\L(s,y,{u})\ge\dfrac{|u|}{{R}}-1\quad \forall |{u}|\ge {R}.\end{equation}
\item[b)] {\em Partial differentiable case.}
Assume that $\L$ satisfies \eqref{tag:Gdiff}.
   By applying \eqref{tag:8989} with $w:=r \dfrac{u}{|u|}$ we obtain
\[\forall r \in [R, |u|]\qquad \L\Big(s, y,  r \dfrac{u}{|u|} \Big)- D_{r \frac{u}{|u|}}\L\Big(s, y,  r \dfrac{u}{|u|} \Big)\le -1,\]
or, equivalently, setting $h(r ):=\L\Big(s, y,  r \dfrac{u}{|u|} \Big)$,
\[\forall r \in [R, |u|]\qquad h(r )+1\le r  \dfrac{d}{dr }h(r )=r  \dfrac{d}{dr }(h(r )+1).\]
It follows that
$\dfrac{d}{dr }\dfrac{h(r )+1}{r }\ge 0\text{ on } [R,|u|]$, whence
\[\dfrac{\L(s, y,  u)+1}{|u|}=\dfrac{h(|u|)+1}{|u|}\ge \dfrac{h(|R|)+1}{R}\ge \dfrac1R.\]
Thus
\begin{equation}\label{tag:eq55pd}\L(s,y,{u})\ge\dfrac{|u|}{{R}}-1\quad \forall |{u}|\ge {R},\, u\in\mathcal U.\end{equation}
It follows from  \eqref{tag:eq55}, \eqref{tag:eq55pd}  and the positivity of $\L$ that
\[\L(s,y,{u})\ge\dfrac{1}{R}|{u}|-2\quad \forall{u}\in\mathcal U.\]
\end{itemize}
\end{proof}
\subsection{Proof of Proposition~\ref{prop:whenM}}
\begin{proof}[Proof of Proposition~\ref{prop:whenM}] Let $\delta\in [0,T[$ and $B\ge 0$. Fix $K\ge 0$, $c>c_{\delta}(B)$ and $\rho>0$.
Let $Q(s,y,u)\in \partial_r\L(s,y,ru)_{r=1}$ whenever $(s,y,u)\in \Dom(\L)$.
Lemma~\ref{tag:boundabovewell} shows that there is $c_1(\rho, K)\in\mathbb R$ such that
\begin{equation}\label{tag:infbound}c_1(\rho, K)\le \inf_{\substack{s\in [0,T],|y|\le K\\|u|<c, u\in \mathcal U\\ \L(s, y, u)<+\infty\\ \dist((s,y,u),\partial\Dom(\L))\ge\rho}}
\left\{\L(s, y, u)-Q(s,y,u)\right\},\end{equation}
proving the validity of {\em i)} in Condition (M$_B^{\delta}$).\\
Choose $0<{\overline\nu}< r_K$ and let $(s,y,u)\in\Dom(\L)$ with $|y|\le K$ and $|u|\ge {\overline\nu}, u\in\mathcal U$.
Then  the fact that $\Dom(\L)$ is star-shaped in the control variable implies that $\Big(s, y, {\overline\nu}\dfrac{u}{|u|}\Big)\in\Dom(\L)$; radial convexity along the direction $u$ gives
\begin{equation}
\L\Big(s, y, {\overline\nu}\dfrac{u}{|u|}\Big)-\L(s, y, u)\ge Q(s,y,u)\Big(\dfrac{{\overline\nu}}{|u|}-1\Big),
\end{equation}
from which we deduce that
\begin{equation}\label{tag:allconvM}
\L\Big(s, y, {\overline\nu}\dfrac{u}{|u|}\Big)-\dfrac{{\overline\nu}}{|u|}Q(s,y,u)\ge
\L(s, y, u)-Q(s,y,u).
\end{equation}
Now, the monotonicity of the convex subdifferential implies that
\[Q(s,y,u)\ge Q\Big(s,y,\dfrac{{\overline\nu}}{|u|}u\Big).\]
Inequality \eqref{tag:allconvM} then yields
\begin{equation}
\L\Big(s, y, {\overline\nu}\dfrac{u}{|u|}\Big)-\dfrac{{\overline\nu}}{|u|}Q\Big(s,y,{\overline\nu}\dfrac{u}{|u|}\Big)\ge
\L(s, y, u)-Q(s,y,u).
\end{equation}
Again, Lemma~\ref{tag:boundabovewell} and the assumption that $\L$ is bounded on bounded sets that are well-inside the effective domain imply  that there is a constant $c_2(K)$ satisfying
\begin{equation}\label{tag:supbound}\forall s\in [0, T], |y|\le K, u\in\mathcal U, |u|\ge {\overline\nu}\qquad  \L(s, y, u)-Q(s,y,u)\le c_2(K),\end{equation}
proving the validity of {\em ii)}  of Condition (M$^{\delta}_B$).
\end{proof}
\section{Proof of the main result}\label{sect:proofs}
This section is devoted to the proof of Theorem~\ref{thm:main2}. Many technical points derive from the fact that the Lagrangian is allowed to take the value $+\infty$. Due to its length, it may be convenient to illustrate what are the main arguments. We will often write, for the sake of clarity, $\displaystyle\int_t^T\L(\phi,y,u)\,ds$ instead of $\displaystyle\int_t^T\L(\phi(s),y(s),u(s))\,ds$.\\
\subsection{Strategy of the proof of Theorem~\ref{thm:main2}}
\begin{proof}We fix $\delta\in [0,T[$, $x_*\in\R^n$, $t\in [0, \delta], x\in B_{\delta_*}^n(x_*)$ and consider an admissible pair $(y,u)$ for ($P_{t,x}$). We  build a suitable admissible pair $(\overline y, \overline u)$ with  $J_t(\overline y, \overline u)\le J_t(y,u)$ if (H$_B^{\delta}$) holds, $J_t(\overline y, \overline u)\le J_t(y,u)+\eta$, where $\eta>0$ is arbitrary, if (M$_B^{\delta}$) holds:
\begin{itemize}[leftmargin=*]
\item Let $c$ be the constant that appears in the growth conditions.
 For $\nu\ge c$ let \[S_{\nu}:=\{s\in [t,T]:\, |u(s)|>\nu\}.\]
 Then the ``excess'' function
 \[\varepsilon_{\nu}:=\displaystyle\int_{S_{\nu}}\left(\dfrac{|y'(s)|}{\nu}-1\right)\,ds\]
 tends to 0 as $\nu\to +\infty$.
 \item  There are $\mu\in ]0,1[$, and $\rho>0$ such that, for a.e. $s$ on a  non negligible set $\Omega$ and a.e. $\tilde s\in [t, T]$,
 \[\L\left(\tilde s, y(s), \dfrac{u(s)}{\mu}\right)<+\infty,\, \dfrac{|u(s)|}{\mu}<c,\, \dist \left(\left(\tilde s, y(s), \dfrac{u(s)}{\mu}\right), \partial\Dom(\L)\right)\ge \rho.\]
 We use here the fact that $c>c_{\delta}(B)$ (see Proposition~\ref{lemma:newC}).
 In the extended valued case this is where  Hypotheses  h$_1$) and h$_2$) play a role, in which case Steps {\em v)} -- {\em vi)} of the proof are rather technical.
 \item For any $\nu>0$ big enough, we can choose $\Sigma_{\nu}\subseteq\Omega$ such that $|\Sigma_{\nu}|=\dfrac{\varepsilon_{\nu}}{1-\mu}$ and $S_{\nu}\cap \Sigma_{\nu}$ is negligible.
\item  Define an absolutely continuous, bijective, increasing  function $\varphi:[t,T]\to [t,T]$ such that
  \[\varphi(t):=t,\qquad\text{for a.e. }\tau\in [t,T]\quad \varphi'(\tau)=\begin{cases}\dfrac{|u(\tau)|}{{\nu}} &\text{ if }\tau \in S_{\nu},\\
\,\,\mu &\text{ if }\tau\in\Sigma_{\nu},\\
\,\,1 &\text{ otherwise};
\end{cases}\] let $\psi=\varphi^{-1}$ be its inverse.
\item The condition $\varphi(T)=T$ is ensured from the fact that $\varepsilon_{\nu}
    =(1-\mu)|\Sigma_{\nu}|$.
    \item Let $\overline y(s):=y(\psi(s)); \overline u(s):=\dfrac{u(\psi(s))}{\varphi'(\psi(s))}$. Then $(\overline y, \overline u)$ is admissible and $|\overline u|\le \nu$ is bounded.
\item The change of variable $\tau=\psi(s)$ gives
\[\begin{aligned}
J_t(\overline y, \overline u)&=\int_t^T\L\left(\varphi(\tau), y(\tau), \dfrac{u(\tau)}{\varphi'(\tau)}\right)\varphi'(\tau)\,d\tau+g(y(T))\\
&=\int_{S_{\nu}}*+\int_{\Sigma_{\nu}}*+\int_{[t,T]\setminus (S_{\nu}\cup \Sigma_{\nu})}\L(\varphi(\tau), y(\tau), u(\tau))\,d\tau+g(y(T)),
\end{aligned}\]
where, above,  $*$ replaces $\left(\varphi(\tau), y(\tau), \dfrac{u(\tau)}{\varphi'(\tau)}\right)\varphi'(\tau)\,d\tau$.
\item Let, in the radial convex case,  $P(s,z,v)\in \partial_{\mu}\Big(\L\Big(s, y, \dfrac{{u}}{\mu}\Big ){\mu}\Big)_{{\mu}=1}$ be as in the growth conditions or set $P(s,z,v):=\L(s,z,v)-v\cdot \nabla_v\L(s,z,v)$ in the differentiable case. There is $K=K( B, \delta_*, x_*)$ such that $\|y\|_{\infty}\le K$. We obtain
\[\int_{S_{\nu}}*\le \int_{S_{\nu}}\L(\varphi, y, u)\,d\tau+\varepsilon_{\nu}\,\Xi(\nu),\quad
\Xi(\nu):=\sup_{\substack{s\in [0,T],|z|\le K\\|v|\ge \nu, v\in \mathcal U\\ \L(s, z, v)<+\infty}}
P(s,z,v).\]
Assumption h$_1$) yields
\[\int_{\Sigma_{\nu}}*\le \int_{\Sigma_{\nu}}\L(\varphi, y, u)\,d\tau -(1-\mu)|\Sigma_{\nu}| \Upsilon,\quad \Upsilon:=\!\!\!\!\!\!\!\!\!\!\displaystyle\inf_{\substack{s\in [0,T],|z|\le K\\|v|<c, v\in \mathcal U\\ \L(s, z, v)<+\infty\\ \dist((s,z,v),\partial\Dom(\L))\ge\rho}}
\!\!\!\!\!\!\!\!\!\!P(s,z,v).\]
\item As a consequence we obtain
\begin{equation}J_t(\overline y, \overline u)\le \int_{t}^T\L(\varphi, y, u)\,d\tau + \varepsilon_{\nu}(\Xi(\nu)-\Upsilon)+ g(y(T)).\end{equation}
\item For $\nu$ large enough, $\|\varphi(\tau)-\tau\|_{\infty}\le\varepsilon_*$: Condition (S) implies
\[\int_{t}^T\L(\varphi, y, u)\,d\tau\le \int_{t}^T\L(\varphi, y, u)\,d\tau+2\varepsilon_{\nu}\Phi(B),\]
so that
\begin{equation}\label{tag:FINAL}J_t(\overline y, \overline u)\le J_t(y,u)+\varepsilon_{\nu}(\Xi(\nu)+2\Phi(B)-\Upsilon).\end{equation}
\item If Condition (H$_B^{\delta}$) holds, we choose $\nu>0$ satisfying  $\Xi(\nu)+2\Phi(B)-\Upsilon<0$. From \eqref{tag:FINAL} we deduce that $J_t(\overline y, \overline u)\le J_t(y,u)$, the inequality being strict if $\varepsilon_{\nu}>0$, i.e., $|u|>\nu$  on a non negligible set.
\item     If Condition (M$_B^{\delta}$) holds, since $\varepsilon_{\nu}\to 0$ as $\nu\to +\infty$ we choose $\nu$ in such a way that $\varepsilon_{\nu}(\Xi(\nu)+2\Phi(B)-\Upsilon)\le \eta$. From \eqref{tag:FINAL} we deduce that $J_t(\overline y, \overline u)\le J_t(y,u)+\eta$.
\end{itemize}
\item Finally, though $S_{\nu}, \Sigma_{\nu}, \varepsilon_{\nu}$ might depend on the chosen admissible pair $(y,u)$, the constant $\nu$ -- therefore the  bound of $\|\overline u\|_{\infty}$ and the Lipschitz constant of $\overline y$ -- depend in fact only on $\delta, B,   \delta_*, x_*$ (and possibly on $\eta$ if one assumes just Condition (M$_B^{\delta}$) instead of (H$_B^{\delta}$)).
\end{proof}
\subsection{Proof of Theorem~\ref{thm:main2}}
\begin{proof} Fix $\delta\in [0,T[$, $x_*\in\R^n$ and let $t\in [0, \delta], x\in B_{\delta_*}^n(x_*)$.
Fix a generic pair $(y,u)\in\mathcal A$ that is admissible for {\rm (P}$_{t,x}${\rm )}. We consider the following steps.
\begin{itemize}[leftmargin=*]
\item[{\em i)}]\emph {Let $\alpha, d$ be as in \eqref{tag:lingrowth}. Then
 \[\displaystyle\int_{t}^T|u(s)|\,ds\le (T-t)c_t(B)\le  R=R(B):=\dfrac{B+d\,T}{\alpha}.\]}
 The claim follows immediately from  Point (1) of Proposition~\ref{lemma:newC}.
\item[{\em ii)}] \emph{There is $K:=K(B,  \delta_*, x_*)$ such that $|y(s)|\le K$ for every $s\in [t,T]$.}
Indeed, for a.e.  $s\in [t,T]$,
\begin{equation}\label{tag:boundyprime}|y'(s)|\le \theta(1+|y(s)|)|u(s)|.\end{equation}
Gronwall's Lemma (see \cite[Theorem 6.41]{ClarkeBook}) and Claim {\em i)} imply that, for all $s\in [t,T]$,
\[\begin{aligned}|y(s)-x|&\le \int_t^s\exp\left(\theta\int_{\tau}^s|u(r)|\,dr\right)
\theta|u(\tau)|(|x|+
1)\,d\tau
\\
&\le \theta TRe^{R\theta}(|x|+1)\le \theta T Re^{R\theta}(|x_*|+\delta_*+1).
\end{aligned}\]
The claim follows from  the fact that $R$ depends on $B$, with \[K=\delta_*+\theta TRe^{R\theta}(|x_*|+\delta_*+1).\]
\item[{\em iii)}]{\em  Definition of $\Xi(\nu)$, $\Upsilon(\rho)$, choice of $\overline\nu, \overline\rho$.}
In the {\em radial convex case}, let $P(s,z,v)$, ${c}>c_{\delta}(B)$  be as in Condition (H$_B^{\delta}$) or  (M$_B^{\delta}$) corresponding to the value $K$ defined in Claim {\em ii)}; in the {\em partial differentiable case} we set
\[P(s,z,v):=\L(s,z,v)-D_{v}\L(s,z,v).\]
For $\rho>0$ and  $\nu>0$ we define
\[\Xi(\nu):=\sup_{\substack{s\in [0,T],|z|\le K\\|v|\ge \nu, v\in \mathcal U\\ \L(s, z, v)<+\infty}}
P(s,z,v),\quad\Upsilon(\rho):=\displaystyle\inf_{\substack{s\in [0,T],|z|\le K\\|v|<c, v\in \mathcal U\\ \L(s, z, v)<+\infty\\ \dist((s,z,v),\partial\Dom(\L))\ge\rho}}
P(s,z,v).\]
We may assume that $\Xi(\nu)>-\infty$ for all $\nu> 0$, otherwise there is $\overline\nu>0$ such that $|v|\le\overline\nu$ for all $(s,y,v)\in\Dom(\L), v\in\mathcal U$ and the result follows trivially.
In view of Remark~\ref{rem:infsup100}, there is $\overline\rho>0$ in such a way that
$\Upsilon(\rho)<+\infty$  for all $0<\rho\le \overline\rho$.
Fix $\overline\nu=\overline\nu(\delta, B, K)=\overline\nu(\delta, B, \delta_*, x_*)\ge 0$ in  such a way that:
\begin{itemize}
\item[$\diamond$] If Condition {\rm (H$_B^{\delta}$)}  holds,
\begin{equation}\label{tag:choicenu50}\forall \rho\in]0, \overline\rho[\quad\forall\nu\ge\overline\nu\qquad
-\infty<\Xi(\nu)+2\Phi(B)<\Upsilon(\rho)<+\infty;\end{equation}
notice in particular that $\Xi(\nu), \Upsilon(\rho)$ are both finite.
\item[$\diamond$] If Condition {\rm (M$_B^{\delta}$)}  holds,
\begin{equation}\label{tag:choicenu500}\forall \rho\in ]0, \overline\rho[\quad\forall\nu\ge\overline\nu\qquad
\Upsilon(\rho)\in\R,\quad \,\Xi(\nu)\in \R.\end{equation}
\end{itemize}
\item[{\em iv)}] {\em For any $\dfrac{c_{\delta}(B)}{c}<\mu<1$ let
$\Omega_{\mu}:=\left\{s\in [t,T]:\, \dfrac{|u(s)|}{\mu}<c\right\}.$
Then
\begin{equation}\label{tag:estimateomegamu}|\Omega_{\mu}|\ge \left(1-\dfrac{c_{\delta}(B)}{\mu c}\right)(T-t).\end{equation}}
Indeed,  it is enough to apply Point (2) of Proposition~\ref{lemma:newC} with $\sigma:=\mu c>c_{\delta}(B)$.
\\
\item[{\em v)}] {\em For every $\rho>0$ let
\[J_{\rho}:=\{s\in [t,T]:\,\dist((s,y(s), u(s)), \partial \Dom(\L))\ge 2\rho\}.\]
Then $\displaystyle\lim_{\rho\to 0}|J_{\rho}|=T-t$  uniformly w.r.t.  $t\in [0, \delta], x\in B_{\delta_*}^n(x_*)$; more precisely, given $\varepsilon>0$ there is $\rho_{\varepsilon}=\rho_{\varepsilon}(B)$ such that $|J_{\rho}|\ge T-t-\varepsilon$ whenever $0<\rho<\rho_{\varepsilon}$. }\\
Fix $\ell>0$;  from Hypothesis h$_{2}$), there exists $\rho_\ell>0$ satisfying
\begin{equation}\label{tag:choicedelta}\forall s\in [t,T],\, \forall v\in\R^m\,\,  \dist ((s,y(s), v), \partial \Dom(\L))<2\rho_\ell\Rightarrow \L(s,y(s),v)\ge\ell.\end{equation}
For $0<\rho\le\rho_{\ell}$, we thus have
\[
B\ge \int_t^T\L(s,y(s),u(s))\,ds\ge
\int_{[t,T]\setminus J_{\rho}}\L(s,y(s),u(s))\,ds\ge \ell|[t,T]\setminus J_{\rho}|,
\]
whence
$B\ge \ell\big(T-t-|J_{\rho}|\big)$, from which we obtain the estimate
\begin{equation}\label{tag:qgiqg}\forall\, 0<\rho\le\rho_{\ell}\qquad |J_{\rho}|\ge \dfrac{\ell(T-t)-B}{\ell}=(T-t)-\dfrac{B}{\ell}.\end{equation}
The claim follows.
\item[{\em vi)}]{\em Let $\overline\rho>0$ be as in Claim iii). There are $\mu=\mu(B, \delta)\in ]0, 1[, \rho\le \overline\rho$, $m:=m(\delta, B)\in ]0, 1]$ and a subset $\Omega$ of $[t, T]$ with $|\Omega|\ge m(T-t)$, such that, for a.e. $s\in\Omega$ and a.e. $\tilde s\in [t, T]$,
    \begin{equation}\label{tag:NEW!}\L\left(\tilde s, y(s), \dfrac{u(s)}{\mu}\right)<+\infty, \,\dfrac{|u(s)|}{\mu}<c,\, \dist \left(\left(\tilde s,y(s), \dfrac{|u(s)|}{\mu}\right), \partial \Dom(\L)\right)\ge\rho.\end{equation}}\\
Choose $\mu_0, \Delta\in ]0,1[$ such that
\begin{equation}\label{tag:mu0Delta} \dfrac{c_{\delta}(B)}{\mu_0 c}<\Delta<1,\end{equation}
and let, from Step {\em v)}, $\rho\in ]0, \overline\rho]$ be such that \begin{equation}\label{tag:chooserho}|J_{\rho}|\ge \Delta (T-t).\end{equation}
Set
\begin{equation}\label{tag:defmu}\mu:=\dfrac{c}{\rho+c}\,\in ]0,1[, \quad \Omega:=\Omega_{\mu}\cap J_{\rho}.\end{equation}
Since $\displaystyle\lim_{\rho\to 0}\dfrac{c}{\rho+c}=1$, we may assume that $\mu\ge \mu_0$.
Since, from Step {\em iv)}, for a.e. $s$ on $\Omega_{\mu}$,   $|u(s)|<  \mu c<c$, we deduce that, for a.e. $s\in \Omega$,
$\left|\dfrac{u(s)}{\mu}-u(s)\right|=\rho\dfrac{|u(s)|}{c}<\rho$,
so that
$\dist\left(\left(s, y(s), \dfrac{u(s)}{\mu}\right), (s, y(s), u(s))\right)\le \rho$ and thus, for a.e. $s\in\Omega$,
\begin{equation}\label{tag:omega33}\L\left(s, y(s),
 \dfrac{u(s)}{\mu}
  \right)<+\infty\,\quad
 \dist\left(\left(s, y(s), \dfrac{u(s)}{\mu}
  \right), \partial\Dom(\L)\right)\ge\rho\,\text{ a.e. }s\in \Omega.\end{equation}
  Hypothesis h$_1$) implies the validity of \eqref{tag:NEW!} for a.e.  $\tilde s\in [t,T]$ and a.e. $s\in \Omega$.
It remains to show that $|\Omega|$ is big enough.  It follows from Step {\em iv)} and \eqref{tag:chooserho} that
\begin{equation}\label{tag:gqigfqe}\begin{aligned}
|\Omega|&=|\Omega_{\mu}|+
|J_{\rho}|-|\Omega\cup J_{\rho}|\\&\ge
 \left(1-\dfrac{c_{\delta}(B)}{\mu c}\right)(T-t)+\Delta (T-t)-(T-t)\\
 &\ge \left(\Delta-\dfrac{c_{\delta}(B)}{\mu c}\right)(T-t)\ge m(T-t),\end{aligned}\end{equation}
  with $m:=\Delta-\dfrac{c_{\delta}(B)}{\mu} c$.
\item[{\em vii)}]
\emph {From  now on we set
$ \Upsilon:=\Upsilon(\rho)$.
For a.e.   $s\in \Omega$ and a.e.  $\tilde s \in [t,T]$,
\begin{equation}\label{tag:C} \L\left(\tilde s , y(s), \dfrac{u(s)}{\mu}\right)\mu-\L(\tilde s , y(s), u(s))\le -(1-\mu) \Upsilon .\end{equation}}
Indeed, it follows from Step {\em vi)} and the basic assumptions on $\Dom(\L)$  that, for a.e. $s\in \Omega$, for a.e.  $\tilde s \in [t, T]$,
\begin{equation}\label{tag:star}\forall 0<r\le \dfrac1{\mu}\qquad
\big(\tilde s , y(s), ru(s)\big)\in\Dom(\L).
\end{equation}
Notice that, for such $s,\tilde s $,
\begin{align}\label{tag:funda}
\L\left(\tilde s , y(s), \dfrac{u(s)}{\mu}\right)\mu-\L(\tilde s , y(s), u(s))=\phantom{AAAAAAAAAAAAAAAAA}\\
\\
\phantom{AAAAAAAAAAA}=-\mu
\left[\L\left(\tilde s , y(s), \dfrac{u(s)}{\mu}\mu\right)\dfrac1{\mu}-\L\left(\tilde s , y(s), \dfrac{u(s)}{\mu}\right)\right].
\end{align}
\noindent
We consider separately the radial convex and the partial differentiable cases.
\begin{itemize}[leftmargin=*]
\item[a)] {\em  Radial convex case.} We have
\[
\L\left(\tilde s , y(s), \dfrac{u(s)}{\mu}\mu\right)\dfrac1{\mu}-\L\left(\tilde s , y(s), \dfrac{u(s)}{\mu}\right)\ge
P\left(\tilde s , y(s), \dfrac{u(s)}{\mu}\right)\dfrac{1-\mu}{\mu}.\]
From Step {\em vi)} and the fact that $\mathcal U$ is a cone we obtain
\begin{align}\label{tag:dist*}\L\left(\tilde s , y(s), \dfrac{u(s)}{\mu}\mu\right)\dfrac1{\mu}-\L\left(\tilde s , y(s), \dfrac{u(s)}{\mu}\right)\ge\phantom{AAAAAAAAAAAAA}\\
\ge\dfrac{1-\mu}{\mu}\!\!\!\!\inf_{\substack{\tilde s \in [0,T],|z|\le K\\|v|<c, v\in \mathcal U\\ \L(\tilde s , z, v)<+\infty\\ \dist((\tilde s ,z,v),\partial\Dom(\L))\ge\rho}}
P(\tilde s ,z,v)=\dfrac{1-\mu}{\mu}\Upsilon ,\end{align}
and thus the conclusion follows from
(\ref{tag:funda}).
\item[b)]
{\em  Partial differentiable case.}
For $\omega>0$ let $h(\omega):=\L\left(\tilde s , y(s), \dfrac{u(s)/\mu}{\omega}\right)\omega$.
Then
\[\begin{aligned}\label{tag:finaldiff}\L&\left(\tilde s , y(s), \dfrac{u(s)}{\mu}\mu\right)\dfrac1{\mu}-\L\left(\tilde s , y(s), \dfrac{u(s)}{\mu}\right)=\\
&\phantom{AAAAAAAAAAAAAAAAA}=h\left(\dfrac1{\mu}\right)-h(1)=
\int_1^{\frac1{\mu}}h'(\omega)\,d\omega.
\end{aligned}\]
It follows from \eqref{tag:changevar2} that, for all $\omega\in \left[1, \dfrac1{\mu}\right]$,
\[h'(\omega)=\dfrac{d}{d\lambda}
\left[\L\left(\tilde s , y(s), \dfrac{u(s)/(\omega\mu)}{\lambda}\right)\lambda\right]_{\lambda=1}.\]
We deduce from  that, for a.e. $s\in\Omega$,
\[\forall \omega\in\left[1, \dfrac1{\mu}\right]\qquad  \left|\dfrac{u(s)}{\omega\mu}\right|\le \left|\dfrac{u(s)}{\mu}\right|<c.\]
It follows  from  Step {\em vi)}  and the fact that $\mathcal U$ is a cone that
\[\begin{aligned}\label{tag:dist**}\L\left(\tilde s , y(s), \dfrac{u(s)}{\mu}\mu\right)&\dfrac1{\mu}-\L\left(\tilde s , y(s), \dfrac{u(s)}{\mu}\right)\ge\phantom{AAAAAAAAAAA}\\
&\ge\dfrac{1-\mu}{\mu}
\inf_{\substack{\tilde s\in [0,T],|z|\le K\\|v|<c, v\in \mathcal U\\ \L(\tilde s, z, v)<+\infty\\ \dist((\tilde s,z,v),\partial\Dom(\L))\ge\rho}}
\dfrac d{d\lambda}\left[\L\left(\tilde s, z, \dfrac{v}{\lambda}\right)\lambda\right]_{\lambda=1}\\
&=\dfrac{1-\mu}{\mu}\Upsilon,  \,\,\end{aligned}\]
where the last equality is a consequence of Remark~\ref{rem:variediff}.
The claim follows from (\ref{tag:funda}).
\end{itemize}
\item[{\em viii)}]
\emph {For every $\nu>0$ define \[S_{\nu}:=\{s\in [t,T]:\, |u(s)|> \nu\},\quad \varepsilon_{\nu}:=
\displaystyle\int_{S_{\nu}}\left(\dfrac{|u(s)|}{\nu}-1\right)\,ds.\]
Then \[|S_{\nu}|\to 0,\quad \varepsilon_{\nu}\le \dfrac{R}{\nu}\rightarrow 0\text{ as }\nu\to +\infty\]
uniformly with respect to $t\in [0, \delta]$ and $x\in B_{\delta_*}^n(x_*)$.
}
Indeed, it follows from Step {\em i)} that
\[\nu |S_{\nu}|\le\int_{S_{\nu}}|u(s)|\,ds\le R.\]
\item[{\em ix)}]
{\em Choice of $\nu\ge \overline\nu$ and of $\Sigma_{\nu}\subseteq \Omega$.}\\
Taking into account Claim {\em viii)}, we choose $\nu=\nu(\delta, B, \delta_*, x_*)\ge \max\{\overline\nu, c\}$  in such a way that \begin{equation}\label{tag:choicenu33}\dfrac{R}{\nu}\le \min\left\{\,(1-\mu)m(T-\delta), \dfrac{\varepsilon_*}2\right\}.\end{equation}
If Condition (M$_B^{\delta}$) holds, we impose moreover that $\nu=\nu(\delta, B, \delta_*, x_*, \eta)$ is large enough in such a way that
    \begin{equation}\label{tag:nularge}\dfrac R{\nu}(2\Phi(B)+\Xi(\nu)-\Upsilon)\le \dfrac R{\nu}(2\Phi(B)+\Xi(\overline\nu)-\Upsilon)\le\eta. \end{equation}
  From now on we set $\Xi:=\Xi(\nu)$.
Choose a measurable subset ${\Sigma_{\nu}}$ of $ \Omega$ in such a way that $|{\Sigma_{\nu}}|=\dfrac{\varepsilon_{\nu}}{1-\mu}$: this is possible since, from  \eqref{tag:choicenu33} and Step {\em vi)}, \[\dfrac{\varepsilon_{\nu}}{1-\mu}\le m(T-\delta)\le  m(T-t)\le |\Omega|.\]
\item[{\em x)}]
{\em  $S_{\nu}\cap \Omega$ is negligible.}\\
Indeed,  if $u(s)$ is defined and $s\in S_{\nu}$ then $|u(s)|> {\nu}$, whereas if $s\in \Omega$ then $|u(s)|< c\le\nu$.
\item[{\em xi)}]
{\em The change of variable $\varphi$.}
We  introduce the following absolutely continuous change of variable $\varphi:[t,T]\to \R$ defined by
\[\varphi(t):=t,\quad\text{for a.e. }\tau\in [t,T]\quad \varphi'(\tau)=\begin{cases}\dfrac{|u(\tau)|}{{\nu}} &\text{ if }\tau \in S_{\nu},\\
\,\,\mu &\text{ if }\tau\in{\Sigma_{\nu}},\\
\,\,1 &\text{ otherwise}.
\end{cases}\]
Notice that $\varphi$ depends on both $y,u$ and is well defined since $S_{\nu}\cap {\Sigma_{\nu}}$, a subset of $S_{\nu}\cap \Omega$, is  negligible.
Clearly $\varphi$ is  strictly increasing and, from steps {\em viii) -- ix)},
\[\begin{aligned}\int_{t}^T\varphi'(\tau)\,d\tau&=
\int_{S_{\nu}}\dfrac{|u(\tau)|}{{\nu}}\,d\tau+
\int_{{\Sigma_{\nu}}}\mu\,d\tau+|[t,T]\setminus S_{\nu}\cup {\Sigma_{\nu}}|\\
&=\left({\varepsilon_{\nu}}+|S_{\nu}|\right)+\mu|{\Sigma_{\nu}}|+\big((T-t)-|S_{\nu}|-|{\Sigma_{\nu}}|\big)\\
&={\varepsilon_{\nu}}-(1-\mu)|{\Sigma_{\nu}}|+(T-t)=T-t.\end{aligned}\]
Therefore the image of $\varphi$ is $[t,T]$ and thus $\varphi:[t,T]\to [t,T]$ is bijective; let us denote by $\psi$ its inverse, which is  absolutely continuous and even Lipschitz, since $\|\psi'\|_{\infty}\le\dfrac1{\mu}$.
\item[{\em xii)}]
{\em
Set
$\overline u(s):=\dfrac{u(\psi(s))}{\varphi'(\psi(s))}$, $\overline y:=y\circ\psi$. Then $(\overline y, \overline u)$ is admissible and  $\overline y(T)=y(T)$.}\\
 It follows from \cite[Corollary 5]{Serrin} and \cite[Chapter IX, Theorem 5]{Natanson} that $\overline y$ is absolutely continuous and that,
 for a.e. $s\in [t, T]$,
\[
\overline y'(s)=y'(\psi(s))\dfrac1{\varphi'(\psi(s))}=b(y(\psi(s))\dfrac{u(\psi(s))}{\varphi'(\psi(s))}=b(\overline y(s))\overline u(s).
\]
Since $\overline y$ is defined via a reparametrization of $y$, we still have that $\overline y(s)\in\mathcal S$ for all $s$. Moreover, $\overline y(t)=y(\psi(t))=y(t)$ and $\overline y(T)=y(\psi(T))=y(T)$.
Notice that \begin{equation}\label{tag:cone!}\overline u(s)=\dfrac{1}{\varphi'(\psi(s))}u(\psi(s))\in \mathcal U \text{ a.e. } s\in [t,T],\end{equation} the set $\mathcal U$ being a cone.
\item[{\em xiii)}]{\em $\overline u$ is bounded, $\overline y$ is Lipschitz;  if Condition} {\rm(}H$_B^{\delta}${\rm)} {\em holds the bound of $\overline u$ and the Lipschitz rank of $\overline y$  depend just on $\delta, B, \delta_*, x_*$; otherwise they might depend also on $\eta$.}\\
It is convenient to write explicitly the function $\overline u(s)$, which is given by
\[\overline u(s)=\begin{cases} \nu\dfrac{\,u(\psi(s))}{|u(\psi(s))|} &\text{ if }\psi(s) \in S_{\nu},\\
\dfrac{u(\psi(s))}{\mu} &\text{ if }\psi(s)\in{\Sigma_{\nu}},\\
u(\psi(s)) &\text{ otherwise.}
\end{cases}\]
Since $|u(s)|\le  \nu$ a.e. out of $S_{\nu}$ it turns out from the fact that $\Sigma_{\nu}\subseteq\Omega$ that
 \begin{equation}\label{tag:estimnu}|\overline u(s)|\le
 \max\left\{\nu, c\right\}=\nu.\end{equation}
%
%
Now
\begin{equation}\label{tag:estimLip}\|\overline y'\|_{\infty}=\|b(\overline y)\overline u\|_{\infty}\le \theta(1+\|\overline y\|_{\infty})\|\overline u\|_{\infty}=\theta(1+\| y\|_{\infty})\|\overline u\|_{\infty}\le
\theta(1+K)\nu;\end{equation}
 the claim follows from Step {\em ii)}.
%
\item[{\em xiv)}]
{\em $\|\varphi(\tau)-\tau\|_{\infty}\le 2\varepsilon_{\nu}\le \varepsilon_*$.}
Indeed, for all $\tau\in [t,T]$ we have
\[\begin{aligned}|\varphi(\tau)-\tau|&\le \int_{t}^{\tau}\left|\varphi'(s)-1\right|\,ds\\
&\le \int_{S_{\nu}}\left(\dfrac{|u(s)|}{\nu}-1\right)\,ds+
\int_{{\Sigma_{\nu}}}(1-\mu)\,ds\\&\le{\varepsilon_{\nu}}
+(1-\mu) |{\Sigma_{\nu}}|=2{\varepsilon_{\nu}}\le\varepsilon_*.
\end{aligned}\]
\item[{\em xv)}]
{\em Estimate of $J_{t}(\overline y, \overline u)$ in terms of $\displaystyle\int_{t}^T\L\big(\varphi(\tau),  y(\tau),u(\tau)\big)\,d\tau$.}\\
 Since $y$ and $\overline y$ share the same boundary values,
we have
\begin{equation}\label{tag:Jm}J_{t}(\overline y, \overline u)=
\int_{t}^T\L(s, \overline y(s), \overline u(s))\,ds+g(y(T)).\end{equation}
The change of variables $s=\varphi(\tau)$ yields
\begin{equation}\begin{aligned}\label{tag:valuenewJ}\int_{t}^T\L(s, \overline y(s), \overline u(s))\,ds&=
\int_{t}^T\L\Big(\varphi(\tau),  y(\tau), \dfrac{u(\tau))}{\varphi'(\tau)}\Big)\varphi'(\tau)\,d\tau\\
&=I_{S_{\nu}}+I_{{\Sigma_{\nu}}}+I_1,\end{aligned}\end{equation}
where we set
\[\begin{aligned}I_{S_{\nu}}&:=
\int_{S_{\nu}}\L\Big(\varphi(\tau),  y(\tau), \nu\dfrac{u(\tau)}{|u(\tau)|}\Big)\dfrac{|u(\tau)|}{\nu}\,d\tau,\\
I_{{\Sigma_{\nu}}}&:=
\int_{{\Sigma_{\nu}}}\L\Big(\varphi(\tau),  y(\tau), \dfrac{u(\tau)}{\mu}\Big)\mu\,d\tau,\\
 I_{1}&:=
\int_{[t,T]\setminus ({\Sigma_{\nu}}\cup S_{\nu})} \L\left(\varphi(\tau),  y(\tau), u(\tau)\right)\,d\tau.\end{aligned}\]
In what follows, for brevity, we set $\displaystyle\int_*\L(\varphi, y, u)\,d\tau :=\displaystyle\int_*\L(\varphi(\tau), y(\tau), u(\tau))\,d\tau$.
\begin{itemize} [leftmargin=*]
\item[$\bullet$] {\em Estimate of $I_{S_{\nu}}$:}
\begin{equation}\label{tag:IS}I_{S_{\nu}}\le \int_{S_{\nu}}\L\big(\varphi,  y, u\big)\,d\tau+\Xi\varepsilon_{\nu}.\end{equation}
Indeed, recall that $(s, y(s), u(s))\in \Dom(\L)$ for a.e.  $s\in [t,T]$.
Since $|u(\tau)|>\nu$  a.e.  $\tau\in S_{\nu}$ and $\mathcal U$ is a cone then, for a.e. $\tau\in S_{\nu}$,
\[\quad \L\Big(\varphi(\tau),  y(\tau), \nu\dfrac{u(\tau)}{|u(\tau)|}\Big)<+\infty, \quad \nu\dfrac{u(\tau)}{|u(\tau)|}\in\mathcal U.\]
We consider the two following frameworks.
\begin{itemize}[leftmargin=*]
\item[a)]  {\em  Radial convex case.} For a.e. $\tau\in S_{\nu}$ we have
\begin{align}
\L\Big(\varphi(\tau),  y(\tau),u(\tau)\Big)\dfrac{\nu}{|u(\tau)|}-\L\Big(\varphi(\tau),  y(\tau), \nu\dfrac{u(\tau)}{|u(\tau)|}\Big)\ge \phantom{AAAAAAAAAAAAA}\\ \ge
P\Big(\varphi(\tau), y(\tau), \nu\dfrac{u(\tau)}{|u(\tau)|}\Big)\Big(\dfrac{\nu}{|u(\tau)|}
-1\Big)
\end{align} which implies
\begin{equation}\label{tag:gyegqyg}\begin{aligned}
\L\Big(\varphi(\tau),  y(\tau), \nu\dfrac{u(\tau)}{|u(\tau)|}\Big)&\dfrac{|u(\tau)|}{\nu}-
\L\Big(\varphi(\tau),  y(\tau),u(\tau)\Big)\le \\\phantom{AAAAAAAAA}&\phantom{AAAAA} \le
P\Big(\varphi(\tau), y(\tau), \nu\dfrac{u(\tau)}{|u(\tau)|}\Big)\Big(\dfrac{|u(\tau)|}{\nu}-1
\Big).
\end{aligned}\end{equation}
Since, for a.e. $\tau\in S_{\nu}$, $|u(\tau)|>\nu$ and
 $\left|\nu\dfrac{u(\tau)}{|u(\tau)|}\right|=\nu$,
 we deduce  that
 \[\text{a.e. }\tau\in S_{\nu}\quad P\Big(\varphi(\tau), y(\tau), \nu\dfrac{u(\tau)}{|u(\tau)|}\Big)\Big(\dfrac{|u(\tau)|}{\nu}-1
\Big)\le\Big(\dfrac{|u(\tau)|}{\nu}-1
\Big)\Xi.\]
Therefore,   for a.e. $\tau\in S_{\nu}$ inequality \eqref{tag:gyegqyg} yields
\[\L\Big(\varphi(\tau),  y(\tau), \nu\dfrac{u(\tau)}{|u(\tau)|}\Big)\dfrac{|u(\tau)|}{\nu}\le
\L\big(\varphi(\tau),  y(\tau),u(\tau)\big)+\Big(\dfrac{|u(\tau)|}{\nu}-1
\Big)\Xi,\]
whence \eqref{tag:IS}.
\item[b)]
{\em Partial differentiable case.}
For a.e. $\tau\in S_{\nu}$ we set
\begin{equation}\label{tag:tau}h_{\tau}(\mu):=\L\Big(\varphi(\tau),  y(\tau), \dfrac{u(\tau)}{\mu}\Big)\mu,
\quad \mu\in \left[1, \dfrac {|u(\tau)|}{\nu}\right].
\end{equation}
Fix  such a value of $\tau$. Since
\[\L\Big(\varphi(\tau),  y(\tau), \nu\dfrac{u(\tau)}{|u(\tau)|}\Big)\dfrac{|u(\tau)|}{\nu}-
\L\Big(\varphi(\tau),  y(\tau),u(\tau)\Big)=h_{\tau}\left(\dfrac{|u(\tau)|}{\nu}\right)-h_{\tau}(1),\]
there is $\mu_{\tau}\in \left]1, \dfrac{|u(\tau)|}{\nu}\right[$ satisfying
\[\L\Big(\varphi(\tau),  y(\tau), \nu\dfrac{u(\tau)}{|u(\tau)|}\Big)\dfrac{|u(\tau)|}{\nu}-
\L\Big(\varphi(\tau),  y(\tau),u(\tau)\Big)=h_{\tau}'(\mu_{\tau})\left(\dfrac{|u(\tau)|}{\nu}-1 \right).\]
Now,  from  (3) of Remark~\ref{rem:variediff} we have
\[h_{\tau}'(\mu_{\tau})=\dfrac{d}{d\lambda}\left[\L\Big(\varphi(\tau),  y(\tau), \dfrac{u(\tau)/\mu_{\tau}}{\lambda}\Big)\lambda\right]_{\lambda=1}.\]
Since $\dfrac{|u(\tau)|}{\mu_{\tau}}>\nu$ on $S_{\nu}$, we deduce that $|h_{\tau}'(\mu_{\tau})|\le \Xi$. It follows that, on $S_{\nu}$,
\[\L\Big(\varphi(\tau),  y(\tau), \nu\dfrac{u(\tau)}{|u(\tau)|}\Big)\dfrac{|u(\tau)|}{\nu}-
\L\Big(\varphi(\tau),  y(\tau),u(\tau)\Big)\le \Xi\left(\dfrac{|u(\tau)|}{\nu}-1 \right),\]
from which  we obtain \eqref{tag:IS}.
\end{itemize}
\item[$\bullet$] {\em Estimate of $I_{{\Sigma_{\nu}}}$.}
The function $\psi$ being Lipschitz, the set of $\tau\in [t,T]$ such that $\tilde s =\varphi(\tau)$ satisfies  \eqref{tag:C} is of full measure.
Since ${\Sigma_{\nu}}\subseteq \Omega$ and $|{\Sigma_{\nu}}|=\dfrac{\varepsilon_{\nu}}{1-\mu}$,  it is immediate  from Step {\em vii)}  that
\begin{equation}\label{tag:ISigma}I_{{\Sigma_{\nu}}}\le \int_{{\Sigma_{\nu}}}\L\big(\varphi,  y, u\big)\,d\tau-(1-\mu)\Upsilon |{\Sigma_{\nu}}|
\le \int_{{\Sigma_{\nu}}}\L\big(\varphi,  y, u\big)\,d\tau-\Upsilon {\varepsilon_{\nu}}.
\end{equation}
\end{itemize}
Therefore, from \eqref{tag:Jm}, \eqref{tag:valuenewJ}, \eqref{tag:IS} and \eqref{tag:ISigma} we deduce the required estimate
\begin{equation}\label{tag:uqfuqyfdque}\begin{aligned}J_{t}(\overline y, \overline u)&
=I_{S_{\nu}}+I_{{\Sigma_{\nu}}}+I_1+g(y(T))\\
&\le
\int_{t}^T\L\big(\varphi,  y, u\big)\,d\tau+
\varepsilon_{\nu}\Big(\Xi-\Upsilon \Big)+g( y(T)).
\end{aligned}\end{equation}
\item[{\em xvi)}]
{\em Estimate of $\displaystyle\int_{t}^T\L\big(\varphi,  y, u\big)\,d\tau$:}
\begin{equation}\label{tag:qigdiuqgdf}
\int_{t}^T\L\big(\varphi,  y, u\big)\,d\tau\le
\int_t^T\L\big(\tau,  y(\tau),u(\tau)\big)\,d\tau+2\Phi(B)\varepsilon_{\nu}.
\end{equation}
Indeed, the assumptions of Theorem~\ref{thm:main2}  ensure that $(\varphi(\tau), y(\tau), u(\tau))\in\Dom(\Lambda)$ for a.e. $\tau\in [t,T]$.
Condition {\rm (S)} and Step {\em xiv)} then imply that, for a.e. $\tau\in [t,T]$,
\[\begin{aligned}\L\big(\varphi(\tau),  y(\tau),u(\tau)\big)
&\le \L\big(\tau,  y(\tau),u(\tau)\big)+k(\tau)|\varphi(\tau)-\tau|\\
&\le \L\big(\tau,  y(\tau),u(\tau)\big)+2k(\tau){\varepsilon_{\nu}},\end{aligned}\]
where we set
\[k(\tau):=\kappa\L(\tau,  y(\tau),u(\tau))+A|u(\tau)|+\gamma(\tau).\]
Notice that, from Proposition~\ref{lemma:newC},
\[\int_t^Tk(\tau)\,d\tau\le \Phi(B):\]
\eqref{tag:qigdiuqgdf} follows  now immediately.
\item[{\em xvii)}]
{\em Final estimate of $J_{t}(\overline y, \overline u)$.}
From \eqref{tag:uqfuqyfdque} and \eqref{tag:qigdiuqgdf} of Steps {\em xv) -- xvi)}, we obtain
\begin{equation}\label{tag:finalX}J_{t}(\overline y, \overline u)\le J_{t}( y, u)+\varepsilon_{\nu}\left(2\Phi(B)+\Xi-\Upsilon \right).
\end{equation}
Two cases may occur.
\begin{itemize}
\item
If Condition (H$_B^{\delta}$) holds true, the choice of $\nu$ in \eqref{tag:choicenu50} of Step {\em iii)} implies that \[2\Phi(B)+\Xi-\Upsilon <0.\]  Thus, from \eqref{tag:finalX}  we  obtain $J_{t}(\overline y, \overline u)\le J_{t}( y, u)$. Actually, the inequality $J_{t}(\overline y, \overline u)< J_{t}( y, u)$ is strict, unless $\varepsilon_{\nu}=0$, and this occurs if and only if $|u|\le \nu$ a.e. on $[t,T]$.
\item If Condition (M$_B^{\delta}$) holds true, from \eqref{tag:finalX} and  \eqref{tag:nularge} of Step {\em ix)} we obtain
    \[J_{t}(\overline y, \overline u)\le J_{t}( y, u)+\dfrac R{\nu}(2\Phi(B)+\Xi-\Upsilon).\]
    The conclusion follows from the choice of $\nu$ in \eqref{tag:nularge}.
\end{itemize}
\end{itemize}
\end{proof}
\begin{remark}[\textbf{Explicit bounds and Lipschitz ranks}] In the {\em real valued case}, the knowledge of $\overline\nu$ and $c$ in  Condition (H$_B^{\delta}$) correspondingly to the value $K$ provided in Step {\em ii)} of the proof of Theorem~\ref{thm:main2} allows to give an explicit bound  of $\|\overline u\|_{\infty}$ and $\|\overline y'\|_{\infty}$) (thus of $K_{\mathcal A}$ in Claim (1) of Theorem~\ref{thm:main2}). Indeed, referring to the proof of the theorem:
\begin{itemize}
\item From Step {\em xiii)},  $\|\overline u\|_{\infty}\le\nu$ and $\|\overline y'\|_{\infty}\le
\theta(1+K)\nu$  where, from Step {\em ix)}, one may take
    \begin{equation}\label{tag:nununu}\nu:=\max\left\{\dfrac{R}{(1-\mu)m(T-\delta)}, c, \overline\nu, \dfrac{2R}{\varepsilon_*}\right\}\end{equation}
and, from Step {\em i)}, $R=\dfrac{B+d\,T}{\alpha}$;
    \item In the real valued case it is enough from Step {\em vi)} to take $\mu$, equal to any real number such that $\dfrac{c_{\delta}(B)}{c}<\mu<1$;
    \item The proof of Step {\em vi)} shows that $m$ in \eqref{tag:nununu} is given by
    \[m=\Delta-\dfrac{c_{\delta}(B)}{\mu}c,\]
    where $\Delta$ is any real number in $\left]\dfrac{c_{\delta}(B)}{\mu} c,1\right[$.
\end{itemize}
\end{remark}
\begin{remark}\label{rem:dist} The proof of Theorem~\ref{thm:main2} shows that one could replace in the growth Conditions (H$_B^{\delta}$), (M$_B^{\delta}$) and in Hypothesis h$_2$) the Euclidean distance $\dist((s,y,v), \partial \Dom(\L))$ with the pseudo-distance
\[\dist_c((s,y,v), \partial \Dom(\L))):=\inf\{|v'-v|:\, (s,y,v')\in\partial\Dom(\L)\}:\]
indeed in this case $\dist_c((s,y,v), \partial \Dom(\L))\ge\rho>0$ whenever $(s,y,v')\in\Dom(\L)$  for every $v'\in B_{\rho}^m(v)$, so that, from Hypothesis h$_1$) it follows  that  \[\dist_c((s,y,v), \partial \Dom(\L))\ge\rho\Rightarrow \dist_c((\tilde s,y,v), \partial \Dom(\L))\ge\rho\quad \forall\tilde s\in [0, T],\]  an essential property in Step {\em vi)} of the proof of Theorem~\ref{thm:main2}.
There are some advantages and drawbacks in replacing the Euclidean distance $\dist$ with the greater $\dist_c$.  Indeed, in doing so, the infima in the growth conditions become smaller, so that (H$_B^{\delta}$), (M$_B^{\delta}$) become more restrictive (i.e., are satisfied by a smaller class of Lagrangians).  Instead, Hypothesis h$_{2}$)  is less restrictive (i.e., satisfied by a wider class of Lagrangians). However,  the notion of being well-inside the domain for $\dist_c$ (Definition~\ref{def:wellinside}) does not correspond anymore to the notion of  relatively compact subset.
\end{remark}

\bibliographystyle{amsplain}

\bibliography{biblio_Slow_Growth_Trans}


\end{document}